\documentclass[twoside]{ucthesis}

\usepackage{amsmath}
\usepackage{amsthm}
\usepackage{amsfonts}
\usepackage{latexsym}
\usepackage[T1]{fontenc}

\makeatletter

\newcommand{\LyX}{L\kern-.1667em\lower.25em\hbox{Y}\kern-.125emX\@}

\theoremstyle{plain}    
\newtheorem{thm}{Theorem}[section]
\theoremstyle{plain}    
\newtheorem{cor}[thm]{Corollary} 
\theoremstyle{plain}    
\newtheorem{lem}[thm]{Lemma} 
\theoremstyle{plain}    
\newtheorem{prop}[thm]{Proposition} 
\theoremstyle{plain}    
\newtheorem{fact}[thm]{Fact} 
\theoremstyle{definition}
\newtheorem{defn}[thm]{Definition}
\theoremstyle{definition}
 \newtheorem{example}[thm]{Example}
\theoremstyle{definition}
 \newtheorem*{example*}{Example}
\theoremstyle{remark}
\newtheorem{rem}[thm]{Remark}
\theoremstyle{remark}
\newtheorem*{rem*}{Remark}
\theoremstyle{remark}    
\newtheorem*{note*}{Note} 
\theoremstyle{remark}    
\newtheorem{notation}[thm]{Notation} 
\theoremstyle{remark}    
\newtheorem*{acknowledgement*}{Acknowledgement} 
\theoremstyle{remark}    
\newtheorem{case}{Case} 

\usepackage[T1]{fontenc}

\makeatletter

\usepackage[T1]{fontenc}

\makeatletter

\usepackage[T1]{fontenc}

\makeatletter

\newenvironment{pf*}[1]{\noindent\mbox{{\em 
{#1}}.}}{\\\hspace*{\fill}$\Box$\\\medskip}

\newcommand{\br}{[\cdot,\cdot ]}
\newcommand{\bform}[2]{\langle #1 , #2 \rangle}
\newcommand{\cinf}{C^{\infty}}
\newcommand{\lon }{\longrightarrow }
\newcommand{\half}{\frac{1}{2}}

\newcommand{\fourth}{\frac{1}{4}}
\newcommand{\sixth}{\frac{1}{6}}
\newcommand{\smalcirc}{\mbox{\tiny{$\circ $}}}
\newcommand{\DD}{{\mathcal D}}
\newcommand{\coker}{\mbox{coker}\;}
\newcommand{\w}{\wedge}
\newcommand{\bw}{\bigwedge}
\newcommand{\g}{\mathfrak{g}}
\newcommand{\X}{\mathfrak{X}}
\newcommand{\R}{\mathbb{R}} 
\newcommand{\A}{{\mathcal A}}
\newcommand{\B}{{\mathcal B}}
\newcommand{\p}{\mathfrak{p}}
\newcommand{\Z}{\mathbb{Z}}
\newcommand{\OO}{{\mathcal O}}
\makeatother

\makeatother
\makeatother
\makeatother

\begin{document}

\title{Courant algebroids, derived brackets and even symplectic supermanifolds}

\author{Dmitry Roytenberg}

\degreeyear{1999}
\degree{Doctor of Philosophy}
\chair{Professor Alan D. Weinstein}
\othermembers{Professor Alexander B. Givental\\
Professor Robert G. Littlejohn}
\numberofmembers{3}
\prevdegrees{B.A. (New York University) 1993}
\field{Mathematics}
\campus{Berkeley}

\maketitle


\begin{abstract}


In this dissertation we study Courant algebroids, objects that first appeared in 
the work of T. Courant on Dirac structures; they were later studied by Liu, 
Weinstein and Xu who used Courant algebroids to generalize the notion of the 
Drinfeld double to Lie bialgebroids. As a first step towards understanding the 
complicated properties of Courant algebroids, we interpret them by associating 
to each Courant algebroid a strongly homotopy Lie algebra in a natural way.

Next, we  propose an alternative construction of the double of a Lie bialgebroid as a 
homological hamiltonian vector field on an even symplectic supermanifold. The 
classical BRST complex and the Weil algebra arise as special cases. We recover 
the Courant algebroid via the derived bracket construction and give a simple 
proof of the doubling theorem of Liu, Weinstein and Xu. We also introduce a 
generalization, quasi-Lie bialgebroids, analogous to Drinfeld's quasi-Lie 
bialgebras; we show that the derived bracket construction in this case also 
yields a Courant algebroid.

Finally, we compute the Poisson cohomology of a one-parameter family of $SU(2)$-covariant 
Poisson structures on $S^2$. As an application, we show that these structures 
are non-trivial deformations of each other, and that they do not admit 
rescaling.

\abstractsignature
\end{abstract}

\begin{frontmatter}

\begin{dedication}
\null\vfil
{\large
\begin{center}
To the memory of Nikolai Afanasievich Pravdin.
\end{center}}
\vfil\null
\end{dedication}

\tableofcontents

\begin{acknowledgements}
First and foremost, I would like to thank my advisor, Professor Alan Weinstein,
for his continual guidance, encouragement and support throughout my graduate
studies. At times he was more patient with me than I deserved. It was from him
that I got my first inspiration to study geometry, and I am still happy with
my choice.

I would also like to thank Professor Alexander Givental with whom I have spent
many hours discussing mathematics. His comments were always lucid, intelligent
and to the point. I have enjoyed our conversations and learned a lot from them.

I am particularly grateful to Professor Theodore Voronov, a close friend and
a valuable colleague, who taught me the theory of supermanifolds during his
visit in Berkeley. My gratitude also goes to Professor Yvette Kosmann-Schwarzbach,
many of whose ideas provided impetus for this work, and to Professor James Stasheff
for helpful discussions and comments on the manuscript.

I am deeply thankful to Marika Zavodovskaya, my fianc\'{e}; her support and
understanding throughout this project has been a great help. I am also indebted
to my parents who were behind me every step of the way.

Last but not least, I would like to thank all my good friends who have kept
me company during my stay in Berkeley, for all the great times we've had together. 
\end{acknowledgements}

\end{frontmatter}

\chapter{Introduction}

The first example of a Courant algebroid appeared in the work of T. Courant
\cite{Cou} on Dirac structures. These structures are a simultaneous generalization
of pre-symplectic and Poisson structures; they appear in Dirac's theory of constrained
mechanical systems. A \emph{Dirac structure} on a manifold \( M \) is a subbundle
\( L\subset TM\oplus T^{*}M \) that is maximally isotropic with respect to
the canonical symmetric bilinear form on \( TM\oplus T^{*}M \), and which satisfies
a certain integrability condition. To formulate the integrability condition,
Courant introduced a bilinear skew-symmetric bracket operation 
\[
[X+\xi ,Y+\eta ]=[X,Y]+(L_{X}\eta -L_{Y}\xi +\half d(i_{Y}\xi -i_{X}\eta ))\]
 on sections of \( TM\oplus T^{*}M \); the condition is that the sections of
\( L \) be closed under this bracket. As one can see, the Courant bracket is
completely natural, in the sense that it does not depend on any additional structure
for its definition, but it has rather complicated properties. In particular,
it does not satisfy the Leibniz rule with respect to multiplication by functions
or the Jacobi identity. The ``defects'' in both cases are differentials of
certain expressions depending on the bracket and the bilinear form; hence they
disappear upon restriction to a Dirac subbundle. A Dirac subbundle transverse
to \( T^{*}M \) is the graph of a 2-form \( \omega  \), whereas one transverse
to \( TM \) is the graph of a bivector field \( \pi  \); the integrability
condition in this case reduces to the familiar \( d\omega =0 \) (resp. \( [\pi ,\pi ]=0 \)).
Dirac structures, as well as the Courant bracket above, were generalized in
the context of formal variational calculus by Dorfman \cite{Dorf}. 

The nature of the Courant bracket itself remained unclear until several years
later when it was observed by Liu, Weinstein and Xu \cite{LWX1} that \( TM\oplus T^{*}M \)
endowed with the Courant bracket plays the role of a ``double'' object, in
the sense of Drinfeld \cite{Dr}, for a pair of Lie algebroids over \( M \).
\emph{Lie algebroids} are structures on vector bundles that combine the features
of both Lie algebras and the tangent bundle, and include foliations, Poisson
manifolds, Lie group actions, Dirac structures and principal bundles as special
cases. Many differential-geometric and Lie-theoretic constructions carry over
to Lie algebroids. For example, a pair of Lie algebras on dual vector spaces
is called a \emph{Lie bialgebra} if a certain compatibility condition between
them is satisfied. Lie bialgebras are linearizations of Poisson-Lie groups and
semi-classical limits of quantum groups; they also provide a tool for generating
classical integrable systems \cite{ChPr} \cite{Dr}. Likewise, one defines
a \emph{Lie bialgebroid} to be a pair of Lie algebroid structures on dual vector
bundles satisfying a compatibility condition. Lie bialgebroids were first introduced
by Mackenzie and Xu \cite{MacXu} as linearizations of Poisson groupoids. Examples
of Lie bialgebroids for which neither of the Lie algebroid structures is trivial
include Lie bialgebras, Poisson manifolds and Poisson-Nijenhuis manifolds \cite{KS4};
Lie bialgebroids were recently found to be the geometric structure behind the
Classical Dynamical Yang-Baxter equation \cite{EtVar} \cite{KSB}. 

A very useful tool for studying Lie bialgebras is the \emph{Drinfeld double},
which is the Lie algebra structure on the direct sum of the two dual Lie algebras
constituting the bialgebra, uniquely characterized by the requirement that the
two Lie algebras be subalgebras and that the canonical inner product be ad-invariant.
In fact, to find Lie bialgebras, one looks for so-called \emph{Manin triples}:
a Lie algebra with an invariant inner product, together with a pair of complementary
isotropic subalgebras\emph{.} Unfortunately, when one tries to construct a Drinfeld
double for a Lie bialgebroid, it quickly becomes clear that it cannot be a Lie
algebroid if it is to satisfy the characterizing property of the double. Instead,
given a pair \( (A,A^{*}) \) of Lie algebroids in duality, Liu, Weinstein and
Xu \cite{LWX1} build a skew-symmetric bracket on sections of the direct sum
\( A\oplus A^{*} \) similar to the Courant bracket above. Then, they prove
that if \( (A,A^{*}) \) is a Lie bialgebroid, \( A\oplus A^{*} \) becomes
a \emph{Courant algebroid}, a notion they define by emulating the properties
of the original Courant bracket; conversely, they show that any Courant algebroid
which admits a pair of transverse Dirac subbundles (maximally isotropic subbundles
whose sections are closed under the bracket) is of this form, thus extending
the theory of Manin triples to Lie bialgebroids.

In this dissertation we solve several of the problems posed in \cite{LWX1}.
First, the properties of a Courant algebroid are rather complicated; in particular,
there are anomalies in the Jacobi identity and the Leibniz rule. We show that
a Courant algebroid is a resolution of a Lie algebra. It is known \cite{BFLS}
that resolutions of Lie algebras inherit the structure of a \emph{strongly homotopy
Lie algebra}, also known as an \( L_{\infty } \)-\emph{algebra} \cite{SHLA}\emph{,}
though in a non-canonical way. We construct an \( L_{\infty } \)-algebra explicitly
out of the Courant algebroid; the anomalies then appear as the structure identities.
This work appeared in \cite{RoyWe}.

Next, it turns out that one can twist the bracket in a Courant algebroid by
adding a symmetric term. The new operation, which we denote by \( \circ  \)
is, in general, not skew-symmetric but all the anomalies disappear. This was
conjectured in \cite{LWX1}, and we supply a proof. Sacrificing skew-symmetry
has proved worthwhile: the equivalent definition of a Courant algebroid we get
is not only much nicer than the old one, but also more natural, as it turns
out. The Jacobi identity in the non skew-symmetric setting looks rather like
a Leibniz rule: it says that \( a\circ \cdot  \) is a derivation of \( \circ  \).
Such structures were studied by Loday, under the name of \emph{Leibniz algebras}
\cite{Loday1}, and by Kosmann-Schwarzbach \cite{KS5}, under the name of \emph{Loday
algebras.} After the modification, the original Courant bracket becomes 
\[
(X+\xi )\circ (Y+\eta )=[X,Y]+(L_{X}\eta -i_{Y}d\xi ).\]
 This is the form used by Dorfman in \cite{Dorf}. Very recently, \v{S}evera
\cite{Sev} showed that this Courant algebroid provides a natural geometric
framework for studying the symmetries of two-dimensional variational problems.
We use the new definition in all that follows.

Next, in what we regard as the most important part of this work, we develop
an alternative approach to the construction of a Drinfeld double for Lie bialgebroids.
It is based on viewing Lie bialgebroids as homological vector fields on supermanifolds.
To each pair of Lie algebroids in duality we associate a pair of odd self-commuting
hamiltonian functions on an even symplectic supermanifold (in fact, a cotangent
bundle) and prove that the compatibility condition for a Lie bialgebroid is
equivalent to the vanishing of the Poisson bracket of these two hamiltonians.
The hamiltonian vector field of the sum is then homological, and we propose
to call this sum the Drinfeld double. This approach was suggested by the work
of Kosmann-Schwarzbach \cite{KS3} who carried it out for Lie bialgebras in
a purely algebraic language, without mentioning supermanifolds. However, supermanifolds
provide a natural framework even in this case; moreover, the general case cannot
be reduced to pure algebra or ``classical'' geometry, so supermanifolds are
unavoidable. 

The advantage of this approach is its clarity and simplicity. Moreover, several
well-known objects in homological algebra arise in this setting. Thus, applying
this construction to the action Lie algebroid associated to a Lie algebra action
on a manifold, we get the classical BRST complex \cite{KoSt}, whereas applying
it to the Lie bialgebroid associated to the canonical linear Poisson structure
on the dual of a Lie algebra yields the Weil algebra \cite{AtBott}. So far
as we know, this is the only ``geometric'' construction of the Weil algebra
to date. 

To recover the Courant algebroid of Liu, Weinstein and Xu, we use the \emph{derived
bracket construction} of Kosmann-Schwarzbach \cite{KS5}: starting with a differential
Leibniz (in particular, Lie) superalgebra, it generates a new Leibniz superalgebra
of the opposite parity. In particular, Poisson and Schouten brackets arise in
this way. That Courant algebroids may also arise in this way was first suggested
by Kosmann-Schwarzbach, who showed, in a private discussion with the author,
that if one considers the differential Lie superalgebra generated by exterior
multiplications by 1-forms, contractions by vector fields and the de Rham differential,
the derived bracket one gets is the original Courant bracket. What we do here
is a ``semiclassical'' version of this, for an arbitrary Lie bialgebroid.
The Lie superalgebra structure is given by the Poisson bracket on the even symplectic
supermanifold, and the differential is the homological hamiltonian vector field,
the Drinfeld double. The derived bracket we get is precisely the (non skew-symmetric)
Courant bracket of \cite{LWX1}. This enables us to give a very simple proof
of the doubling theorem of Liu, Weinstein and Xu mentioned above.\footnote{
When this research was carried out, we learned that the picture of Lie bialgebroids
as a pair of Poisson-commuting hamiltonians on a symplectic supermanifold was
also considered by A. Vaintrob who studied representations of Lie algebroids;
however, the relation with the Courant algebroids was not elucidated. Our work
is completely independent of his. 
} 

Furthermore, using this approach we are also able to generalize the notion of
a \emph{quasi-Lie bialgebra}, introduced by Drinfeld \cite{Dr2} and studied
by Kosmann-Schwarzbach \cite{KS3}, to the Lie algebroid setting simply by adding
cubic terms to our hamiltonian, thus answering another question posed in \cite{LWX1}.
This also gives a Courant algebroid via the derived bracket construction, thus
answering in the affirmative the question of the existence of nontrivial Courant
algebroids which do not come from Lie bialgebroids. As a special case, we look
at exact Courant algebroids recently classified by \v{S}evera \cite{Sev}. The
cubic term in this case is just the closed 3-form whose cohomology class is
the characteristic class of the Courant algebroid.

This thesis is organized as follows. In Chapter \ref{chapter:courant} we recall
the notions of Lie bialgebra, Lie bialgebroid and Courant algebroid and prove
that Courant algebroids can be considered as strongly homotopy Lie algebras;
we then give a new definition of a Courant algebroid based on the non skew-symmetric
operation and prove its equivalence to the old one. 

In Chapter \ref{chapter:supermanifolds} we develop the theory of Lie bialgebroids
and quasi-bialgebroids in terms of even symplectic supermanifolds, give the
derived bracket construction of the Courant algebroid and re-prove the doubling
theorem of Liu, Weinstein and Xu, generalizing it also for quasi-bialgebroids. 

In the final Chapter \ref{chapter:Poisson}, somewhat disjoint from the rest,
we study a one-parameter family of Poisson structures on \( S^{2} \) covariant
with respect to the action of \( SU(2) \) with its standard Poisson-Lie group
structure. We compute the Poisson cohomology of these structures and show, as
an application, that they do not admit rescaling, and also that they are non-trivial
deformations of each other. 

Throughout this dissertation, a manifold will always mean a smooth real manifold,
and all vector spaces, algebras, etc. are over the field of real numbers, unless
otherwise specified. The Einstein summation convention is used consistently.

\chapter{Courant algebroids and strongly homotopy Lie algebras\label{chapter:courant}}

In this chapter we recall the definition of a Courant algebroid first given
in \cite{LWX1} and some of the results obtained therein. We then make the first
step toward explaining the anomalies of Courant algebroids by showing that they
can be considered as strongly homotopy Lie algebras. This is essentially the
content of \cite{RoyWe}. In the last section we propose an equivalent definition
of a Courant algebroid which has the advantage of being anomaly-free (except
for lack of skew-symmetry), and will be useful in what follows. To begin, we
recall the notions of a Lie bialgebra, Lie algebroid and bialgebroid and give
some examples.

\section{Lie bialgebras\label{sec:bialgebras}}

\begin{defn}
\label{def:LieBialgebra} A \emph{Lie bialgebra} is a vector space \( \g  \)
together with a bilinear skew-\\ symmetric map \( \mu =\br :\wedge ^{2}\g \rightarrow \g  \)
(the bracket) and a linear map \( \gamma :\g \rightarrow \wedge ^{2}\g  \)
(the cobracket) such that the following properties are satisfied:
\begin{itemize}
\item \( \g  \) together with \( \br  \) is a Lie algebra;
\item \( \g ^{*} \) together with \( \br _{*}=\gamma ^{*}:\wedge ^{2}\g ^{*}\rightarrow \g ^{*} \)
is a Lie algebra;
\item \( \gamma  \) is a 1-cocycle on the Lie algebra \( (\g ,\mu ) \) with values
in the (exterior square of the) adjoint module \( \wedge ^{2}\g  \), i.e. 
\[
\gamma ([a,b])=ad_{a}\gamma (b)-ad_{b}\gamma (a)\]
 holds for all \( a,b\in \g  \).
\end{itemize}
\end{defn}
One sometimes calls the pair \( (\g ,\g ^{*}) \) a Lie bialgebra with the underlying
structures implicitly understood. Lie bialgebras are the infinitesimal objects
corresponding to \emph{Poisson-Lie groups} (see Appendix); they are also the
semi-classical limits of \emph{quantum groups} (see \cite{ChPr} for a thorough
treatment and numerous examples).

\begin{defn}
Given a Lie bialgebra \( (\g ,\g ^{*}) \), its \emph{double} (or \emph{Drinfeld
double)} is the vector space direct sum \( \g \oplus \g ^{*} \) together with
the bracket 
\begin{equation}
\label{eqn:doubleofbialgebra}
[X+\xi ,Y+\eta ]=([X,Y]+ad^{*}_{\xi }Y-ad^{*}_{\eta }X)+(ad^{*}_{X}\eta -ad^{*}_{Y}\xi +[\xi ,\eta ]_{*})
\end{equation}
 
\end{defn}
This bracket is completely characterized by the property that both \( \g  \)
and \( \g ^{*} \) be subalgebras of \( \g \oplus \g ^{*} \) and that the canonical
inner product 
\begin{equation}
\label{eqn:CanBform}
\bform {X+\xi }{Y+\eta }=\xi (Y)+\eta (X)
\end{equation}
 be \( ad \)-invariant; it satisfies the Jacobi identity if \( (\g ,\g ^{*}) \)
is a Lie bialgebra. In fact, the notion of a Lie bialgebra is equivalent to
that of a \emph{Manin triple} which is a triple \( (\p ,\p _{+},\p _{-}) \),
where \( \p  \) is a Lie algebra with an invariant symmetric bilinear form,
and \( \p _{+} \) and \( \p _{-} \) are complementary isotropic subalgebras.
Manin triples abound in nature: for example, every complex semisimple Lie algebra
gives rise to a Manin triple via the Iwasawa decomposition (see \cite{LuWe}).

\section{Lie algebroids and bialgebroids\label{sec:algebroids}}

\begin{defn}
\label{def:liealg} A \emph{Lie algebroid} is a vector bundle \( A\rightarrow M \)
together with a Lie algebra bracket \( \br _{A} \) on the space of sections
\( \Gamma (A) \) and a bundle map \( a:A\rightarrow TM \), called the \emph{anchor},
satisfying the following conditions: 
\begin{enumerate}
\item For any \( X,Y\in \Gamma (A) \), \( a[X,Y]_{A}=[aX,aY] \)
\item For any \( X,Y\in \Gamma (A),\; f\in \cinf (M) \), \( [X,fY]_{A}=f[X,Y]_{A}+(a(X)f)Y \)
\end{enumerate}
\end{defn}
In other words, the sections of the bundle act on smooth functions by derivations
via the anchor in such a way that brackets act as commutators, and the behavior
of the bracket with respect to multiplication by functions is governed by the
Leibniz rule. Thus, Lie algebroids are a straightforward generalization of the
tangent bundle. They are also the infinitesimal objects corresponding to Lie
groupoids \cite{Mac}; when the base manifold is a point, a Lie groupoid reduces
to a Lie group, while a Lie algebroid is just a Lie algebra.

A Lie algebroid structure on \( A\rightarrow M \) gives rise to the following
structures, dual to one another. The \emph{generalized Schouten bracket} is
defined as the unique extension \( \br _{A} \) of the Lie bracket on \( \Gamma (A) \)
and the action of \( \Gamma (A) \) on functions to \( \Gamma (\bigwedge ^{*}A) \)
such that:

\begin{enumerate}
\item \( [X,Y]_{A}=-(-1)^{pq}[Y,X]_{A} \), for \( X\in \Gamma (\bigwedge ^{p+1}A) \),
\( Y\in \Gamma (\bigwedge ^{q+1}A) \), 
\item \( [X,f]_{A}=a(X)f \) for \( X\in \Gamma (A) \), \( f\in \cinf (M) \), 
\item For \( X\in \Gamma (\bigwedge ^{p+1}A) \), \( [X,\cdot ]_{A} \) is a derivation
of degree \( p \) of the exterior multiplication on \( \Gamma (\bigwedge ^{*}A) \). 
\end{enumerate}
One checks that this bracket satisfies the graded Jacobi identity with respect
to the grading shifted down by one, and the resulting structure is a type of
graded Poisson algebra called a \emph{Gerstenhaber algebra}.

Dually, one gets a derivation \( d_{A} \) of degree 1 on the graded commutative
algebra \( \Gamma (\bigwedge ^{*}A^{*}) \), defined by a formula identical
to the Cartan formula for the de Rham differential: 
\[
\begin{array}{ccc}
d_{A}\omega (X_{0},\ldots ,X_{p}) & = & \sum ^{p}_{i=0}(-1)^{i}a(X_{i})(\omega (X_{0},\ldots ,\hat{X}_{i},\ldots ,X_{p}))+\\
 & + & \sum _{0\leq i<j\leq p}(-1)^{i+j}\omega ([X_{i},X_{j}]_{A},X_{0},\ldots ,\hat{X}_{i},\ldots ,\hat{X}_{j},\ldots ,X_{p}),
\end{array}\]
 where \( \omega \in \Gamma (\bigwedge ^{p}A^{*}) \), and satisfying \( d^{2}_{A}=0 \).
The space \( \Gamma (\wedge ^{*}A^{*}) \) thereby acquires the structure of
a differential graded commutative algebra. \( d_{A} \) is uniquely determined
by its action on \( \cinf (M) \) and \( \Gamma (A^{*}) \): 
\begin{equation}
\label{eqn:CartanFla}
\begin{array}{rcl}
d_{A}f(X) & = & a(X)f\\
d_{A}\xi (X,Y) & = & a(X)\xi (Y)-a(Y)\xi (X)-\xi ([X,Y]_{A})
\end{array}
\end{equation}
 It is clear that, conversely, the Lie algebroid structure is completely determined
by either \( d_{A} \) (all the structural identities are encoded in \( d_{A}^{2}=0 \)),
or the generalized Schouten bracket \( \br _{A} \).

Many notions of the usual calculus on manifolds carry over without change to
Lie algebroids. In particular, for every \( X\in \Gamma (A) \) there is a contraction
(interior derivative) operator \( i_{X} \) acting on \( \Gamma (\bigwedge ^{*}A^{*}) \)
by derivations of degree \( -1 \), and the ``Lie derivative'' operator \( L_{X}^{A}=[d_{A},i_{X}] \)
acting by derivations of degree \( 0 \) (here \( \br  \) denotes the supercommutator).
These derivations satisfy the usual (super)commutation relations: 
\begin{equation}
\label{eqn:canonical_{r}elations}
\begin{array}{ccc}
[d_{A},d_{A}]=0, & [d_{A},L_{X}^{A}]=0, & [d_{A},i_{X}]=L_{X}^{A},\\
{[L_{X}^{A},L_{Y}^{A}]}=L_{[X,Y]_{A}}^{A}, & [i_{X},i_{Y}]=0, & [L_{X}^{A},i_{Y}]=i_{[X,Y]_{A}}
\end{array}
\end{equation}
 Now suppose that we are given a pair \( (A,A^{*}) \) of Lie algebroids over
\( M \) which are in duality as vector bundles. Then the Lie algebroid structure
of \( A \) induces a Schouten bracket on \( \Gamma (\bigwedge ^{*}A) \) and
a differential \( d_{A} \) on \( \Gamma (\bigwedge ^{*}A^{*}) \); on the other
hand, from \( A^{*} \) we get a Schouten bracket on \( \Gamma (\bigwedge ^{*}A^{*}) \)
and a differential \( d_{A^{*}} \) on \( \Gamma (\bigwedge ^{*}A) \).

\begin{defn}
\label{def:liebialg} A pair \( (A,A^{*}) \) of Lie algebroids
in duality is called a \emph{Lie bialgebroid} if the induced differential \( d_{A} \)
is a derivation of the Schouten bracket \( \br _{A^{*}} \) on \( \Gamma (\wedge ^{*}A^{*}) \). 
\end{defn}
Thus, Lie bialgebroids correspond to \emph{differential Gerstenhaber algebras}
\cite{KS1}. The notion of a Lie bialgebroid is due to Mackenzie and Xu \cite{MacXu}
who studied them and the corresponding global objects, Poisson groupoids (although the definition we quoted is
an equivalent one from \cite{KS1}). It
can be shown that this notion is self-dual, i.e. if \( (A,A^{*}) \) is a Lie
bialgebroid, so is \( (A^{*},A) \) (Corollary\ref{cor:self-dual} below).

\begin{rem}
Any Lie algebroid is a Lie bialgebroid with the zero anchor and bracket on the
dual bundle. 
\end{rem}
\begin{example}
\label{eg:LieAlg-standard} Let \( M \) be a manifold. Then its tangent bundle
\( TM \) is a Lie algebroid whose bracket is the Jacobi-Lie bracket of vector
fields, and the anchor is \( \rho =\textrm{Id}:TM\rightarrow TM \). The corresponding
extended bracket is the (original) Schouten bracket of multivector fields, while
the differential is just the de Rham differential. 
\end{example}

\begin{example}
\label{eg:LieAlg-action} Consider a (right) action of a Lie algebra \( \g  \)
on a manifold \( M \), i.e. a Lie algebra homomorphism \( \rho :\g \rightarrow \X (M) \).
This gives rise to a Lie algebroid structure on the trivial bundle \( M\times \g \rightarrow M \)
whose anchor is given on constant sections by \( \rho  \) and extended to all
sections by linearity over \( \cinf (M) \), while the bracket of constant sections
is just the bracket in \( \g  \) taken pointwise, which is then extended to
all sections by the Leibniz rule. This Lie algebroid is called the \emph{action
Lie algebroid} associated to \( \rho  \). If \( \rho  \) is a left action
(a Lie algebra \emph{anti}homomorphism), then we must take \( -\rho  \) as
the anchor. 
\end{example}

\begin{example}
\label{eg:LieBialgebra} If the base manifold \( M \) is a point, a Lie bialgebroid
\( (A,A^{*}) \) over \( M \) is just a Lie bialgebra (we shall see later that
Definition \ref{def:liebialg} is equivalent to Definition \ref{def:LieBialgebra}
in this case). 
\end{example}

\begin{example}
\label{eg:LieBialg-Poisson} Let \( M \) be a Poisson manifold with Poisson
tensor \( \pi  \) and the corresponding bundle map \( \tilde{\pi }:T^{*}M\rightarrow TM \)
given by \( \langle \tilde{\pi }\alpha ,\beta \rangle =\pi (\alpha ,\beta ) \).
Let \( A=TM \), the tangent bundle Lie algebroid, \( A^{*}=T^{*}M \) with
anchor \( \tilde{\pi } \) and the bracket of 1-forms given by the \emph{Koszul
bracket:}
\begin{equation}
[\alpha ,\beta ]_{A^{*}}={\mathcal{L}}_{\tilde{\pi }\alpha }\beta -{\mathcal{L}}_{\tilde{\pi }\beta }\alpha -d(\pi (\alpha ,\beta ))
\end{equation}
 Then \( d_{A} \) is the usual de Rham differential of forms, \( d_{A^{*}}=[\pi ,\cdot ]_{A} \),
where \( \br _{A} \) is the Schouten bracket, and it is straightforward to
verify that \( (A^{*},A) \) is a Lie bialgebroid. 
\end{example}

Detailed discussion and more examples of Lie bialgebroids and Gerstenhaber algebras
from geometry and physics can be found in \cite{KS1},\cite{KS2} and \cite{KS4}.

\section{Courant algebroids\label{sec:courant}}

\begin{defn}
\label{def:jac} Given a bilinear, skew-symmetric operation \( \br  \) on a
vector space \( V \), its \emph{Jacobiator} \( J \) is the trilinear operator
on \( V \): 
\[
J(e_{1},e_{2},e_{3})=[[e_{1},e_{2}],e_{3}]+[[e_{2},e_{3}],e_{1}]+[[e_{3},e_{1}],e_{2}],\]
 \( e_{1},e_{2},e_{3}\in V \). 
\end{defn}
The Jacobiator is obviously skew-symmetric. Of course, in a Lie algebra \( J\equiv 0 \).

\begin{defn}
\label{def:quasi-algebroid} A \emph{Courant algebroid} is a vector bundle \( E\lon M \)
equipped with a nondegenerate symmetric bilinear form \( \bform {\cdot }{\cdot } \)
on the bundle, a skew-symmetric bracket \( \br  \) on \( \Gamma (E) \), and
a bundle map \( \rho :E\lon TM \) such that the following properties are satisfied: 
\begin{enumerate}
\item For any \( e_{1},e_{2},e_{3}\in \Gamma (E) \), \( J(e_{1},e_{2},e_{3})={\mathcal{D}}T(e_{1},e_{2},e_{3}); \)
\item for any \( e_{1},e_{2}\in \Gamma (E) \), \( \rho [e_{1},e_{2}]=[\rho e_{1},\rho e_{2}]; \)
\item for any \( e_{1},e_{2}\in \Gamma (E) \) and \( f\in \cinf (M) \), \( [e_{1},fe_{2}]=f[e_{1},e_{2}]+(\rho (e_{1})f)e_{2}-\half \langle e_{1},e_{2}\rangle {\mathcal{D}}f; \)
\item \( \rho \smalcirc {\mathcal{D}}=0 \), i.e., for any \( f,g\in \cinf (M) \),
\( \langle {\mathcal{D}}f,{\mathcal{D}}g\rangle =0 \); 
\item for any \( e,h_{1},h_{2}\in \Gamma (E) \), \( \rho (e)\langle h_{1},h_{2}\rangle =\langle [e,h_{1}]+\half {\mathcal{D}}\langle e,h_{1}\rangle ,h_{2}\rangle +\langle h_{1},[e,h_{2}]+\half {\mathcal{D}}\langle e,h_{2}\rangle \rangle  \), 
\end{enumerate}
where \( T(e_{1},e_{2},e_{3}) \) is the function on the base \( M \) defined
by: 
\begin{equation}
\label{eq:T0}
T(e_{1},e_{2},e_{3})=\sixth \langle [e_{1},e_{2}],e_{3}\rangle +c.p.,
\end{equation}
 (``c.p.'' denotes the cyclic permutations of the \( e_{i} \)'s) and \( {\mathcal{D}}:\cinf (M)\lon \Gamma (E) \)
is the map defined by \( {\mathcal{D}}=\rho ^{*}d \), where \( E \) is identified
with \( E^{*} \) by the bilinear form and \( d \) is the deRham differential.
In other words, 
\begin{equation}
\label{eq:D}
\langle {\mathcal{D}}f,e\rangle =\rho (e)f.
\end{equation}

\end{defn}
\begin{note*}
In our convention, the bilinear form \( \bform {\cdot }{\cdot } \) is two times
the one in \cite{LWX1}. 
\end{note*}
In a Courant algebroid \( E \), a \emph{Dirac structure}, or \emph{Dirac subbundle},
is a subbundle \( L \) that is maximally isotropic under \( \bform {\cdot }{\cdot } \)
and whose sections are closed under \( \br  \). It is immediate from the definition
that a Dirac subbundle is a Lie algebroid under the restrictions of the bracket
and anchor.

Suppose now that both \( A \) and \( A^{*} \) are Lie algebroids over the
base manifold \( M \), with anchors \( a \) and \( a_{*} \) respectively.
Let \( E \) denote their vector bundle direct sum: \( E=A\oplus A^{*} \).
On \( E \), there exist two natural nondegenerate bilinear forms, one symmetric
and another antisymmetric:

\begin{equation}
\label{eq:pairing}
(X_{1}+\xi _{1},X_{2}+\xi _{2})_{\pm }=(\langle \xi _{1},X_{2}\rangle \pm \langle \xi _{2},X_{1}\rangle ).
\end{equation}

On \( \Gamma (E) \), we introduce a bracket by

\begin{equation}
\label{eq:double}
\begin{array}{rcl}
[e_{1},e_{2}] & = & ([X_{1},X_{2}]_{A}+L^{A^{*}}_{\xi _{1}}X_{2}-L^{A^{*}}_{\xi _{2}}X_{1}-\half d_{A^{*}}(e_{1},e_{2})_{-})+\\
 & + & ([\xi _{1},\xi _{2}]_{A^{*}}+L^{A}_{X_{1}}\xi _{2}-L^{A}_{X_{2}}\xi _{1}+\half d_{A}(e_{1},e_{2})_{-}),
\end{array}
\end{equation}
 where \( e_{1}=X_{1}+\xi _{1} \) and \( e_{2}=X_{2}+\xi _{2} \).\\
\\

Finally, we let \( \rho :E\lon TM \) be the bundle map defined by \( \rho =a+a_{*} \).
That is, 
\begin{equation}
\rho (X+\xi )=a(X)+a_{*}(\xi ),\, \, \, \, \forall X\in \Gamma (A)\mbox {and}\xi \in \Gamma (A^{*})
\end{equation}
 It is easy to see that in this case the operator \( {\mathcal{D}} \) as defined
by Equation (\ref{eq:D}) is given by 
\[
{\mathcal{D}}=d_{A^{*}}+d_{A}\]
 The following results, which we quote from \cite{LWX1}, show that the notion
of Courant algebroid permits us to generalize the double construction to Lie
bialgebroids:

\begin{thm}
\label{thm:bialg->courant} If \( (A,A^{*}) \) is a Lie bialgebroid, then \( E=A\oplus A^{*} \)
together with \\ \( (\br ,\rho ,(\cdot ,\cdot )_{+}) \) is a Courant algebroid. 
\begin{thm}
\label{thm:courant->bialg} In a Courant algebroid \( (E,\rho ,\br ,\bform {\cdot }{\cdot }) \),
suppose that \( L_{1} \) and \( L_{2} \) are Dirac subbundles transversal
to each other, i.e., \( E=L_{1}\oplus L_{2} \). Then, \( (L_{1},L_{2}) \)
is a Lie bialgebroid, where \( L_{2} \) is considered as the dual bundle of
\( L_{1} \) under the pairing \( \bform {\cdot }{\cdot } \). 
\end{thm}
\end{thm}
An immediate consequence of the theorems above is the following duality property
of Lie bialgebroids, which was first proved in \cite{MacXu} and then by Kosmann-Schwarzbach
\cite{KS1} using a simpler method.

\begin{cor}
\label{cor:self-dual} If \( (A,A^{*}) \) is a Lie bialgebroid, so is \( (A^{*},A) \). 
\end{cor}
The theorems above are proved in \cite{LWX1} by rather laborious computations;
in the next chapter we shall give a new, simple proof of Theorem \ref{thm:bialg->courant}
and Corollary \ref{cor:self-dual}.

\begin{example}
Given a manifold \( M \), consider \( TM \) with its standard Lie algebroid
structure and \( T^{*}M \) with zero anchor and bracket. Then \( (TM,T^{*}M) \)
is a Lie bialgebroid, and the double bracket (\ref{eq:double}) reduces to 
\[
\textstyle [X_{1}+\xi _{1},X_{2}+\xi _{2}]=[X_{1},X_{2}]+(L_{X_{1}}\xi _{2}-L_{X_{2}}\xi _{1}+d({{\frac{1}{2}}}(\xi _{1}(X_{2})-\xi _{2}(X_{1}))).\]
 This is the bracket originally introduced by Courant in \cite{Cou}. The anchor
\( \rho  \) in this case is the projection to \( TM \), and \( \DD =d \),
the deRham differential. 
\begin{example}
When \( M \) is a point, \( (A,A^{*}) \) is a Lie bialgebra and the bracket
\ref{eq:double} on \( E \) becomes the famous Drinfeld double bracket. 
\end{example}
\end{example}

\section{Strongly homotopy Lie algebras and Courant algebroids\label{sec:SHLA}}

Let \( V \) be a graded vector space. Let \( T(V) \) denote the tensor algebra
of \( V \) in the category of graded vector spaces, and let \( \bw (V) \)
denote its exterior algebra in the same category; i.e. \( \bw (V)=T(V)/<v\otimes w+(-1)^{{\tilde{v}}{\tilde{w}}}w\otimes v> \),
where \( {\tilde{v}} \) denotes the degree of \( v \). \( \bw (V) \) has
a natural Hopf algebra structure with the comultiplication \( \Delta :\bw (V)\rightarrow \bw (V)\otimes \bw (V) \)
uniquely defined by the requirement that the elements of \( V \) be primitive
(i.e. \( \Delta v=1\otimes v+v\otimes 1 \) for \( v\in V \)) and that \( \Delta  \)
be a homomorphism of algebras (see \cite{SHLA} for details).

\begin{defn}
\label{def:SHLA} A \emph{strongly homotopy Lie algebra} (SHLA, \( L_{\infty } \)-algebra)
is a graded vector space \( V \) together with a collection of linear maps
\( l_{k}:\bw ^{k}V\rightarrow V \) of degree \( k-2 \), \( k\geq 1 \), satisfying
the following relation for each \( n\geq 1 \) and for all homogeneous \( x_{1},\dots ,x_{n}\in V \):

\begin{equation}
\label{eqn:SHLA}
\sum _{i+j=n+1}(-1)^{i(j-1)}\sum _{\sigma }(-1)^{\sigma }\epsilon (\sigma )l_{j}(l_{i}(x_{\sigma (1)},\dots ,x_{\sigma (i)}),x_{\sigma (i+1)},\dots ,x_{\sigma (n)})=0,
\end{equation}
 where \( \sigma  \) runs over all \( (i,n-i) \)-unshuffles (permutations
satisfying \( \sigma (1)<\dots <\sigma (i) \) and \( \sigma (i+1)<\dots <\sigma (n) \))
with \( i\geq 1 \), and \( \epsilon (\sigma ) \) is the \emph{Koszul sign}
(arising from the fundamental convention of supermathematics that a minus sign
is introduced whenever two consecutive odd elements are permuted) . 
\end{defn}
For \( n=1 \) this means simply that \( l_{1} \) is a differential on \( V \);
for \( n=2 \), \( l_{2} \) is a superbracket on \( V \) of which \( l_{1} \)
is a derivation (equivalently, \( l_{2}:\bw ^{2}(V)\rightarrow V \) is a chain
map of complexes); \( n=3 \) gives the Jacobi identity for \( l_{2} \) satisfied
up to chain homotopy given by \( l_{3} \), and higher \( l_{k} \)'s can be
interpreted as higher homotopies. The algebraic theory of \( L_{\infty } \)-algebras
is studied in \cite{HinSch} and \cite{SHLA} .

We shall write the equation (\ref{eqn:SHLA}) in the more succinct equivalent
form: 
\begin{equation}
\label{eqn:SHLA'}
\sum _{i+j=n+1}(-1)^{i(j-1)}l_{j}l_{i}=0,
\end{equation}
 where we have extended each \( l_{i} \) to all of \( \bw (V) \) as a coderivation
of the coalgebra structure on \( \bw (V) \). This accounts for the permutations
and signs in (\ref{eqn:SHLA}).

We are interested in \( L_{\infty } \)-algebras for the following reason: it
is shown in \cite{BFLS} that, given a resolution \( (X_{*},d) \) of a vector
space \( H \) (graded or not), any Lie algebra structure on \( H \) can be
lifted to an \( L_{\infty } \)-algebra structure on the total resolution space
\( X \) with \( l_{1}=d \). The starting point of this construction is the
observation that Lie brackets on \( H \) correspond to bilinear skew-symmetric
brackets \( \br  \) on \( X_{0} \) for which the boundaries form an ideal
and the Jacobi identity is satisfied up to a boundary. This correspondence is
in no way unique or canonical, as it requires a choice of a homotopy inverse
to the quasi-isomorphism \( (X_{*},d)\rightarrow (H,0) \)). But it is this
bracket \( \br  \) on \( X_{0} \) that provides the starting point for constructing
the SHLA structure on \( X \), hence, if it is given, no choice is required
at this stage, and we need never mention \( H \). We shall presently see that
with Courant algebroids we are in precisely this situation.

Let \( E \) be a Courant algebroid over a manifold \( M \). We know from the
definition that the Courant bracket on \( \Gamma (E) \) satisfies Jacobi up
to a \( \DD  \)-exact term. It turns out that, moreover, \( Im(\DD ) \) is
an ideal in \( \Gamma (E) \) with respect to the bracket. More precisely, the
following identity holds:

\begin{lem}
\label{lemma:ideal1}For any \( e\in \Gamma (E) \), \( f\in C^{\infty }(M) \)
one has 
\[
[e,\DD f]=\half \DD \bform {e}{\DD f}\]
 
\end{lem}
\begin{proof}
Use axiom 5 in the definition of Courant algebroid with \( e={\mathcal{D}}f \)
and arbitrary \( h_{1} \) and \( h_{2} \), and then cyclically permute \( e \),
\( h_{1} \) and \( h_{2} \):

\begin{eqnarray*}
\rho (\DD f)\bform {h_{1}}{h_{2}} & = & \bform {[\DD f,h_{1}]+\half \DD \bform {\DD f}{h_{1}}}{h_{2}}+\bform {h_{1}}{[\DD f,h_{2}]+\half \DD \bform {\DD f}{h_{2}}}\\
\rho (h_{1})\bform {h_{2}}{\DD f} & = & \bform {[h_{1},h_{2}]+\half \DD \bform {h_{1}}{h_{2}}}{\DD f}+\bform {h_{2}}{[h_{1},\DD f]+\half \DD \bform {h_{1}}{\DD f}}\\
\rho (h_{2})\bform {\DD f}{h_{1}} & = & \bform {[h_{2},\DD f]+\half \DD \bform {h_{2}}{\DD f}}{h_{1}}+\bform {\DD f}{[h_{2},h_{1}]+\half \DD \bform {h_{2}}{h_{1}}}.
\end{eqnarray*}
 Now add the first two identities and subtract the third. Using Courant algebroid
axioms 2, 4 and the definition of \( \DD  \), we get: 
\[
\rho ([h_{1},h_{2}])f=\bform {\DD f}{2[h_{1},h_{2}]}+\bform {h_{1}}{2[\DD f,h_{2}]}+\bform {h_{2}}{\DD \bform {\DD f}{h_{1}}}.\]
 Using the definition of \( \DD  \) again, we can rewrite this as: 
\begin{eqnarray*}
0 & = & \rho ([h_{1},h_{2}])f+\bform {h_{1}}{2[\DD f,h_{2}]}+\rho (h_{2})\bform {h_{1}}{\DD f}=\\
 & = & \rho ([h_{1},h_{2}])f+\bform {h_{1}}{2[\DD f,h_{2}]}+\rho (h_{2})(\rho (h_{1})f)=\\
 & = & \rho (h_{1})(\rho (h_{2})f)+\bform {h_{1}}{2[\DD f,h_{2}]}=\\
 & = & \bform {h_{1}}{\DD (\rho (h_{2})f)+2[\DD f,h_{2}]}=\\
 & = & \bform {h_{1}}{2(\half \DD \bform {h_{2}}{\DD f}-[h_{2},\DD f])}.
\end{eqnarray*}
 The statement follows from the nondegeneracy of \( \bform {\cdot }{\cdot } \). 
\end{proof}
It will follow that we can extend the Courant bracket to an \( L_{\infty } \)-structure
on the total space of the following resolution of \( H=\coker \DD  \): 
\begin{equation}
\label{eqn:res}
\cdots \lon 0\lon X_{2}\stackrel{d_{2}}{\lon }X_{1}\stackrel{d_{1}}{\lon }X_{0}\lon H\lon 0,
\end{equation}
 where \( X_{0}=\Gamma (E) \), \( X_{1}=\cinf (M) \), \( X_{2}=\ker \DD  \),
\( d_{1}=\DD  \) and \( d_{2} \) is the inclusion \( \iota :\ker \DD \hookrightarrow \cinf (M) \).
Remarkably, it turns out that, owing to the properties of Courant algebroids,
the choices in the extension procedure can be made in a natural and simple way.

Let us fix some notation: we will denote elements of \( X_{0} \) by \( e \),
elements of \( X_{1} \) by \( f \) or \( g \), and elements of \( X_{2} \)
by \( c \).

\begin{thm}
\label{thm:CA->SHLA}A Courant algebroid structure on a vector bundle \( E\lon M \)
gives rise naturally to a SHLA structure on the total space \( X \) of (\ref{eqn:res})
with \( l_{1}=d \) and the higher structure maps given by the following explicit
formulas: 
\[
\begin{array}{lccl}
l_{2}(e_{1}\wedge e_{2}) & = & [e_{1},e_{2}] & \textrm{in degree }0\\
l_{2}(e\wedge f) & = & \half \bform {e}{\DD f} & \textrm{in degree }1\\
l_{2} & = & 0 & \textrm{in degree }>1\\
l_{3}(e_{1}\wedge e_{2}\wedge e_{3}) & = & -T(e_{1},e_{2},e_{3}) & \textrm{in degree }0\\
l_{3} & = & 0 & \textrm{in degree }>0\\
l_{n} & = & 0 & \textrm{for }n>3
\end{array}\]

\end{thm}
\begin{proof}
\noindent Starting with the Courant bracket on \( X_{0} \), we shall, following
\cite{BFLS}, extend it to an \( l_{2} \) on all of \( X \) satisfying (\ref{eqn:SHLA'})
for \( n=2 \). The extension will proceed, essentially, by induction on the
degree of the argument: for each degree \( l_{2} \) will be a primitive of
a certain cycle depending on the values of \( l_{2} \) on elements of lower
degree. Higher \( l_{k} \)'s will be introduced and extended in a similar fashion,
as primitives of cycles (using the acyclicity of (\ref{eqn:res})). The main
work will consist in calculating these cycles, in particular, showing that most
of them vanish; these computations are mostly relegated to the technical lemmas
of the next section. \emph{Step 1: \( n=2 \)}. In degree 0, we are given \( l_{2}(e_{1}\w e_{2})=[e_{1},e_{2}] \).
Consider now an element \( e\w f \) of degree 1. Then \( l_{2}l_{1}(e\w f)\in X_{0} \)
is defined and is, in fact, a boundary by Lemma \ref{lemma:ideal1}: 
\[
l_{2}l_{1}(e\w f)=l_{2}(l_{1}e\w f+e\w l_{1}f)=[e,\DD f]=\half \DD \bform {e}{\DD f},\]
 so we set \( l_{2}(e\w f)=\half \bform {e}{\DD f} \) so that the SHLA identity
(\ref{eqn:SHLA'}) for \( n=2 \), 
\begin{equation}
\label{eqn:SHLA2}
l_{1}l_{2}-l_{2}l_{1}=0,
\end{equation}
 holds in degree 1. Now, \( \bw ^{2}(X)_{2} \) is spanned by elements of the
form \( f\w g \) or \( c\w e \). As above, \( l_{2}l_{1} \) is defined on
elements of degree 2, and is, in fact, a cycle (cf. \cite{BFLS}). We have 
\[
l_{2}l_{1}(f\w g)=l_{2}(l_{1}f\w g-f\w l_{1}g)=l_{2}(\DD f\w g-f\w \DD g)=\half (\bform {\DD f}{\DD g}+\bform {\DD g}{\DD f})=0\]
 by Courant algebroid axiom 4, whereas 
\[
l_{2}l_{1}(c\w e)=l_{2}(l_{1}c\w e+c\w l_{1}e)=l_{2}(\iota c\w e)=-\half \bform {e}{\DD \iota c}=0,\]
 so we set \( l_{2}(f\w g)=l_{2}(c\w e)=0 \). Now observe that, since \( l_{2}=0 \)
in degree 2, we can define \( l_{2} \) to be zero on elements of degree higher
than 2 as well and still have (\ref{eqn:SHLA2}). We have thus defined an \( l_{2} \)
that satisfies (\ref{eqn:SHLA2}) by construction. \emph{Step 2: \( n=3 \)}.
In degree 0, by Courant algebroid axiom 1 we have 
\[
l_{2}l_{2}(e_{1}\w e_{2}\w e_{3})=J(e_{1},e_{2},e_{3})=\DD T(e_{1},e_{2},e_{3}),\]
 where \( J \) is the Jacobiator. So we set \( l_{3}(e_{1}\w e_{2}\w e_{3})=-T(e_{1},e_{2},e_{3}) \),
so that the homotopy Jacobi identity identity (\ref{eqn:SHLA'}) for \( n=3 \),
\begin{equation}
\label{eqn:SHLA3}
l_{1}l_{3}+l_{2}l_{2}+l_{3}l_{1}=0,
\end{equation}
 holds on \( \bw ^{3}(X)_{0} \) (as \( l_{1}(X_{0})=0 \)). Consider now an
element \( e_{1}\w e_{2}\w f\in \bw ^{3}(X)_{1} \) . The expression \( (l_{2}l_{2}+l_{3}l_{1})(e_{1}\w e_{2}\w f) \)
is defined and is a cycle in \( X_{1} \) (cf. \cite{BFLS}), hence we can define
\( l_{3}(e_{1}\w e_{2}\w f) \) to be some primitive of this cycle, so that
(\ref{eqn:SHLA3}) holds. But in our particular situation we in fact have (see
the next section for a proof): 
\begin{lem}
\label{lemma:t1} \( (l_{2}l_{2}+l_{3}l_{1})(e_{1}\w e_{2}\w f)=0 \) \( \forall e_{1},e_{2},f \). 
\end{lem}
\noindent Therefore, we can define \( l_{3}(e_{1}\w e_{2}\w f)=0 \). Now observe
that on elements of degree \( >1 \) \( l_{3} \) has to be 0 because deg\( (l_{3})=1 \),
whereas \( X_{k}=0 \) for \( k>2 \). We now have \( l_{3} \) defined on all
of \( \bw ^{3}(X) \) and satisfying (\ref{eqn:SHLA3}) by construction. \emph{Step
3: \( n=4 \) and higher}. Proceeding in a similar fashion, we look at the expression\\
 \( (l_{3}l_{2}-l_{2}l_{3})(e_{1}\w e_{2}\w e_{3}\w e_{4}) \) (always a cycle
in \( X_{1} \)) and define \( l_{4}(e_{1}\w e_{2}\w e_{3}\w e_{4}) \) to be
its primitive in \( X_{2} \), so as to satisfy (\ref{eqn:SHLA'}). However,
it turns out that (see the next section for a proof) 
\begin{lem}
\label{lemma:t2} \( (l_{3}l_{2}-l_{2}l_{3})(e_{1}\w e_{2}\w e_{3}\w e_{4})=0 \)
\( \forall e_{1},e_{2},e_{3},e_{4} \). 
\end{lem}
Hence we can set \( l_{4}(e_{1}\w e_{2}\w e_{3}\w e_{4})=0 \) and observe that
\( l_{4} \) has to vanish on elements of degree \( >0 \) as deg\( (l_{4})=2 \),
while \( X_{k}=0 \) for \( k>2 \). By similar degree counting, all \( l_{n} \),
\( n>4 \), have to vanish identically. This finishes the proof modulo Lemmas
\ref{lemma:t1} and \ref{lemma:t2}. 
\end{proof}
\begin{rem}
If the base \( M \) is a point, a Courant algebroid reduces to a Lie algebra
\( \p  \) with an invariant inner product; however, even though the differential
\( \DD  \) is trivial in this case and all the anomalies vanish, the homotopy
Lie algebra we get is not ``just a Lie algebra'': in addition to the Lie algebra
bracket there is also a trilinear operation \( T \), the structure tensor of
the Lie algebra: 
\[
T(X,Y,Z)=\half \bform {[X,Y]}{Z}\]
 for \( X,Y,Z\in \p  \).
\end{rem}

\section{Proofs of technical lemmas\label{sec:lemmas}}

Let \( (E,\bform {}{},\br ,\rho ) \) be a Courant algebroid over \( M \).
Given \( e\in \Gamma (E) \), \( f\in \cinf (M) \), we will denote \( \rho (e)f \)
simply by \( ef \), for short. Let us first prove two auxiliary lemmas.

\begin{lem}
\label{lemma:a1} The identity 
\[
T(e_{1},e_{2},\DD f)=\fourth [e_{1},e_{2}]f\]
 holds in any Courant algebroid. 
\end{lem}
\begin{proof}
Using Courant algebroid axiom 2 and Lemma\ref{lemma:ideal1}, we have 
\begin{eqnarray*}
T(e_{1},e_{2},\DD f) & = & \sixth (\bform {[e_{1},e_{2}]}{\DD f}+\bform {[\DD f,e_{1}]}{e_{2}}+\bform {[e_{2},\DD f]}{e_{1}})=\\
 & = & \sixth (\bform {[e_{1},e_{2}]}{\DD f}-\half \bform {\DD \bform {e_{1}}{\DD f}}{e_{2}}+\half \bform {\DD \bform {e_{2}}{\DD f}}{e_{1}})=\\
 & = & \sixth ([e_{1},e_{2}]f-\half e_{2}(e_{1}f)+\half e_{1}(e_{2}f))=\\
 & = & \sixth ([e_{1},e_{2}]f+\half [e_{1},e_{2}]f)=\fourth [e_{1},e_{2}]f.
\end{eqnarray*}

\end{proof}
\begin{lem}
\label{lemma:a2} Given \( e_{1},e_{2},e_{3},e_{4}\in \Gamma (E) \), let 
\begin{eqnarray*}
{\textbf {J}} & = & \bform {J(e_{1},e_{2},e_{3})}{e_{4}}-\bform {J(e_{1},e_{2},e_{4})}{e_{3}}+\bform {J(e_{1},e_{3},e_{4})}{e_{2}}-\bform {J(e_{2},e_{3},e_{4})}{e_{1}}\\
{\textbf {K}} & = & \bform {[e_{1},e_{2}]}{[e_{3},e_{4}]}-\bform {[e_{1},e_{3}]}{[e_{2},e_{4}]}+\bform {[e_{1},e_{4}]}{[e_{2},e_{3}]},
\end{eqnarray*}
 where \( J \) is the Jacobiator (cf. Def \ref{def:jac}). Then \( {\textbf {K}}+2{\textbf {J}}=0 \). 
\end{lem}
\begin{proof}
Using Courant algebroid axioms 1 and 5, we can rewrite \( {\textbf {J}} \)
as follows: 
\begin{eqnarray*}
\bform {J(e_{1},e_{2},e_{3})}{e_{4}}  = \bform {\DD T(e_{1},e_{2},e_{3})}{e_{4}}=e_{4}T(e_{1},e_{2},e_{3})=\sixth e_{4}(\bform {[e_{1},e_{2}]}{e_{3}}+c.p.)=\\
  =  \sixth (\bform {[e_{4},[e_{1},e_{2}]]+\half \DD \bform {e_{4}}{[e_{1},e_{2}]}}{e_{3}}+\bform {[e_{1},e_{2}]}{[e_{4},e_{3}]+\half \DD \bform {e_{4}}{e_{3}}})+c.p.
\end{eqnarray*}
 Expressing the other summands of \( {\textbf {J}} \) in this form and collecting
like terms in the parentheses, we find that the terms of the form \( \bform {[e_{i},e_{j}]}{\DD \bform {e_{k}}{e_{l}}} \)
cancel out, terms of the form \( \bform {[e_{i},e_{j}]}{[e_{k},e_{l}]} \) add
up to \( -4{\textbf {K}} \), those of the form \( \bform {[e_{i},[e_{j},e_{k}]]}{e_{l}} \)
add up to \( {\textbf {J}} \), and finally, terms of the form \( \bform {\DD \bform {e_{i}}{[e_{j},e_{k}]}}{e_{l}} \)
add up to \( -3{\textbf {J}} \) after we use Courant algebroid axiom 1. Thus,
\[
{\textbf {J}}=\sixth ({\textbf {J}}-3{\textbf {J}}-4{\textbf {K}}),\]
 and the statement of the lemma follows immediately. 
\end{proof}
\begin{pf*}{Proof of Lemma \ref{lemma:t1}} In the notation of the previous section, we have, using Lemma \ref{lemma:a1} and Courant algebroid axiom 2: 

\begin{eqnarray*}
&&(l_{2}l_{2}+l_{3}l_{1})(e_{1}\w e_{2}\w f)=\\
&=&l_{2}(l_{2}(e_{1}\w e_{2})\w f+l_{2}(e_{2}\w f)\w e_{1}+l_{2}(f\w e_{1})\w 
e_{2})+\\
&+&l_{3}(l_{1}e_{1}\w e_{2}\w f+e_{1}\w l_{1}e_{2}\w f+e_{1}\w e_{2}\w 
l_{1}f)=\\
&=&l_{2}([e_{1},e_{2}]\w f+\half\bform{e_{2}}{\DD f}\w e_{1}-\half\bform{\DD f}{e_{1}}\w 
e_{2})+l_{3}(e_{1}\w e_{2}\w \DD f)=\\
&=&\half\bform{[e_{1},e_{2}]}{\DD f}-\fourth\bform{e_{1}}{\DD\bform{e_{2}}{\DD f}}+
\fourth\bform{e_{2}}{\DD\bform{e_{1}}{\DD f}}-T(e_{1},e_{2},\DD f)=\\
&=&\half[e_{1},e_{2}]f-\fourth
e_{1}(e_{2}f)+\fourth e_{2}(e_{1}f)-\fourth[e_{1},e_{2}]f=0
\end{eqnarray*}
\end{pf*}

\begin{pf*}{Proof of Lemma \ref{lemma:t2}}In the notation of the previous section we have
\begin{eqnarray*}
l_{2}l_{3}(e_{1}\w e_{2}\w e_{3}\w e_{4})&=&
l_{2}(l_{3}(e_{1}\w e_{2}\w e_{3})\w e_{4}\pm(3,1)-unshuffles)=\\
&=&-l_{2}(T(e_{1},e_{2},e_{3})\w e_{4}\pm(3,1)-unshuffles)=\\
&=&\half\bform{\DD T(e_{1},e_{2},e_{3})}{e_{4}}\pm(3,1)-unshuffles=\\
&=&\half\bform{J(e_{1},e_{2},e_{3})}{e_{4}}\pm(3,1)-unshuffles=\half{\textbf{J}}.
\end{eqnarray*}
On the other hand,
\begin{eqnarray*}
l_{3}l_{2}(e_{1}\w e_{2}\w e_{3}\w e_{4})&=&l_{3}(l_{2}(e_{1}\w e_{2})\w e_{3}\w 
e_{4})\pm(2,2)-unshuffles=\\
&=&-T([e_{1},e_{2}],e_{3},e_{4})\mp(2,2)-unshuffles=\\
&=&-\sixth(\bform{[e_{1},e_{2}],e_{3}]}{e_{4}}+\bform{[e_{3},e_{4}]}{[e_{1},e_{2
}]}+
\bform{[e_{4},[e_{1},e_{2}]]}{e_{3}})\\
&\pm&\cdots=-\sixth({\textbf{J}}+2{\textbf{K}}),
\end{eqnarray*}
after collecting like terms. An application of Lemma \ref{lemma:a2} 
immediately yields \(l_{2}l_{3}=l_{3}l_{2}\).
\end{pf*}

\section{Alternative definition of Courant algebroid\label{sec:newdef}}

\begin{defn}
\label{def:quasi-algebroid2}A \emph{Courant algebroid} is a vector bundle \( E\rightarrow M \)
together with a nondegenerate symmetric bilinear form \( \bform {\cdot }{\cdot } \)
on the bundle, a bilinear operation \( \circ  \) on \( \Gamma (E) \), and
a bundle map \( \rho :E\rightarrow TM \) satisfying the following properties: 
\begin{enumerate}
\item \( e_{1}\circ (e_{2}\circ e_{3})=(e_{1}\circ e_{2})\circ e_{3}+e_{2}\circ (e_{1}\circ e_{3})\; \; \forall e_{1},e_{2},e_{3}\in \Gamma (E) \); 
\item \( \rho (e_{1}\circ e_{2})=[\rho (e_{1}),\rho (e_{2})]\; \; \forall e_{1},e_{2}\in \Gamma (E) \); 
\item \( e_{1}\circ fe_{2}=f(e_{1}\circ e_{2})+(\rho (e_{1})\cdot f)e_{2}\; \; \forall e_{1},e_{2}\in \Gamma (E),\; f\in \cinf (M) \); 
\item \( e\circ e=\half \DD \bform {e}{e}\; \; \forall e\in \Gamma (E) \); 
\item \( \rho (e)\cdot \bform {h_{1}}{h_{2}}=\bform {e\circ h_{1}}{h_{2}}+\bform {h_{1}}{e\circ h_{2}}\; \; \forall e,h_{1},h_{2}\in \Gamma (E) \), 
\end{enumerate}
where \( \DD :\cinf (M)\rightarrow \Gamma (E) \) is given by (\ref{eq:D}).
\end{defn}
Notice that all the anomalies of Definition \ref{def:quasi-algebroid} are absent
in this one, but the skew-symmetric bracket \( \br  \) is replaced by a not
necessarily skew-symmetric operation \( \circ  \).\footnote{
This was first proposed in \cite{LWX1}, but the properties of this operation
were then an open question; this definition of a Courant algebroid was also
used in a note of P. \v{S}evera \cite{Sev}, without a proof of its equivalence
to the original one.
} Property 1 above is to be interpreted as the ``Jacobi identity'' for \( \circ  \)
in the sense that \( e\circ \cdot  \) is a derivation of \( \circ  \) for
any \( e\in \Gamma (E) \);\footnote{
It looks more like a Leibniz rule; in fact, a vector space with a bilinear operation
satisfying this property was called a Leibniz algebra by Loday \cite{Loday1},
and a Loday algebra by Kosmann-Schwarzbach \cite{KS5}.
} if \( \circ  \) is skew-symmetric, this is equivalent to the usual Jacobi
identity. On the other hand, Property 4 is equivalent to saying, by a polarization
identity, that \( \circ  \) is skew-symmetric ``up to a coboundary'', i.e.
the symmetric part of it is \( \DD  \) of something. More precisely, we have
\begin{equation}
\label{eqn:symplusskewsym}
e_{1}\circ e_{2}=[e_{1},e_{2}]+\half \DD \bform {e_{1}}{e_{2}}
\end{equation}
for all \( e_{1},e_{2}\in \Gamma (E) \), where 
\begin{equation}
\label{eqn:skewsym}
[e_{1},e_{2}]=\half (e_{1}\circ e_{2}-e_{2}\circ e_{1})
\end{equation}
is the skew-symmetrization of \( \circ . \) We shall now prove that the new
Definition \ref{def:quasi-algebroid2} is equivalent to the old Definition \ref{def:quasi-algebroid}.
We need a couple of lemmas first. The first one is the non-skew-symmetric version
of Lemma \ref{lemma:ideal1}. 

\begin{lem}
\label{lemma:ideal2} If \( (E,\bform {\cdot }{\cdot },\circ ,\rho ) \) satisfies
properties 2-5 of Definition \ref{def:quasi-algebroid2}, then \( \forall e\in \Gamma (E)\), \(f\in \cinf (M) \)
one has 
\[
\begin{array}{ccl}
e\circ \DD f & = & \DD \bform {e}{\DD f}\\
\DD f\circ e & = & 0
\end{array}\]

\end{lem}
\begin{proof}
Let \( h\in \Gamma (E) \) be arbitrary. Then, using Properties 2 and 5 we have
\[
\begin{array}{ccl}
\rho (e)(\rho (h)f) & = & \rho (e)\bform {\DD f}{h}=\bform {e\circ \DD f}{h}+\bform {\DD f}{e\circ h}=\\
 & = & \bform {e\circ \DD f}{h}+\rho (e\circ h)f=\\
 & = & \bform {e\circ \DD f}{h}+\rho (e)(\rho (h)f)-\rho (h)(\rho (e)f)
\end{array}\]
 Hence, 
\[
\bform {e\circ \DD f}{h}=\rho (h)(\rho (e)f)=\bform {h}{\DD \bform {e}{\DD f}}\]
 The first statement follows by the nondegeneracy of \( \bform {\cdot }{\cdot }. \)
On the other hand, by (\ref{eqn:symplusskewsym}), 
\[
\DD f\circ e=\DD f\circ e+e\circ \DD f-e\circ \DD f=\DD \bform {e}{\DD f}-\DD \bform {e}{\DD f}=0\]

\end{proof}
\begin{rem}
\label{rem:(lem:ideal1=3Dlem:ideal2)}Observe that the statement of Lemma \ref{lemma:ideal2}
is equivalent to the statement of Lemma \ref{lemma:ideal1} for the skew-symmetrization
(\ref{eqn:skewsym}), in view of (\ref{eqn:symplusskewsym}).
\end{rem}
\begin{lem}
\label{lemma:jacisskewsym} If \( (E,\bform {\cdot }{\cdot },\circ ,\rho ) \)
satisfies properties 2-5 of definition \ref{def:quasi-algebroid2}, then the
expression
\[
K(e_{1},e_{2},e_{3})=(e_{1}\circ e_{2})\circ e_{3}+e_{2}\circ (e_{1}\circ e_{3})-e_{1}\circ (e_{2}\circ e_{3})\]
 is completely skew-symmetric in \( e_{1},e_{2},e_{3} \).
\end{lem}
\begin{proof}
We have to show that \( K \) vanishes if any two of the entries coincide. But
\[
K(e_{1},e_{1},e_{3})=(e_{1}\circ e_{1})\circ e_{3}+e_{1}\circ (e_{1}\circ e_{3})-e_{1}\circ (e_{1}\circ e_{3})=\half \DD \bform {e_{1}}{e_{1}}\circ e_{3}=0\]
 by property 4 and Lemma \ref{lemma:ideal2}. On the other hand, 
\[
\begin{array}{ccl}
K(e_{1},e_{2},e_{2}) & = & (e_{1}\circ e_{2})\circ e_{2}+e_{2}\circ (e_{1}\circ e_{2})-e_{1}\circ (e_{2}\circ e_{2})=\\
 & = & \DD (\bform {e_{1}\circ e_{2}}{e_{2}}-\bform {e_{1}}{e_{2}\circ e_{2}})=\\
 & = & \DD (\bform {e_{1}\circ e_{2}}{e_{2}}+\bform {e_{2}\circ e_{1}}{e_{2}}-\bform {e_{2}}{\DD \bform {e_{1}}{e_{2}}})=\\
 & = & \DD (\bform {\DD \bform {e_{1}}{e_{2}}}{e_{2}}-\bform {e_{2}}{\DD \bform {e_{1}}{e_{2}}})=0,
\end{array}\]
 where we have used properties 4 and 5 and Lemma \ref{lemma:ideal2}. And finally,
\[
\begin{array}{l}
K(e_{1},e_{2},e_{1})=(e_{1}\circ e_{2})\circ e_{1}+e_{2}\circ (e_{1}\circ e_{1})-e_{1}\circ (e_{2}\circ e_{1})=\\
=(e_{1}\circ e_{2})\circ e_{1}+(e_{2}\circ e_{1})\circ e_{1}+e_{2}\circ (e_{1}\circ e_{1})-(e_{2}\circ e_{1})\circ e_{1}-e_{1}\circ (e_{2}\circ e_{1})=\\
=\DD \bform {e_{1}}{e_{2}}\circ e_{1}-\DD \bform {e_{2}\circ e_{1}}{e_{1}}+\DD \bform {e_{2}}{e_{1}\circ e_{1}}=\\
=-\DD (\bform {e_{2}\circ e_{1}}{e_{1}}-\bform {e_{2}}{e_{1}\circ e_{1}})=0,
\end{array}\]
 just as in the previous calculation; we have used again properties 4 and 5
and Lemma \ref{lemma:ideal2}.
\end{proof}
We are now ready to prove the equivalence of the two definitions of Courant
algebroid.

\begin{prop}
\label{prop:1=3D2}Let \( (E,\bform {\cdot }{\cdot },\br ,\rho ) \) be a Courant
algebroid in the sense of Definition \ref{def:quasi-algebroid}. Let the operation
\( \circ  \) be given by (\ref{eqn:symplusskewsym}). Then \( (E,\bform {\cdot }{\cdot },\circ ,\rho ) \)
is a Courant algebroid in the sense of Definition \ref{def:quasi-algebroid2}.

Conversely, let \( (E,\bform {\cdot }{\cdot },\circ ,\rho ) \) be a Courant
algebroid in the sense of Definition \ref{def:quasi-algebroid2}. Let \( \br  \)
be the skew-symmetrization of \( \circ  \), as in (\ref{eqn:skewsym}). Then
\( (E,\bform {\cdot }{\cdot },\br ,\rho ) \) is a Courant algebroid in the
sense of Definition \ref{def:quasi-algebroid}.
\end{prop}
\begin{proof}
It is not hard to show the equivalence of all of the properties, except for
the Jacobi identity which will take a bit more work. So we shall first show
the equivalence of properties 2-5 for both definitions, and then prove the equivalence
of the two versions of Jacobi using Lemmas \ref{lemma:ideal2} and \ref{lemma:jacisskewsym}. 

Now, in view of (\ref{eqn:symplusskewsym}) it is obvious that Property 5 is
equivalent for the two definitions. Property 2 in the new definition implies
immediately that \( \rho (e\circ e)=0 \) for all \( e\in \Gamma (E) \), hence
\[
[\rho (e_{1}),\rho (e_{2})]=\rho (e_{1}\circ e_{2})=\rho ([e_{1},e_{2}]),\]
 and we have the old Property 2. Moreover, by the new Property 4, 
\[
0=\rho (e_{1}\circ e_{2}+e_{2}\circ e_{1})=\half \DD \bform {e_{1}}{e_{2}}\; \; \forall e_{1},e_{2}\in \Gamma (E),\]
 hence we have the old Property 4 (\( \rho \smalcirc \DD =0 \)) by the nondegeneracy
of \( \bform {\cdot }{\cdot } \). Conversely, if we start with the old definition,
the new Property 4 is immediate by (\ref{eqn:symplusskewsym}), while the old
Properties 2 and 4 combine to give 
\[
\rho (e_{1}\circ e_{2})=\rho ([e_{1},e_{2}]+\half \DD \bform {e_{1}}{e_{2}})=\rho ([e_{1},e_{2}])=[\rho (e_{1}),\rho (e_{2})],\]
 the new Property 2. 

As for the Leibniz rule, one has 
\[
\begin{array}{l}
e_{1}\circ fe_{2}=[e_{1},fe_{2}]+\half \DD \bform {e_{1}}{fe_{2}}=\\
=[e_{1},fe_{2}]+\half f\DD \bform {e_{1}}{e_{2}}+\half \bform {e_{1}}{e_{2}}\DD f
\end{array}\]
 for all \( e_{1},e_{2}\in \Gamma (E) \), \( f\in \cinf (M) \), hence it follows
immediately that the new and old Properties 3 are equivalent.

Now for the Jacobi identity. In view of (\ref{eqn:symplusskewsym}), it is clear
that one has 
\[
K(e_{1},e_{2},e_{3})=J(e_{1},e_{2},e_{3})+R(e_{1},e_{2},e_{3}),\]
 where \( K \) is as in Lemma \ref{lemma:jacisskewsym}, \( J \) is the Jacobiator
(Def. \ref{def:jac}), and 
\[
\begin{array}{ccl}
R(e_{1},e_{2},e_{3}) & = & \half ([\DD \bform {e_{1}}{e_{2}},e_{3}]+[e_{2},\DD \bform {e_{1}}{e_{3}}]-[e_{1},\DD \bform {e_{2}}{e_{3}}])+\\
 & + & \half \DD (\bform {e_{1}\circ e_{2}}{e_{3}}+\bform {e_{2}}{e_{1}\circ e_{3}}-\bform {e_{1}}{e_{2}\circ e_{3}}).
\end{array}\]
 To show the equivalence of the old and new Properties 1, we only need to show
that \( R(e_{1},e_{2},e_{3})=-\DD T(e_{1},e_{2},e_{3}) \), where \( T \) is
as in (\ref{eq:T0}). But, by Lemma \ref{lemma:ideal1} and Remark \ref{rem:(lem:ideal1=3Dlem:ideal2)},
\[
\begin{array}{l}
\half ([\DD \bform {e_{1}}{e_{2}},e_{3}]+[e_{2},\DD \bform {e_{1}}{e_{3}}]-[e_{1},\DD \bform {e_{2}}{e_{3}}])=\\
=-\fourth \DD (\bform {\DD \bform {e_{1}}{e_{2}}}{e_{3}}-\bform {e_{2}}{\DD \bform {e_{1}}{e_{3}}}+\bform {e_{1}}{\DD \bform {e_{2}}{e_{3}}}),
\end{array}\]
 whereas 
\[
\begin{array}{l}
\half \DD (\bform {e_{1}\circ e_{2}}{e_{3}}+\bform {e_{2}}{e_{1}\circ e_{3}}-\bform {e_{1}}{e_{2}\circ e_{3}})=\\
=\half \DD (\bform {[e_{1},e_{2}]}{e_{3}}+\bform {e_{2}}{[e_{1},e_{3}]}-\bform {e_{1}}{[e_{2},e_{3}]})+\\
+\fourth \DD (\bform {\DD \bform {e_{1}}{e_{2}}}{e_{3}}+\bform {e_{2}}{\DD \bform {e_{1}}{e_{3}}}-\bform {e_{1}}{\DD \bform {e_{2}}{e_{3}}})
\end{array}\]
 by (\ref{eqn:symplusskewsym}). Therefore, 
\begin{equation}
\label{eqn:R}
\begin{array}{ccl}
R(e_{1},e_{2},e_{3}) & = & \half \DD (\bform {[e_{1},e_{2}]}{e_{3}}-\bform {[e_{3},e_{1}]}{e_{2}}-\bform {[e_{2},e_{3}]}{e_{1}})+\\
 & + & \half \DD (\bform {e_{2}}{\DD \bform {e_{1}}{e_{3}}}-\bform {e_{1}}{\DD \bform {e_{2}}{e_{3}}})
\end{array}
\end{equation}
 However, since both \( J \) and \( K \) are completely antisymmetric (Lemma
\ref{lemma:jacisskewsym}), so is \( R \); therefore, \( R \) is equal to
its skew-symmetrization. But it is obvious that the skew-symmetrization of the
first term on the right hand side of (\ref{eqn:R}) is \( -\DD T(e_{1},e_{2},e_{3}) \),
whereas the skew-symmetrization of the second term is easily seen to be zero.
Hence \\ \( R(e_{1},e_{2},e_{3})=-\DD T(e_{1},e_{2},e_{3}) \), and we are done.
\end{proof}
\begin{rem}
Notice that the notion of a Dirac subbundle remains unchanged when we switch
to the new definition of a Courant algebroid, thanks to (\ref{eqn:symplusskewsym}),
and that the restrictions of the two operations to any Dirac subbundle are identical.
\end{rem}
\begin{example}
\label{eg:double}Let \( (A,A^{*}) \) be a pair of Lie algebroids in duality,
with anchors \( a \) and \( a_{*} \), respectively. Then on \( E=A\oplus A^{*} \)
we define 

\begin{equation}
\label{eqn:double2}
\begin{array}{rcl}
\bform {X+\xi }{Y+\eta } & = & \xi (Y)+\eta (X)\\
(X+\xi )\circ (Y+\eta ) & = & ([X,Y]_{A}+L_{\xi }^{A^{*}}Y-i_{\eta }d_{A^{*}}X)+\\
 & + & ([\xi ,\eta ]_{A^{*}}+L^{A}_{X}\eta -i_{Y}d_{A}\xi )\\
\rho (X+\xi ) & = & a(X)+a_{*}(\xi )
\end{array}
\end{equation}

\end{example}
If \( (A,A^{*}) \) is a Lie bialgebroid, then by Theorem \ref{thm:bialg->courant}
and Proposition \ref{prop:1=3D2} \( (E,\bform {\cdot }{\cdot },\circ ,\rho ) \)
is a Courant algebroid in the sense of Definition \ref{def:quasi-algebroid2}.

\begin{example}
As a special case of Example \ref{eg:double}, consider \( TM \) with its standard
Lie algebroid structure and \( T^{*}M \) with the zero anchor and bracket.
Then on sections of \( TM\oplus T^{*}M \) the operation \( \circ  \) reduces
to 
\begin{equation}
\label{eqn:standardCA}
(X+\xi )\circ (Y+\eta )=[X,Y]+L_{X}\eta -i_{Y}d\xi 
\end{equation}
 whose skew-symmetrization was originally considered by Courant in his study
of Dirac manifolds \cite{Cou}. This Courant algebroid was also considered by
P. \v{S}evera \cite{Sev} as the natural geometric framework for two-dimensional
variational problems.
\end{example}

\chapter{The double of a Lie bialgebroid as a homological vector field on an even symplectic
supermanifold\label{chapter:supermanifolds}}

\newcommand{\der}[2]{\frac{\partial #1 }{\partial #2 }}

We shall now present an alternative construction of the double of a Lie bialgebroid.
It is based on an interpretation of a Lie algebroid as an odd self-commuting
vector field on a supermanifold, which we then lift as a hamiltonian on its
cotangent bundle. Adding the two hamiltonians coming from the dual Lie algebroids,
we get a third one which Poisson-commutes with itself if and only if the compatibility
condition of a Lie bialgebroid is satisfied; the corresponding hamiltonian vector
field is interpreted as the Drinfeld double. As an application, we show that
the Weil differential and the classical BRST differential arise in this way.
The Courant algebroid of Example \ref{eg:double} is recovered via the derived
bracket construction; this allows us to give a simple proof of the doubling
theorem of Liu, Weinstein and Xu \cite{LWX1}. Finally, we consider quasi-Lie
bialgebroids which one gets by adding cubic terms to the hamiltonian and show
that exact Courant algebroids, recently classified by \v{S}evera \cite{Sev},
arise in this way. 

The starting point for us is a picture of Lie bialgebras due to Lecomte, Roger and 
Kosmann-Schwarzbach.

\section{An alternative  picture of Lie bialgebras\label{sec:KSpicture}}

There is an elegant way to express the structure
relations of a Lie bialgebra by embedding it into a larger space endowed with
a canonical Poisson superalgebra structure \cite{LecRog} \cite{KS3}. By viewing this construction from
an appropriate angle we shall be able to generalize it to Lie bialgebroids,
obtain a new notion of the Drinfeld double and recover the old one.

Consider a Lie bialgebra \( (\g ,\mu ,\gamma ) \) (see Definition \ref{def:LieBialgebra});
view the bracket \( \mu  \) and the cobracket \( \gamma  \) as elements \( \mu \in \wedge ^{2}\g ^{*}\otimes \g  \)
and \( \gamma \in \g ^{*}\otimes \wedge ^{2}\g  \). 

The basic idea is to embed \( \mu  \) and
\( \gamma  \) into the full exterior algebra \( \wedge (\g \oplus \g ^{*})=(\wedge \g )\otimes (\wedge \g ^{*}) \),
and take advantage of a natural Poisson superalgebra structure on this space,
defined as follows. The commutative superalgebra structure is given by the exterior
multiplication, whereas the (even) Poisson bracket 
\[
\{\cdot ,\cdot \}:\wedge ^{k}(\g \oplus \g ^{*})\times \wedge ^{l}(\g \oplus \g ^{*})\longrightarrow \wedge ^{k+l-2}(\g \oplus \g ^{*})\]
(called \emph{the big bracket} in \cite{LecRog} and \cite{KS3}, although it goes back to \cite{KoSt}) is uniquely determined
by the following properties: 

\begin{itemize}
\item For any \( a,b\in \g  \), \( \{a,b\}=0 \);
\item For any \( \xi ,\eta \in \g ^{*} \), \( \{\xi ,\eta \}=0 \);
\item For any \( a\in \g  \), \( \xi \in \g ^{*} \), \( \{\xi ,a\}=\xi (a) \) ;
\item \( \{\cdot ,\cdot \} \) is skew-symmetric, i.e. for any \( e_{1}\in \wedge ^{k}(\g \oplus \g ^{*}) \),
\( e_{2}\in \wedge ^{l}(\g \oplus \g ^{*}) \), 
\[
\{e_{1},e_{2}\}=-(-1)^{kl}\{e_{2},e_{1}\}\]
 
\item For every \( e\in \wedge ^{k}(\g \oplus \g ^{*}) \), \( \{e,\cdot \} \) is
a derivation of the exterior algebra \( \wedge (\g \oplus \g ^{*}) \) of degree
\( k-2 \). 
\end{itemize}
In other words. \( \{\cdot ,\cdot \} \) is the unique extension of the canonical
symmetric bilinear form \( \bform {\cdot }{\cdot } \) on \( \g \oplus \g ^{*} \)
(\ref{eqn:CanBform}) to an even Poisson superalgebra structure on \( \wedge (\g \oplus \g ^{*}) \):
it is easy to show that the super Jacobi identity 
\[
\{e_{1},\{e_{2},e_{3}\}\}=\{\{e_{1},e_{2}\},e_{3}\}+(-1)^{kl}\{e_{2},\{e_{1},e_{3}\}\}\]
 holds for all \( e_{1}\in \wedge ^{k}(\g \oplus \g ^{*}) \), \( e_{2}\in \wedge ^{l}(\g \oplus \g ^{*}) \),
\( e_{3}\in \wedge ^{m}(\g \oplus \g ^{*}) \). 

Using this operation, it can be shown without difficulty that \( (\g ,\mu ,\gamma ) \)
is a Lie bialgebra if and only if 
\begin{equation}
\label{eqn:LieBialgebra}
\{\mu ,\mu \}=\{\gamma ,\gamma \}=\{\mu ,\gamma \}=0.
\end{equation}
Here \( \{\mu ,\mu \}=0 \) (resp. \( \{\gamma ,\gamma \}=0 \)) is equivalent
to the Jacobi identity for \( \br  \) (resp. \( \br _{*} \), while \( \{\mu ,\gamma \}=0 \)
is equivalent to the cocycle condition. The brackets \( \br  \) and \( \br _{*} \)
can be recovered by the formulas 
\begin{equation}
\label{eqn:brackets}
\begin{array}{ccc}
{[a,b]} & = & \{\{\mu ,a\},b\}\\
{[\xi ,\eta ]_{*}} & = & \{\{\gamma ,\xi \},\eta \}
\end{array}
\end{equation}
 where \( a,b\in \g  \), \( \xi ,\eta \in \g ^{*} \). Furthermore, if we set
\( \theta =\mu +\gamma  \), then (\ref{eqn:LieBialgebra}) is equivalent to
\begin{equation}
\label{eqn:master1}
\{\theta ,\theta \}=0,
\end{equation}
 and if \( e_{1},e_{2}\in \g \oplus \g ^{*}\subset \wedge (\g \oplus \g ^{*}) \),
\begin{equation}
\label{eqn:double1}
[e_{1},e_{2}]=\{\{\theta ,e_{1}\},e_{2}\}
\end{equation}
 is precisely the Drinfeld double bracket on \( \g \oplus \g ^{*} \) (\ref{eqn:doubleofbialgebra}):
it is skew-symmetric, since the symmetric part is proportional to \( \{\theta ,\{e_{1},e_{2}\}\} \)
which is zero because \( \{e_{1},e_{2}\}\in \R  \); it is easy to see from
(\ref{eqn:brackets}) that \( \g  \) and \( \g ^{*} \)are subalgebras, and
that the canonical inner product \( \bform {\cdot }{\cdot } \) is ad-invariant;
therefore, it must coincide with (\ref{eqn:doubleofbialgebra}). The Jacobi
identity is a consequence of (\ref{eqn:master1}). 

From the point of view of \cite{KS3}, the main advantage of this approach is
that it affords an elegant treatment of Drinfeld's \emph{quasi-Lie bialgebras},
a generalization of Lie bialgebras in which the Jacobi identity for one of the
brackets is satisfied only up to a coboundary. This amounts to adding a \( \phi \in \wedge ^{3}\g  \)
or \( \psi \in \wedge ^{3}\g ^{*} \) (or both) to \( \theta  \) so that \( \theta  \)
still satisfies (\ref{eqn:master1}). Note that the double (\ref{eqn:double1})
is still a Lie algebra, even though \( \g  \) or \( \g ^{*} \) are not.

For our purposes, however, the chief value of this approach is
that it generalizes to Lie bialgebroids, if interpreted correctly; this will
be our next order of business.

\section{Lie algebroids revisited\label{sec:algebroidsagain}}

As a naive attempt to generalize the above construction to Lie algebroids,
we might try, given a vector bundle \( A\rightarrow M \), to consider the exterior
algebra \( \Gamma (\bigwedge (A\oplus A^{*})) \) and build the ``big bracket''
\( \{\cdot ,\cdot \} \) on this space from the canonical inner product \( \bform {\cdot }{\cdot } \)
on \( A\oplus A^{*} \) (see Example \ref{eg:double}), as above. However, one
quickly realizes that this is not enough to encode a Lie algebroid structure
on \( A \) or \( A^{*} \). First, the anchor: \( a:A\rightarrow TM \) can
be viewed as a section of \( A^{*}\otimes TM \), so there is no room for it
in \( \Gamma (\bigwedge (A\oplus A^{*})) \); on the other hand, the Lie algebroid
bracket (say, on \( \Gamma (A) \)) can no longer be viewed as a section of
\( \bigwedge ^{2}A^{*}\otimes A \), since it is not linear over \( \cinf (M) \)
unless \( a \) is trivial. Furthermore, the structural identities of a Lie
algebroid (e.g. the Jacobi identity) are not algebraic but differential equations,
so they cannot be encoded by \( \{\cdot ,\cdot \} \) which is \( \cinf (M) \)-linear.
Thus it is clear that we need a bigger space with an even Poisson superalgebra
structure in which \( (\Gamma (\bigwedge (A\oplus A^{*})),\{\cdot ,\cdot \}) \)
can be embedded. In order to find this space, we must shift our point of view
from an algebraic to a geometric one and use the language of supermanifolds. 

Recall from Chapter \ref{chapter:courant} that a Lie algebroid structure on
a vector bundle \( A \) over \( M \) is equivalent to a derivation \( d_{A} \)
of the exterior algebra \( \Gamma (\bigwedge A^{*}) \) of degree one and square
zero. Just as above, we will view \( \Gamma (\bigwedge A^{*}) \) as the algebra
of functions on the supermanifold \( \Pi A \), where \( \Pi  \) here denotes
the change of parity functor applied to each fibre (see Appendix). If \( \{x^{i}\}_{i=1,\ldots ,\dim M} \)
is a coordinate chart on \( U\subset M \), and \( \{e^{a}\}_{a=1,\ldots ,\textrm{rk}A} \)
is a local basis of sections of \( A^{*} \) over \( U \) (dual to a basis
\( \{e_{a}\} \) of sections of \( A \)), denote by \( \xi ^{a} \) the corresponding
generators of the Grassman algebra \( \Gamma (U,\bigwedge A^{*}) \); then \( \{(x^{i},\xi ^{a})\} \)
give a coordinate chart on \( \Pi A \) with the transformation law inherited
from the vector bundle \( A^{*} \). The derivation \( d_{A} \) can then be
viewed as an (odd) vector field on \( \Pi A \), satisfying 
\begin{equation}
\label{eq:homvf}
[d_{A},d_{A}]=2d^{2}_{A}=0,
\end{equation}
 where the bracket denotes the (super)commutator. Such vector fields are called
\emph{homological} for an obvious reason. This motivates the following 

\begin{defn}
A \emph{Lie algebroid structure} on a vector bundle \( A\rightarrow M \) is
the supermanifold \( \Pi A \) together with a homological vector field \( d_{A} \)
of degree 1.
\end{defn}
\begin{rem}
This interpretation of Lie algebroids was proposed by Kontsevich \cite{Ko1}
and Vaintrob \cite{Vain}. It is important that \( d_{A} \) be of degree 1
with respect to the natural \( \Z  \)-grading on functions on \( \Pi A \),
rather than merely be odd, in order to define a Lie algebroid structure on \( A \).
Arbitrary odd homological vector fields on supermanifolds lead to strongly homotopy
Lie algebras \cite{AKSZ}. 
\end{rem}
In local coordinates, we have, according to the Cartan formula (\ref{eqn:CartanFla}),
\begin{equation}
\label{eqn:dA}
d_{A}=\xi ^{a}A^{i}_{a}(x)\der{}{x^{i}}-\half C^{c}_{ab}(x)\xi ^{a}\xi ^{b}\der{}{\xi ^{c}}
\end{equation}
 where 
\[
\begin{array}{rcl}
a(e_{a}) & = & A^{i}_{a}(x)\der{}{x^{i}}\\
{[e_{a},e_{b}]_{A}} & = & C_{ab}^{c}(x)e_{c}
\end{array}\]
 are the local expressions for the anchor and the bracket on the Lie algebroid
\( A \). Similarly, a Lie algebroid structure on the dual bundle \( A^{*} \)
is equivalent to a homological vector field \( d_{A^{*}} \) on the supermanifold
\( \Pi A^{*} \) given in local coordinates \( (x^{i},\theta _{a}) \) by 
\begin{equation}
\label{eqn:dA*}
d_{A^{*}}=\theta _{a}\bar{A}^{ai}(x)\der{}{x^{i}}-\half \bar{C}^{ab}_{c}(x)\theta _{a}\theta _{b}\der{}{\theta _{c}}
\end{equation}
 where 
\[
\begin{array}{rcl}
a_{*}(e^{a}) & = & \bar{A}^{ai}(x)\der{}{x^{i}}\\
{[e^{a},e^{b}]_{A^{*}}} & = & \bar{C}^{ab}_{c}(x)e^{c}
\end{array}\]
 are the local expressions for the anchor and the bracket on \( A^{*} \).

\section{The cotangent bundle\label{sec:cotbundle}}

Once in a ``supermathematical'' frame of mind, one immediately realizes that
the exterior algebra \( \wedge (\g \oplus \g ^{*}) \) is to be interpreted as the algebra of functions on the (purely
odd) superspace \( \Pi (\g \oplus \g ^{*}) \). The crucial observation, however,
is that \( \Pi (\g \oplus \g ^{*}) \) is naturally isomorphic to the cotangent
bundle \( T^{*}\Pi \g  \), while the big bracket \( \{\cdot ,\cdot \} \) is
nothing but the canonical symplectic Poisson bracket on \( T^{*}\Pi \g  \).
Indeed, if \( \{e_{a}\}_{a=1,\ldots ,\dim \g } \) is a basis of \( \g  \),
\( \{e^{a}\} \) the dual basis, denote by \( \{\theta _{a}\} \) the corresponding
generators of the Grassman algebra \( \wedge \g =\cinf (\Pi \g ^{*}) \), and
\( \{\xi ^{a}\} \) the corresponding generators of \( \wedge \g ^{*}=\cinf (\Pi \g ) \).
Then the \( \theta _{a} \) can be viewed as the momenta conjugate to the Grassman
coordinates \( \xi ^{a} \) on \( \Pi \g  \), and the defining relations of
the big bracket (Section \ref{sec:KSpicture}) can be rewritten as 
\[
\{\xi ^{a},\xi ^{b}\}=0;\; \{\theta _{a},\theta _{b}\}=0;\; \{\xi ^{a},\theta _{b}\}=\delta ^{a}_{b}\]
 which one immediately recognizes as the canonical Poisson brackets between
coordinates and momenta on \( T^{*}\Pi \g  \). This bracket is nondegenerate
and even, in the sense that the bracket of two Grassman polynomials of parity
\( \epsilon _{1},\epsilon _{2}\in \Z _{2} \) is of parity \( \epsilon _{1}+\epsilon _{2} \);
the corresponding symplectic 2-form is 
\[
\omega =d\theta _{a}d\xi ^{a}\]
 Thus, \( T^{*}\Pi \g  \) is an even symplectic supermanifold (see Appendix).
This is completely analogous to the canonical symplectic structure on \( T^{*}V\simeq V\oplus V^{*} \),
where \( V \) is a vector space, except now the Poisson bracket on linear functions
(which are odd) is symmetric rather than skew-symmetric; in fact, the matrix
of \( \omega  \) coincides with the matrix of \( \bform {\cdot }{\cdot } \)
(\ref{eqn:CanBform}) in the basis \( \{e_{a},e^{b}\} \) of \( \g \oplus \g ^{*} \).
The advantage of this point of view is that it generalizes immediately to vector
bundles.

Just as in the purely odd or even case, given any supermanifold \( Q \), its
cotangent bundle \( T^{*}Q \) is an even symplectic supermanifold. If \( \{x^{\alpha }\} \)
is a coordinate chart for \( Q \), the corresponding Darboux chart for \( T^{*}Q \)
is \( \{x^{\alpha },x_{\alpha }^{*}\} \), where \( x_{\alpha }^{*} \) is of
the same parity as \( x^{\alpha } \) and the canonical Poisson brackets are
given by 
\begin{equation}
\label{eqn:CanBrackets}
\begin{array}{ccc}
\{x^{\alpha },x^{\beta }\}=0;\;  & \{x_{\alpha }^{*},x_{\beta }^{*}\}=0;\;  & \{x_{\alpha }^{*},x^{\beta }\}=\delta _{\alpha }^{\beta }
\end{array}
\end{equation}
 Any vector field \( v \) on a \( Q \) gives rise to a fibrewise-linear function
\( h_{v} \) on the cotangent bundle \( T^{*}Q \) in an obvious manner: in
local coordinates, if \( v=v^{\alpha }(x)\der{}{x^{\alpha }} \), then \( h_{v}=v^{\alpha }(x)x_{\alpha }^{*} \)
in the corresponding Darboux coordinates on \( T^{*}Q \). This is well defined
since the momenta \( x_{\alpha }^{*} \) transform in the same way as the derivations
\( \der{}{x^{\alpha }} \) under changes of coordinates on \( Q \). This ``hamiltonian
lift'' has the following properties:

\begin{lem}
\label{lemma:cotbundle}Let \( v,w \) be vector fields on \( Q \), \( f\in \cinf (Q) \),
and let \( \pi :T^{*}Q\rightarrow Q \) denote the canonical projection. Then
\begin{enumerate}
\item \( \{h_{v},\pi ^{*}f\}=\pi ^{*}(vf) \)
\item \( \{h_{v},h_{w}\}=h_{[v,w]} \)
\end{enumerate}
\end{lem}
\begin{proof}
This is best done by direct computation in local coordinates, just as for ordinary
manifolds. (1) is obvious by (\ref{eqn:CanBrackets}) and the definition of
\( h_{v} \), whereas for (2) we can easily deduce from (\ref{eqn:CanBrackets})
that 
\[
\begin{array}{l}
\{h_{v},h_{w}\}=\{v^{\alpha }(x)x^{*}_{\alpha },w^{\beta }(x)x^{*}_{\beta }\}=\\
=(v^{\alpha }\der{w^{\beta }}{x^{\alpha }}-(-1)^{\tilde{v}\tilde{w}}w^{\alpha }\der{v^{\beta }}{x^{\alpha }})x^{*}_{\beta }=h_{[v,w]},
\end{array}\]
 where \( \tilde{v} \) denotes the parity of the vector field \( v \). 
\end{proof}
Now let \( Q=\Pi A \), and let \( \mu =h_{d_{A}}\in \cinf (T^{*}\Pi A) \).
Then by Lemma \ref{lemma:cotbundle} and (\ref{eq:homvf}) we immediately get
\begin{equation}
\label{eq:master}
\{\mu ,\mu \}=0
\end{equation}
 The formula (\ref{eqn:dA}) leads to the following local expression for \( \mu  \):
\begin{equation}
\label{eqn:mu}
\mu =\xi ^{a}A^{i}_{a}(x)x^{*}_{i}-\half C^{c}_{ab}(x)\xi ^{a}\xi ^{b}\xi ^{*}_{c}
\end{equation}
 Thus, a Lie algebroid structure on \( A \) gives rise to an odd linear function
\( \mu  \) on \( T^{*}\Pi A \) satisfying \( \{\mu ,\mu \}=0 \), but how
do we characterize those that come from Lie algebroids? Two remarks are in order.

\begin{rem}
\label{rem:notsimple}Unless the bundle \( A\rightarrow M \) is trivial, the
supermanifold \( T^{*}\Pi A \) is \emph{not} of the form \( \Pi V \) in any
natural way, where \( V \) is some vector bundle over \( T^{*}M \), the support
(``even part'') of \( T^{*}\Pi A \). That is to say, there is no canonical
projection from \( T^{*}\Pi A \) to \( T^{*}M \) and a natural \( \Z _{+} \)-grading
on \( \cinf (T^{*}\Pi A) \) inducing the \( \Z _{2} \)-grading and respecting
the projection. The reason is that under the natural fiberwise linear coordinate
changes 
\[
\begin{array}{rcl}
x^{i} & = & x^{i}(x')\\
\xi ^{a} & = & T^{a}_{a'}(x')\xi ^{a'}
\end{array}\]
 on \( \Pi A \) the momenta transform by 
\[
\begin{array}{rcl}
x^{*}_{i} & = & \der{x^{i'}}{x^{i}}(x(x'))x^{*}_{i'}+\der{T_{b}^{a'}}{x^{i}}(x(x'))T^{b}_{b'}(x')\xi ^{b'}\xi ^{*}_{a'}\\
\xi ^{*}_{a} & = & T^{a'}_{a}(x(x'))\xi ^{*}_{a'}
\end{array}\]
 so the total degree in the odd variables \( (\xi ^{a},\xi ^{*}_{b}) \) is
not preserved because of the second term in the transformation law for \( x_{i}^{*} \).
Because of this fact, our constructions cannot be recast in the ``classical''
framework of manifolds and vector bundles - one cannot get around using supermanifolds. 
\begin{rem}
\label{rem:Gradings}Nevertheless, there are three \( \Z _{+} \)-gradings on
fiberwise-polynomial functions on \( T^{*}\Pi A \) which are preserved under
the natural transformations above. The first one, which we denote by \( \epsilon  \),
exists by virtue of the fact that \( T^{*}\Pi A \) is a vector bundle over
\( \Pi A \) - it is simply the fiberwise degree, i.e the total degree of a
polynomial in the momenta \( (x_{i}^{*},\xi _{a}^{*}) \); the second, \( \delta  \),
is the total degree of a polynomial in \( (x^{*}_{i},\xi ^{a}) \). These gradings
are not compatible with the \( \Z _{2} \)-grading (parity), since the even
variables \( x^{*}_{i} \) have \( \epsilon (x^{*}_{i})=\delta (x^{*}_{i})=1 \),
but their sum \( \kappa  \) is. The total grading \( \kappa  \) assigns degree
2 to \( x_{i}^{*} \) and degree 1 to \( \xi ^{a} \) and \( \xi _{a}^{*} \).
Those functions \( \mu  \) that come from Lie algebroid structures on \( A \)
can be characterized by their \( (\epsilon ,\delta ) \)-bidegree, namely, 
\[
\epsilon (\mu )=1;\; \delta (\mu )=2;\; \kappa (\mu )=3\]
 The Poisson bracket \( \{\cdot ,\cdot \} \) has \( (\epsilon ,\delta ) \)-bidegree
\( (-1,-1) \), and hence total \( \kappa  \)-degree \( -2 \).
\end{rem}
\end{rem}
Similarly, a Lie algebroid structure on the dual bundle \( A^{*} \) gives rise
to a linear function \( \gamma =h_{d_{A^{*}}} \) on \( T^{*}\Pi A^{*} \) satisfying
\( \{\gamma ,\gamma \}=0 \). By (\ref{eqn:dA*}), it is given in local coordinates
\( (x^{i},\theta _{a},x^{*}_{i},\theta ^{a}_{*}) \) by 
\begin{equation}
\label{eqn:gamma}
\gamma =\theta _{a}\bar{A}^{ai}(x)x^{*}_{i}-\half \bar{C}^{ab}_{c}(x)\theta _{a}\theta _{b}\theta ^{c}_{*}
\end{equation}
Since \emph{a priori} the functions \( \mu  \) and \( \gamma  \) live on different
manifolds, it seems unclear at this point how to express the compatibility condition
between them in case \( (A,A^{*}) \) is a Lie bialgebroid. It is also not clear
how to express the Schouten brackets \( [\cdot ,\cdot ]_{A} \) and \( \br _{A^{*}} \)
in this formalism. Fortunately, it turns out that the supermanifolds \( T^{*}\Pi A \)
and \( T^{*}\Pi A^{*} \) are canonically symplectomorphic, via the \emph{Legendre
transform.}

\section{The Legendre transform\label{sec:legendre}}

The Legendre transform is widely known in classical mechanics for its crucial
role in providing a transition from the Lagrangian to the Hamiltonian formalism.
Recall that if \( M \) is the configuration space of a classical mechanical
system, \( l=l(q,\dot{q})\in \cinf (TM) \) the Lagrangian function, then the
dynamics of the system are given by the \emph{Euler-Lagrange equations} 
\[
\der{l}{q^{i}}(q(t),\dot{q}(t))-\frac{d}{dt}\der{l}{\dot{q}^{i}}(q(t),\dot{q}(t))=0\]
 satisfied by a classical trajectory \( q=q(t) \). Then one introduces the
momenta \( p_{i}\in \cinf (T^{*}M) \) by 
\[
p_{i}=\der{l}{\dot{q}^{i}}(q,\dot{q})\]
 Suppose the Lagrangian \( l \) is strongly nondegenerate in the sense that
the above equations have a unique solution \( \dot{q}^{i}=\dot{q}^{i}(q,p) \).
Then one can define the Hamiltonian function \( h=h(q,p)\in \cinf (T^{*}M) \)
as the Legendre transform of \( l \), i.e. 
\[
h(q,p)=\dot{q}^{i}p_{i}-l(q,\dot{q})\]
 where we substitute \( \dot{q}^{i}=\dot{q}^{i}(q,p) \). The Euler-Lagrange
equations are equivalent to \emph{Hamilton's equations} 
\[
\frac{dq^{i}}{dt}=\der{h}{p_{i}};\; \; \; \frac{dp_{i}}{dt}=-\der{h}{q^{i}}\]
 In 1977, W.M. Tulczyjew \cite{Tul} gave a geometric interpretation of the
Legendre transform as a canonical symplectomorphism between the cotangent bundles
\( T^{*}(TM) \) and \( T^{*}(T^{*}M) \). It turns out that in Tulczyjew's
construction one can replace \( TM \) with an arbitrary vector bundle \( A \)
or a supermanifold \( \Pi A \). We shall now describe this construction.

Let \( P \) be manifold, \( Q\subset P \) a submanifold, \( f\in \cinf (Q) \).
This data gives rise to a Lagrangian submanifold of \( L\subset T^{*}P \) as
follows: 
\[
L=\{\xi \in T^{*}P|\pi _{P}(\xi )=x\in Q;\: \xi (v)=df(v)\: \forall v\in T_{x}Q\},\]
where \( \pi _{P} \) denotes the canonical projection \( T^{*}P\rightarrow P \).
If \( f=0 \), \( L \) is just the conormal bundle to \( Q \); if \( Q=P \),
\( L \) is just the image of \( df \). 

We are interested in the following special case. Let \( A \) be the total space
of a vector bundle \( A\rightarrow M \), \( A^{*} \) that of the dual bundle.
Consider their fibre product, i.e. the total space of the Whitney sum \( Q=A\oplus A^{*}\subset A\times A^{*}=P \).
On \( Q \), there is a canonical evaluation function \( f=-ev\in \cinf (Q) \)
given by \( f(v,\xi )=-v(\xi ) \). This gives rise to a Lagrangian submanifold
\( L\subset T^{*}(A\times A^{*})\simeq T^{*}A\times \overline{T^{*}A^{*}} \),
where the bar denotes the opposite symplectic structure and the isomorphism
is given by the ``Schwartz transform'' 
\[
S((x,y),(\xi ,\eta ))=((x,\xi ),(y,-\eta ))\]
(see \cite{WeinsteinBates} for an explanation of this name). This \( L\subset T^{*}A\times \overline{T^{*}A^{*}} \)
is the graph of a symplectomorphism \( L:T^{*}A\rightarrow T^{*}A^{*} \) that
can be interpreted as a geometric version of the Legendre transform. It is given
in local coordinates simply by 
\[
((x,v),(p,\xi ))\longmapsto ((x,-\xi ),(p,v))\]

\begin{example}
\label{eg:LegTransf1}Let \( A=TM \), \( l\in \cinf (TM) \) a strongly nondegenerate
Lagrangian. Then the image of \( TM \) under \( -dl:TM\rightarrow T^{*}TM \)
followed by \( L:T^{*}(TM)\rightarrow T^{*}(T^{*}M) \) coincides with the image
of \( dh:T^{*}M\rightarrow T^{*}(T^{*}M) \), where \( h \) is the classical
Legendre transform of \( l \).
\end{example}
What is important for our purposes is that all of the above carries over to
supermanifolds without change, as long as the function \( f \) is even, otherwise
\( df \) is a section not of \( T^{*}Q \) but of \( \Pi T^{*}Q \). On \( \Pi (A\oplus A^{*}) \)
there is a canonical even function \( ev\in \cinf (\Pi (A\oplus A^{*})) \)
given in local coordinates by 
\[
ev(x,\xi ,\theta )=\xi ^{a}\theta _{a},\]
 giving rise to a canonical symplectomorphism \( L:T^{*}\Pi A\rightarrow T^{*}\Pi A^{*} \)
via the above construction. In local coordinates, 
\begin{equation}
\label{eqn:LegTransf}
L(x,\xi ,x^{*},\xi ^{*})=(x,\xi ^{*},x^{*},\xi )
\end{equation}
 In other words, the fibre coordinates \( \xi ^{a} \) on \( \Pi A \) become
conjugate to the fibre coordinates \( \theta _{a} \) on \( \Pi A^{*} \), and
vice versa. In a way, this local description is more illuminating than the canonical
geometric construction above, but we choose to present the geometric construction
rather than go through a proof that (\ref{eqn:LegTransf}) does not depend on
a choice of local coordinates.

\begin{example}
\label{eg:LegTransf2}Consider the supermanifold \( \Pi TM \) and a 2-form
\( \omega  \) on \( M \) viewed as a quadratic function on \( \Pi TM \).
If \( \omega  \) is nondegenerate, then the image of \( \Pi TM \) under \( d\omega :\Pi TM\rightarrow T^{*}\Pi TM \)
followed by the Legendre transform \( L:T^{*}\Pi TM\rightarrow T^{*}\Pi T^{*}M \)
coincides with the image of \( d\pi :\Pi T^{*}M\rightarrow T^{*}\Pi T^{*}M \),
where \( \pi \in \cinf (\Pi T^{*}M) \) is the bivector field given by the inverse
of \( \omega  \). Indeed, if \( \omega =\half \omega _{ab}(x)\xi ^{a}\xi ^{b} \),
then the image of \( d\omega  \) in \( T^{*}\Pi TM \) is given by 
\[
\begin{array}{rcl}
x^{*}_{c} & = & \der{\omega }{x^{c}}=\half \der{\omega _{ab}}{x^{c}}\xi ^{a}\xi ^{b}\\
\xi ^{*}_{c} & = & \der{\omega }{\xi ^{c}}=\omega _{cb}(x)\xi ^{b}
\end{array}\]
 The second equation can be solved for \( \xi  \) if and only if \( \omega  \)
is nondegenerate, in which case 
\[
\xi ^{a}=\pi ^{ab}(x)\xi ^{*}_{b}\]
where \( \pi ^{ab}(x)\omega _{bc}(x)=\delta ^{a}_{c} \). Applying the Legendre
transform (\ref{eqn:LegTransf}), we get 
\[
\begin{array}{rcl}
x_{c}^{*} & = & \half \der{\omega _{ab}}{x^{c}}\theta ^{a}_{*}\theta _{*}^{b}\\
\theta _{c} & = & \omega _{cb}(x)\theta _{*}^{b}
\end{array}\]
 and 
\[
\theta ^{a}_{*}=\pi ^{ab}(x)\theta _{b}\]
 Setting 
\[
-\pi (x,\theta )=\xi ^{a}\theta _{a}-\omega (x,\xi )=-\half \pi ^{ab}(x)\theta _{a}\theta _{b}\]
 we get 
\[
\begin{array}{rcl}
x^{*}_{c} & = & \der{\pi }{x^{c}}\\
\theta _{*}^{c} & = & \der{\pi }{\theta _{c}}
\end{array}\]
Note that here \( d \) denotes the deRham differential of functions on \( \Pi TM \)
(or \( \Pi T^{*}M \)), not of forms on \( M \)! Notice also how this example
parallels Example \ref{eg:LegTransf1}.
\end{example}
Let us now derive some properties of the Legendre transform that will be useful
in what follows. We begin by drawing the following diagram: 
\begin{equation}
\label{diagram:Legendre}
\begin{array}{ccc}
T^{*}\Pi A & \stackrel{{L}}{\longrightarrow } & T^{*}\Pi A^{*}\\
\left\downarrow \pi \right.  &  & \left. \bar{\pi }\right\downarrow \\
\Pi A &  & \Pi A^{*}
\end{array}
\end{equation}
 where \( \pi  \) and \( \bar{\pi } \) are the canonical projections. It is
obvious that \( \pi _{A}\smalcirc \pi =\pi _{A^{*}}\smalcirc \bar{\pi }\smalcirc L \),
where \( \pi _{A}:\Pi A\rightarrow M \) and \( \pi _{A^{*}}:\Pi A^{*}\rightarrow M \)
are the canonical projections; therefore, by abstract nonsense the diagram above
gives rise to a projection to the fibered product 
\begin{equation}
\label{eqn:proj}
T^{*}\Pi A\stackrel{p}{\longrightarrow }\Pi (A\oplus A^{*})
\end{equation}
 More specifically, if \( \xi \in \Gamma (A^{*}) \) viewed as a linear function
on \( \Pi A \) or on \( \Pi (A\oplus A^{*}) \), then \( p^{*}\xi =\pi ^{*}\xi  \);
on the other hand, \( \xi  \) also gives rise to a vector field \( i_{\xi } \)
on \( \Pi A^{*} \), the interior derivative, hence a linear hamiltonian \( h_{i_{\xi }} \)
on \( T^{*}\Pi A^{*} \). Similarly, an \( X\in \Gamma (A) \) can be viewed
as either a linear function on \( \Pi A^{*} \), pulled back to \( T^{*}\Pi A^{*} \),
or a vector field \( i_{X} \) on \( \Pi A \) lifted to the hamiltonian \( h_{i_{X}} \)
on \( T^{*}\Pi A \). These functions are related by the Legendre transform,
according to the following 

\begin{lem}
\label{lemma:intertwine}For \( X\in \Gamma (A), \) \( \xi \in \Gamma (A^{*}) \),
\[
\begin{array}{ccc}
L^{*}\bar{\pi }^{*}X & = & h_{i_{X}}\\
L^{*}h_{i_{\xi }} & = & \pi ^{*}\xi 
\end{array}\]
 
\end{lem}
\begin{proof}
If locally \( X=X^{a}(x)\theta _{a} \), \( \xi =f_{a}(x)\xi ^{a} \), then
\( i_{X}=X^{a}(x)\der{}{\xi ^{a}} \), hence
\[
h_{i_{X}}=X^{a}(x)\xi _{a}^{*}=L^{*}(X^{a}(x)\theta _{a})=L^{*}\bar{\pi }^{*}X\]
 by (\ref{eqn:LegTransf}), and similarly for \( \xi  \).
\end{proof}
In particular, \( p^{*}\xi =\pi ^{*}\xi  \), \( p^{*}X=L^{*}\bar{\pi }^{*}X=h_{i_{X}} \)
and, of course, if \( f\in \cinf (M) \), \( p^{*}f=\pi ^{*}\pi _{A}^{*}f=\bar{\pi }^{*}\pi ^{*}_{A^{*}}f \).
This can be interpreted as follows. On \( \Pi (A\oplus A^{*}) \) there is a
natural even Poisson structure given by the canonical inner product on \( A\oplus A^{*} \).
The corresponding Poisson bracket is just the fiberwise big bracket described
in the beginning of this section. The symplectic leaves are the fibres that
inherit the big bracket, so the pullbacks of functions on \( M \) are the Casimir
functions. We have the following 

\begin{cor}
\label{cor:PoissonMap}The projection \( p \) (\ref{eqn:proj}) is a Poisson
map.
\end{cor}
\begin{proof}
Immediate from (\ref{eqn:canonical_{r}elations}) and Lemmas \ref{lemma:cotbundle}
and \ref{lemma:intertwine}.
\end{proof}
In other words, \( T^{*}\Pi A \) is a symplectic realization of the Poisson
supermanifold \( \Pi (A\oplus A^{*}) \). This fact will be useful in dealing
with Courant algebroids. But for now, we need one more construction to be able
to deal with the Schouten brackets and Lie bialgebroids.

\section{Derived brackets\label{sec:derived}}

Let \( (\A ,\br _{\A },d) \) be an even or odd differential Lie superalgebra.
That is, 

\begin{itemize}
\item \( \A =\A _{0}\oplus \A _{1} \) is a \( \Z _{2} \)-graded vector space; 
\item \( \br _{\A } \) is skew-symmetric of parity \( \epsilon \in \Z _{2} \), i.e.
\( [\A _{i},\A _{j}]_{\A }\subset \A _{i+j+\epsilon } \) and 
\[
[a,b]_{\A }=-(-1)^{(\tilde{a}+\epsilon )(\tilde{b}+\epsilon )}[b,a]_{\A }\]
 for all \( a\in \A _{\tilde{a}} \), \( b\in \A _{\tilde{b}} \);
\item \( \br _{\A } \) satisfies the Jacobi identity 
\[
[a,[b,c]_{\A }]_{\A }=[[a,b]_{\A },c]_{\A }+(-1)^{(\tilde{a}+\epsilon )(\tilde{b}+\epsilon )}[b,[a,c]_{\A }]_{\A }\]
 
\item \( d:\A \rightarrow \A  \) is an odd derivation of \( \br _{\A } \): 
\[
d[a,b]_{\A }=[da,b]_{\A }+(-1)^{\tilde{a}+\epsilon }[a,db]_{\A }\]
 satisfying \( d^{2}=0 \).
\end{itemize}
One defines the \emph{derived bracket} on \( \A  \) as follows: 
\[
a\circ _{d}b=[ad_{a},d]b=(-1)^{\tilde{a}+1}[da,b]_{\A }\]
 where the bracket in the middle is the (super)commutator of derivations of
\( \A  \) and \( ad_{a}=[a,\cdot ]_{\A } \). The derived bracket has parity
\( \epsilon +1 \). It is not necessarily skew-symmetric; its skew-symmetrization
\[
[a,b]_{d}=\half (a\circ _{d}b-(-1)^{(\tilde{a}+1)(\tilde{b}+1)}b\circ _{d}a)\]
 is also sometimes called the derived bracket, but the non skew-symmetric version
is, in some sense, more fundamental and enjoys many nice properties.\footnote{
It is somewhat unfortunate that what we call the derived bracket is not denoted
by a bracket; nevertheless, we feel that the bracket notation ought to be reserved
for skew-symmetric operations. In \cite{KS5} \( f_{d} \) was used to denote
the derived bracket regardless of skew-symmetry. 
} They are summarized in the following 

\begin{lem}
\label{lemma:PropOfDer}The derived bracket has the following properties:
\begin{enumerate}
\item \( a\circ _{d}(b\circ _{d}c)=(a\circ _{d}b)\circ _{d}c+(-1)^{(\tilde{a}+\epsilon +1)(\tilde{b}+\epsilon +1)}b\circ _{d}(a\circ _{d}c) \)
\( \forall a,b,c\in \A  \).
\item \( d(a\circ _{d}b)=(da)\circ _{d}b+(-1)^{\tilde{a}+\epsilon +1}a\circ _{d}(db) \)
\( \forall a,b\in \A  \).
\item \( a\circ _{d}b=[a,b]_{d}+\frac{(-1)^{\tilde{a}+\epsilon +1}}{2}d[a,b]_{\A } \)
\end{enumerate}
\end{lem}
\begin{proof}
A straightforward computation, carried out in \cite{KS5}. Property 1 depends
both on the Jacobi identity for \( \br _{\A } \) and \( d^{2}=0 \).
\end{proof}
The first two properties imply that \( (\A ,\circ _{d},d) \) is a \emph{differential
Leibniz superalgebra} in the sense of Loday \cite{Loday1} (they are called
\emph{Loday algebras} in \cite{KS5}). The third property implies that \( \circ _{d} \)
is skew-symmetric up to a \( d \)-coboundary. Notice how these properties resemble
some of the properties of Courant algebroids (Definition \ref{def:quasi-algebroid2}
and (\ref{eqn:symplusskewsym})). 

\begin{cor}
\label{cor:AbSub}Let \( \B \subset \A  \) be an abelian subalgebra (with respect
to \( \br _{\A } \)) closed under \( \circ _{d} \); then the restriction of
\( \circ _{d} \) to \( \B  \) is skew-symmetric and \( (\B ,\br _{d}) \)
is a Lie superalgebra of parity \( \epsilon +1 \); if, moreover, \( d\B \subset \B  \),
then \( (\B ,\br _{d},d) \) is a differential Lie superalgebra.
\end{cor}
\begin{rem*}
Sometimes \( \A  \) also has a (super)commutative algebra structure such that
\( \br _{\A } \) is a derivation of the multiplication in each argument. If
\( \epsilon =0 \) (i.e. \( \A  \) is an even Poisson superalgebra), then \( \B  \)
becomes an odd Poisson superalgebra, i.e. a Gerstenhaber algebra. On the other
hand, if \( \A  \) is a Gerstenhaber algebra, \( \B  \) is an even Poisson
superalgebra \cite{KS5}.
\end{rem*}
\begin{example}
\label{eg:EvenPoisson}Consider \( \A =\cinf (\Pi T^{*}M)=\Gamma (\bigwedge TM)=\X ^{*}(M) \),
the multivector fields on a supermanifold \( M \) endowed with the Schouten
bracket \( \br  \); pick a quadratic function (a bivector field) \( \pi  \)
satisfying \( [\pi ,\pi ]=0 \) and consider the inner derivation \( d_{\pi }=[\pi ,\cdot ] \);
\( \A  \) is a differential Gerstenhaber algebra. Then \( \B =\cinf (M) \)
is an abelian subalgebra of \( \A  \) stable under \( d_{\pi } \), and the
derived bracket on \( \B  \)
\[
\{f,g\}=(-1)^{\tilde{f}}[[\pi ,f],g]\]
 is precisely the Poisson bracket generated by the bivector field \( \pi  \). 
\begin{example}
\label{eg:AlgSchouten}If \( (\g ,\mu ,\gamma ) \) is a Lie bialgebra, let
\( \A =\wedge (\g \oplus \g ^{*}) \) with \( \br _{\A }=\{\cdot ,\cdot \} \),
the big bracket. Then \( \wedge \g  \) is an abelian subalgebra stable under
the differential \( \{\mu ,\cdot \} \), while \( \wedge \g ^{*} \) is an abelian
subalgebra stable under \( \{\gamma ,\cdot \} \). The corresponding derived
brackets give the \emph{algebraic Schouten brackets,} generalizing the formulas
(\ref{eqn:brackets}) \cite{KS3}. Notice how the Drinfeld double bracket (\ref{eqn:double1})
is generated by \( \theta =\mu +\gamma  \) as a derived bracket: although \( \g \oplus \g ^{*} \)
is not closed under \( \{\cdot ,\cdot \} \), it is closed under the derived
bracket.
\end{example}
\end{example}
We will show below that the Schouten brackets associated to Lie algebroids,
as well as the Courant bracket (\ref{eqn:double2}), arise in exactly the same
way.

\section{Schouten brackets, Lie bialgebroids and the Drinfeld double\label{sec:schouten}}

The concept of a derived bracket enables us to define the Schouten brackets
and recast the notion of Lie bialgebroid in the supermanifold context. Let \( d_{A^{*}} \)
be a homological vector field on \( \Pi A^{*} \) giving rise to a Lie algebroid
structure on \( A^{*}\rightarrow M \); let \( \gamma =h_{d_{A^{*}}} \) be
the corresponding linear hamiltonian on \( T^{*}\Pi A^{*} \), and consider
its Legendre transform \( L^{*}\gamma \in \cinf (T^{*}\Pi A) \). By (\ref{eqn:gamma})
and (\ref{eqn:LegTransf}), 
\begin{equation}
\label{eqn:L*gamma}
L^{*}\gamma =\bar{A}^{ai}(x)x^{*}_{i}\xi _{a}^{*}-\half \xi ^{c}\bar{C}^{ab}_{c}(x)\xi ^{*}_{a}\xi ^{*}_{b}
\end{equation}

\begin{rem}
\label{rem:gradings2}Notice that \( L^{*}\gamma  \) is fiberwise quadratic,
i.e. \( \epsilon (L^{*}\gamma )=2 \); on the other hand, \( \delta (L^{*}\gamma )=1 \),
so the total degree \( \kappa (L^{*}\gamma ) \) is again 3.  This characterizes
those functions on \( T^{*}\Pi A \) that come from Lie algebroid structures
on \( A^{*} \). In fact, the grading \( \delta  \) is seen to correspond to
the momentum grading \( \epsilon ^{*} \) on \( T^{*}\Pi A^{*} \) under \( L \),
whereas the \( \epsilon  \)-grading corresponds to \( \delta ^{*} \). Thus,
the Legendre transform interchanges the \( \epsilon  \) and \( \delta  \)
gradings and preserves the total grading \( \kappa  \).
\end{rem}
Since \( L \) is a symplectomorphism, we have 
\[
\{L^{*}\gamma ,L^{*}\gamma \}=L^{*}\{\gamma ,\gamma \}=0,\]
 hence \( (\cinf (T^{*}\Pi A),\{\cdot ,\cdot \},\{L^{*}\gamma ,\cdot \}) \)
is a differential Lie superalgebra, and we can consider the derived bracket.
It turns out that the abelian subalgebra \( \pi ^{*}\cinf (\Pi A) \) is closed
under the derived bracket, and the restriction of the derived bracket coincides
with the Schouten bracket \( \br _{A^{*}} \). More precisely, we have 

\begin{lem}
\label{lemma:Schouten}Let \( \xi ,\eta \in \cinf (\Pi A)=\Gamma (\bigwedge A^{*}) \).
Then 
\[
\pi ^{*}[\xi ,\eta ]_{A^{*}}=(-1)^{\tilde{\xi }+1}\{\{L^{*}\gamma ,\pi ^{*}\xi \},\pi ^{*}\eta \}\]
 
\end{lem}
\begin{proof}
The skew-symmetry and derivation property are consequences of Corollary \ref{cor:AbSub}.
Hence, we only need to consider fiberwise constant and fiberwise linear functions,
i.e. elements of \( \cinf (M) \) and \( \Gamma (A^{*}) \). We have:
\[
\begin{array}{l}
\{\{L^{*}\gamma ,\pi ^{*}f\},\pi ^{*}g\}=\{L^{*}\{\gamma ,\bar{\pi }^{*}f\},\pi ^{*}g\}=\{L^{*}\{h_{d_{A^{*}}},\bar{\pi }^{*}f\},\pi ^{*}g\}=\\
=\{L^{*}\bar{\pi }^{*}d_{A^{*}}f,\pi ^{*}g\}=L^{*}\{\bar{\pi }^{*}d_{A^{*}}f,\bar{\pi }^{*}g\}=0=\pi ^{*}[f,g]_{A^{*}}
\end{array}\]
 for all \( f,g\in \cinf (M) \); 
\[
\begin{array}{l}
\{\{L^{*}\gamma ,\pi ^{*}\xi \},\pi ^{*}f\}=\{L^{*}\{h_{d_{A^{*}}},h_{i_{\xi }}\},\pi ^{*}f\}=L^{*}\{h_{[d_{A^{*}},i_{\xi }]},\bar{\pi }^{*}f\}=\\
=L^{*}\bar{\pi }^{*}(L_{\xi }^{A^{*}}f)=\pi ^{*}(a_{*}(\xi )f)=\pi ^{*}[\xi ,f]_{A^{*}}
\end{array}\]
 for all \( f\in \cinf (M) \), \( \xi \in \Gamma (A^{*}) \). And finally,

\[
\begin{array}{l}
\{\{L^{*}\gamma ,\pi ^{*}\xi \},\pi ^{*}\eta \}=\{L^{*}\{h_{d_{A^{*}}},h_{i_{\xi }}\},\pi ^{*}\eta \}=L^{*}\{h_{[d_{A^{*}},i_{\xi }]},h_{i_{\eta }}\}=\\
=L^{*}h_{[[d_{A^{*}},i_{\xi }],i_{\eta }]}=L^{*}h_{[L_{\xi }^{A^{*}},i_{\eta }]}=L^{*}h_{i_{[\xi ,\eta ]_{A^{*}}}}=\pi ^{*}[\xi ,\eta ]_{A^{*}}
\end{array}\]
 for all \( \xi ,\eta \in \Gamma (A^{*}) \). We have made repeated use of the
commutation relations (\ref{eqn:canonical_{r}elations}), Lemma \ref{lemma:cotbundle}
and Lemma \ref{lemma:intertwine}. 
\end{proof}
\begin{rem*}
Of course, the same is true for the Schouten bracket \( \br _{A} \) associated
to a Lie algebroid structure on \( A \), if we use \( (L^{-1})^{*}\mu  \)
where \( \mu =h_{d_{A}} \). 
\end{rem*}
\begin{example}
\label{eg:SchoutenBr}For any (super)manifold \( M \), a canonical Lie algebroid
structure on the bundle \( A^{*}=TM \) (Example \ref{eg:LieAlg-standard})
is given by the de Rham differential \( d \) on \( \Pi TM \). The corresponding
quadratic hamiltonian 
\[
L^{*}h_{d}=\theta _{*}^{a}x^{*}_{a}\]
 on \( T^{*}\Pi T^{*}M \) generates the Schouten bracket of multivector fields
on \( M \) as the derived bracket. 
\end{example}
\begin{rem}
\label{rem:Odd-Even} The Schouten bracket \( \br _{A^{*}} \) on the supermanifold
\( \Pi A \) is an odd Poisson structure (see Appendix). Examples \ref{eg:EvenPoisson},
\ref{eg:AlgSchouten} and \ref{eg:SchoutenBr} are special cases of the following
general phenomenon: even Poisson structures on a supermanifold \( M \) are
generated by bivector fields, i.e. even quadratic hamiltonians on the odd symplectic
supermanifold \( \Pi T^{*}M \), whereas odd Poisson structures are generated
by odd quadratic hamiltonians on the even symplectic supermanifold \( T^{*}M \)
(see the Appendix in \cite{Vor1}, also \cite{KS5}).
\end{rem}
We can now prove the following simple characterization of Lie bialgebroids.

\begin{prop}
\label{prop:Bialgebroid} A pair \( (A,A^{*}) \) of Lie algebroids in duality
is a Lie bialgebroid if and only if 
\begin{equation}
\label{eqn:Bialgebroid}
\{\mu ,L^{*}\gamma \}=0,
\end{equation}
 where \( \mu =h_{d_{A}} \), \( \gamma =h_{d_{A^{*}}} \). 
\end{prop}
\begin{proof}
We must show that \( d_{A} \) is a derivation of \( \br _{A^{*}} \) if and
only if (\ref{eqn:Bialgebroid}) holds. However, by Lemma \ref{lemma:cotbundle}
and Lemma \ref{lemma:Schouten}, we have 
\[
\begin{array}{rcl}
\pi ^{*}d_{A}[\xi ,\eta ]_{A^{*}} & = & \{\mu ,\pi ^{*}[\xi ,\eta ]_{A^{*}}\}=(-1)^{\tilde{\xi }+1}\{\mu ,\{\{L^{*}\gamma ,\pi ^{*}\xi \},\pi ^{*}\eta \}\\
 & = & (-1)^{\tilde{\xi }+1}(\{\{\mu ,\{L^{*}\gamma ,\pi ^{*}\xi \}\},\pi ^{*}\eta \}+\\
 & + & (-1)^{\tilde{\xi }+1}\{\{L^{*}\gamma ,\pi ^{*}\xi \},\{\mu ,\pi ^{*}\eta \}\})=\\
 & = & (-1)^{\tilde{\xi }+1}(\{\{\{\mu ,L^{*}\gamma \},\pi ^{*}\xi \},\pi ^{*}\eta \}-\\
 & - & \{\{L^{*}\gamma ,\{\mu ,\pi ^{*}\xi \}\},\pi ^{*}\eta \}+(-1)^{\tilde{\xi }+1}\{\{L^{*}\gamma ,\pi ^{*}\xi \},\pi ^{*}d_{A}\eta \})=\\
 & = & (-1)^{\tilde{\xi }}\{\{L^{*}\gamma ,\pi ^{*}d_{A}\xi \},\pi ^{*}\eta \}+\{\{L^{*}\gamma ,\pi ^{*}\xi \},\pi ^{*}d_{A}\eta \}+\\
 & + & (-1)^{\tilde{\xi }+1}\{\{\{\mu ,L^{*}\gamma \},\pi ^{*}\xi \},\pi ^{*}\eta \}=\\
 & = & \pi ^{*}([d_{A}\xi ,\eta ]_{A^{*}}+(-1)^{\tilde{\xi }+1}[\xi ,d_{A}\eta ]_{A^{*}})+\\
 & + & (-1)^{\tilde{\xi }+1}\{\{\{\mu ,L^{*}\gamma \},\pi ^{*}\xi \},\pi ^{*}\eta \}
\end{array}\]
Since \( \{\mu ,L^{*}\gamma \} \) is fiberwise quadratic, the second term in
the last expression vanishes if and only if \( \{\mu ,L^{*}\gamma \}=0 \).
The statement follows by the injectivity of \( \pi ^{*} \).
\end{proof}
\begin{cor}
\( (A,A^{*}) \) is a Lie bialgebroid if and only if \( (A^{*},A) \) is.
\end{cor}
\begin{proof}
The Legendre transform \( L \) is a symplectomorphism.
\end{proof}
Now set \( \theta =\mu +L^{*}\gamma  \). Clearly, \( (A,A^{*}) \) is a Lie
bialgebroid if and only if 
\begin{equation}
\label{eqn:master2}
\{\theta ,\theta \}=0
\end{equation}
 This motivates the following 

\begin{defn}
Given a Lie bialgebroid \( (A,A^{*}) \), its \emph{Drinfeld double} is \( T^{*}\Pi A \)
together with the homological vector field \( D=\{\theta ,\cdot \} \).
\end{defn}
\begin{example}
If \( (\g ,\g ^{*}) \) is a Lie bialgebra, \( \cinf (T^{*}\Pi \g )=\bigwedge (\g \oplus \g ^{*})^{*} \)and
\( D \) is the Chevalley-Eilenberg differential in the standard complex of
the Drinfeld double Lie algebra \( \g \oplus \g ^{*} \) (\ref{eqn:doubleofbialgebra}).
\begin{example}
Let \( A=M\times \g  \) be the action Lie algebroid corresponding to a Lie
algebra action \( \rho :\g \rightarrow \X (M) \) (Example \ref{eg:LieAlg-action}).
View \( (A,A^{*}) \) as a Lie bialgebroid with the trivial structure on \( A^{*} \).
Then \( \cinf (T^{*}\Pi A)=\cinf (M)\otimes \bigwedge (\g \oplus \g ^{*}) \),
and \( D=\{\mu ,\cdot \} \) coincides with the classical BRST differential
associated to the hamiltonian lift of \( \rho  \) to \( T^{*}M \). Indeed,
recall that the classical BRST differential \( d \) is the sum of the Chevalley-Eilenberg
differential \( \delta  \) for the Lie algebra \( \g  \) with values in the
module \( \wedge \g \otimes \cinf (T^{*}M) \), and the Koszul differential
\( \partial  \) for the zero level of the momentum map (in this case, the ideal
generated by the linear hamiltonians \( \{h_{\rho (X)}|X\in \g \} \)) \cite{KoSt}.
On generators, we have 
\[
df(Y)=\{h_{\rho (Y)},\pi ^{*}f\}=\pi ^{*}\rho (Y)f=\pi ^{*}d_{A}f(Y)\]
 for \( f\in \cinf (M) \) and for all \( Y\in \g  \). Hence, \( df=\{\mu ,\pi ^{*}f\} \);
\[
d\xi (X,Y)=-\xi ([X,Y])=\pi ^{*}d_{A}\xi (X,Y)\]
 for \( \xi \in \g ^{*} \) (a constant section of \( A^{*} \)). Hence, \( d\xi =\{\theta ,\pi ^{*}\xi \} \);
\[
dh_{v}(Y)=\{h_{\rho (Y)},h_{V}\}=h_{[\rho (Y),v]}\]
 for any vector field \( v \) on \( M \). Hence, \( dh_{v}=h_{[d_{A},v]}=\{\mu ,h_{v}\} \);
and finally, for \( X\in \g  \), 
\[
\delta X(Y)=ad_{Y}X=-[X,Y];\; \partial X=h_{\rho (X)}\]
 so 
\[
dX=h_{\rho (X)+ad^{*}_{X}}=h_{[d_{A},i_{X}]}=\{\mu ,h_{i_{X}}\}\]
 where \( ad^{*}_{X} \) is viewed as a vector field on \( \Pi \g  \). Thus
the BRST differential \( d \) coincides with our differential \( \{\mu ,\cdot \} \).
\begin{example}
Let \( \g  \) be a Lie algebra, then \( \g ^{*} \) is a Poisson manifold,
with the canonical linear Poisson structure. Consider the corresponding Lie
bialgebroid \( A=T\g ^{*}\simeq \g ^{*}\times \g ^{*} \), \( A^{*}=T^{*}\g ^{*}\simeq \g ^{*}\times \g  \).
Then \( \cinf (T^{*}\Pi A)=\cinf (\g ^{*}\oplus \g )\otimes \bigwedge (\g ^{*}\oplus \g ) \).
A choice of a basis \( \{e_{a}\} \) of \( \g  \) and a dual basis \( \{e^{a}\} \)
of \( \g ^{*} \) gives rise to coordinates \( (u_{a},\theta _{a}) \) on \( \Pi A \)
and \( (u_{a},\xi ^{a}) \) on \( \Pi A^{*} \). Then the differentials are
the deRham differential 
\[
d=\theta _{a}\der{}{u_{a}}\]
 on \( \Pi A \), and the Poisson differential 
\[
d_{\pi }=u_{a}C^{a}_{bc}\xi ^{b}\der{}{u_{c}}-\half C_{ab}^{c}\xi ^{a}\xi ^{b}\der{}{\xi ^{c}}\]
 where \( C_{ab}^{c} \) are the structure constants of \( \g  \). Thus, 
\[
\theta =u_{a}C_{bc}^{a}\theta ^{b}_{*}u^{c}_{*}-\half C^{c}_{ab}\theta ^{a}_{*}\theta ^{b}_{*}\theta _{c}+\theta _{a}u^{a}_{*}\]
 on \( T^{*}\Pi A \) and 
\[
\begin{array}{rcl}
D=\{\theta ,\cdot \} & = & (\theta _{c}+u_{a}C^{a}_{bc}\theta ^{b}_{*})\der{}{u_{c}}+(\theta ^{b}_{*}C_{bc}^{a}\theta _{a}+u_{a}C^{a}_{cb}u^{b}_{*})\der{}{\theta _{c}}+\\
 & + & C_{ab}^{c}u_{*}^{a}\theta ^{b}_{*}\der{}{u^{c}_{*}}+(u^{c}_{*}-\half C_{ab}^{c}\theta ^{a}_{*}\theta ^{b}_{*})\der{}{\theta ^{c}_{*}}
\end{array}\]
 Notice that the fibre over the origin, given by the equations \( u_{a}=\theta _{a}=0 \),
is a Lagrangian submanifold \( F \) stable under \( D \). The restriction
of \( D \) to \( F \) is 
\[
D=C_{ab}^{c}u_{*}^{a}\theta ^{b}_{*}\der{}{u_{*}^{c}}+(u_{*}^{c}-\half C^{c}_{ab}\theta ^{a}_{*}\theta ^{b}_{*})\der{}{\theta ^{c}_{*}}\]
 The algebra of polynomial functions on \( F \), isomorphic to \( S\g ^{*}\otimes \wedge \g ^{*}=\R [u^{a}_{*},\theta ^{a}_{*}] \),
is known as the \emph{Weil algebra \( W(\g ) \),} while the restricted differential
\( D \) above is the \emph{Weil differential.} This is Weil's deRham model
for the universal classifying space \( BG \), at least when the group \( G \)
(whose Lie algebra is \( \g  \)) is compact \cite{AtBott}. Its appearance
in this context is a mystery to us. Notice, however, that our \( \kappa  \)-grading,
\( \kappa (u^{a}_{*})=2 \), \( \kappa (\theta ^{a}_{*})=1 \), is consistent
with the grading in the Weil algebra. Notice also that \( T^{*}\Pi A\simeq T^{*}F \)
by a Legendre transform; after this identification, the ``full'' \( D \)
is just the hamiltonian lift of the Weil differential. 
\end{example}
\end{example}
\end{example}

\section{The Courant Algebroid\label{sec:courantagain}}

The Courant algebroid constructed in \cite{LWX1} as the double of a Lie bialgebroid
\( (A,A^{*}) \) (see Example \ref{eg:double}) can be recovered from the supermanifold
double \( (T^{*}\Pi A,D) \) via the derived bracket construction. We shall
view sections of \( A\oplus A^{*} \) and functions on \( M \) as functions
on \( \Pi (A\oplus A^{*}) \) and use the projection \( p \) (see (\ref{eqn:proj})). 

\begin{thm}
\label{thm:Courant=3DDerived}Let \( (E=A\oplus A^{*},\bform {\cdot }{\cdot },\circ ,\rho ) \)
be as in Example \ref{eg:double}. Then, for any \( e_{1},e_{2}\in \Gamma (E) \),
\( f\in \cinf (M) \) we have 
\begin{enumerate}
\item \( p^{*}\bform {e_{1}}{e_{2}}=\{p^{*}e_{1},p^{*}e_{2}\} \)
\item \( p^{*}\DD f=Dp^{*}f=\{\theta ,p^{*}f\} \) 
\item \( p^{*}(e_{1}\circ e_{2})=p^{*}e_{1}\circ _{D}p^{*}e_{2} \) 
\end{enumerate}
\end{thm}
\begin{proof}
(1) is just a restatement of Corollary \ref{cor:PoissonMap}; (2) follows by
computation: 
\[
\begin{array}{l}
p^{*}\DD f=p^{*}(d_{A}f+d_{A^{*}}f)=\pi ^{*}d_{A}f+L^{*}\bar{\pi }^{*}d_{A^{*}}f=\\
=\{\mu ,\pi ^{*}f\}+L^{*}\{\gamma ,\bar{\pi }^{*}f\}=\{\mu +L^{*}\gamma ,\pi ^{*}f\}=Dp^{*}f
\end{array}\]
 (3) takes a bit more work. We have  
\[
\begin{array}{rcl}
p^{*}X\circ _{D}p^{*}Y & = & \{\{\mu +L^{*}\gamma ,h_{i_{X}}\},h_{i_{Y}}\}=\\
 & = & \{h_{[d_{A},i_{X}]},h_{i_{Y}}\}+L^{*}\{\{\gamma ,\bar{\pi }^{*}X\},\bar{\pi }^{*}Y\}=\\
 & = & h_{[L_{X}^{A},i_{Y}]}+L^{*}\{\bar{\pi }^{*}d_{A^{*}}X,\bar{\pi }^{*}Y\}=\\
 & = & h_{i_{[X,Y]_{A}}}=p^{*}[X,Y]_{A}
\end{array}\]
for \( X,Y\in \Gamma (A) \);  

\[
\begin{array}{rcl}
p^{*}\xi \circ _{D}p^{*}\eta  & = & \{\{\mu +L^{*}\gamma ,\pi ^{*}\xi \},\pi ^{*}\eta \}=\\
 & = & \{\{\mu ,\pi ^{*}\xi \},\pi ^{*}\eta \}+L^{*}\{\{\gamma ,h_{i_{\xi }}\},h_{i_{\eta }}\}=\\
 & = & L^{*}\{h_{[d_{A^{*}},i_{\xi }]},h_{i_{\eta }}\}=L^{*}h_{[L^{A^{*}}_{\xi },i_{\eta }]}=\\
 & = & L^{*}h_{i_{[\xi ,\eta ]_{A^{*}}}}=\pi ^{*}[\xi ,\eta ]_{A^{*}}=p^{*}[\xi ,\eta ]_{A^{*}}
\end{array}\]
for \( \xi ,\eta \in \Gamma (A^{*}) \);  

\[
\begin{array}{rcl}
p^{*}X\circ _{D}p^{*}\eta  & = & \{\{\mu +L^{*}\gamma ,h_{i_{X}}\},\pi ^{*}\eta \}=\\
 & = & \{h_{[d_{A},i_{X}]},\pi ^{*}\eta \}+L^{*}\{\{\gamma ,\bar{\pi }^{*}X\},h_{i_{\eta }}\}=\\
 & = & \{h_{L^{A}_{X}},\pi ^{*}\eta \}+L^{*}\{\bar{\pi }^{*}d_{A^{*}}X,h_{i_{\eta }}\}=\\
 & = & \pi ^{*}L_{X}^{A}\eta -L^{*}\bar{\pi }^{*}i_{\eta }d_{A^{*}}X=\\
 & = & p^{*}(L^{A}_{X}\eta -i_{\eta }d_{A^{*}}X)
\end{array}\]
for \( X\in \Gamma (A) \), \( \eta \in \Gamma (A^{*}) \); and finally,  

\[
\begin{array}{rcl}
p^{*}\xi \circ _{D}p^{*}Y & = & \{\{\mu +L^{*}\gamma ,\pi ^{*}\xi \},h_{i_{Y}}\}=\\
 & = & \{\{\mu ,\pi ^{*}\xi \},h_{i_{Y}}\}+L^{*}\{h_{[d_{A^{*}},i_{\xi }]},\bar{\pi }^{*}Y\}=\\
 & = & \{\pi ^{*}d_{A}\xi ,h_{i_{Y}}\}+L^{*}\{h_{L^{A^{*}}_{\xi }},\bar{\pi }^{*}Y\}=\\
 & = & -\pi ^{*}i_{Y}d_{A}\xi +L^{*}\bar{\pi }^{*}L^{A^{*}}_{\xi }Y=\\
 & = & p^{*}(-i_{Y}d_{A}\xi +L^{A^{*}}_{\xi }Y)
\end{array}\]
This proves (3). We have made extensive use of the commutation relations (\ref{eqn:canonical_{r}elations})
and Lemmas \ref{lemma:cotbundle} and \ref{lemma:intertwine}.
\end{proof}
\begin{rem}
The Theorem above is true regardless of whether \( (A,A^{*}) \) is a Lie bialgebroid,
i.e. whether (\ref{eqn:master2}) holds; however, if it is the case, we can
use the differential Lie superalgebra \( (\cinf (T^{*}\Pi A),\{\cdot ,\cdot \},D) \)
and its derived bracket to prove that \( (A\oplus A^{*},\bform {\cdot }{\cdot },\circ ,\rho ) \)
actually is a Courant algebroid, thus recovering the doubling theorem of Liu,
Weinstein and Xu: 
\end{rem}
\begin{thm}
\label{thm:Bialg=3D>Courant}If \( (A,A^{*}) \) is a Lie bialgebroid, then
\( (A\oplus A^{*},\bform {\cdot }{\cdot },\circ ,\rho ) \) is a Courant algebroid. 
\end{thm}
\begin{proof}
We need to verify properties 1-5 of Definition \ref{def:quasi-algebroid2}.
Since \( p \) is a Poisson map, we can embed sections of \( A\oplus A^{*} \)
and functions on \( M \) into \( \cinf (T^{*}\Pi A) \) using \( p^{*} \)
as above and carry out all the computations up in \( \cinf (T^{*}\Pi A) \).
We shall identify \( e_{i}\in \Gamma (A\oplus A^{*}) \) and \( f\in \cinf (M) \)
with their images under \( p^{*} \). 

Now, it follows that properties 1 (the Leibniz-Jacobi identity) and 4 (about
the symmetric part) are just consequences of the properties of the derived bracket
on a differential Lie superalgebra (Lemma \ref{lemma:PropOfDer}). On the other
hand, Property 3, 
\[
e_{1}\circ fe_{2}=f(e_{1}\circ e_{2})+(\rho (e_{1})f)e_{2}\]
 translates, by Theorem \ref{thm:Courant=3DDerived}, into 
\[
\{\{\theta ,e_{1}\},fe_{2}\}=\{\{\theta ,e_{1}\},f\}e_{2}+f\{\{\theta ,e_{1}\},e_{2}\},\]
 but this is obvious. Property 5, 
\[
\rho (e)\bform {e_{1}}{e_{2}}=\bform {e\circ e_{1}}{e_{2}}+\bform {e_{1}}{e\circ e_{2}}\]
 translates into 
\[
\{e,\{\theta ,\{e_{1},e_{2}\}\}\}=\{\{\{\theta ,e\},e_{1}\},e_{2}\}+\{e_{1},\{\{\theta ,e\},e_{2}\}\}\]
 However, by the Jacobi identity for \( \{\cdot ,\cdot \} \), 
\[
\begin{array}{rcl}
\{e,\{\theta ,\{e_{1},e_{2}\}\}\} & = & \{\{e,\theta \},\{e_{1},e_{2}\}\}-\{\theta ,\{e,\{e_{1},e_{2}\}\}\}=\\
 & = & \{\{\{\theta ,e\},e_{1}\},e_{2}\}+\{e_{1},\{\{\theta ,e\},e_{2}\}\}
\end{array}\]
 since \( \{e,\{e_{1},e_{2}\}\}=0 \) because \( \{e_{1},e_{2}\}\in \cinf (M) \)
is a Casimir function for \\
\( \Pi (A\oplus A^{*}) \). Finally, property 2, 
\[
\rho (e_{1}\circ e_{2})=[\rho (e_{1}),\rho (e_{2})]\]
 when both sides are applied to an arbitrary \( f\in \cinf (M) \) translates
into 
\[
\{\{\{\theta ,e_{1}\},e_{2}\},\{\theta ,f\}\}=\{e_{1},\{\theta ,\{e_{2},\{\theta ,f\}\}\}\}-\{e_{2},\{\theta ,\{e_{1},\{\theta ,f\}\}\}\}\]
 for all \( f\in \cinf (M) \). Using Jacobi again, we have 
\[
\begin{array}{rcl}
-\{\{e_{2},\{e_{1},\theta \}\},\{\theta ,f\}\} & = & -\{e_{2},\{\{e_{1},\theta \},\{\theta ,f\}\}\}+\{\{e_{1},\theta \},\{e_{2},\{\theta ,f\}\}\}=\\
 & = & -\{e_{2},\{\theta ,\{e_{1},\{\theta ,f\}\}\}\}-\{e_{2},\{e_{1},\{\theta ,\{\theta ,f\}\}\}\}+\\
 & + & \{\theta ,\{e_{1},\{e_{2},\{\theta ,f\}\}\}\}+\{\{e_{1},\{\theta ,\{e_{2},\{\theta ,f\}\}\}\}=\\
 & = & \{\{e_{1},\{\theta ,\{e_{2},\{\theta ,f\}\}\}\}-\{e_{2},\{\theta ,\{e_{1},\{\theta ,f\}\}\}\}
\end{array}\]
 since \( \{\theta ,\{\theta ,f\}\}=0 \) by (\ref{eqn:master2}), while \( \{e_{1},\{e_{2},\{\theta ,f\}\}\}=0 \)
since \( \{e_{2},\{\theta ,f\}\}\in \cinf (M) \) is a Casimir function on \( \Pi (A\oplus A^{*}) \).

Thus all of the properties of a Courant algebroid are verified. Notice that
(\ref{eqn:master2}) was only needed to derive properties 1 and 2.
\end{proof}

\section{Quasi-bialgebroids\label{sec:quasibialg}}

The hamiltonian \( \theta  \) we constructed above was a sum of two terms,
\( \mu  \) of bidegree \( (1,2) \), and \( L^{*}\gamma  \) of bidegree \( (2,1) \),
so it has total degree \( \kappa (\theta )=3 \). There is nothing to prevent
us from adding a \( \phi  \) of bidegree \( (0,3) \) and/or a \( \psi  \)
of bidegree \( (3,0) \) to \( \theta  \), and require that \( \{\theta ,\theta \}=0 \). 

\begin{defn}
A \emph{proto-bialgebroid} is the supermanifold \( T^{*}\Pi A \) together with
a function \( \theta  \) such that \( \kappa (\theta )=3 \) and \( \{\theta ,\theta \}=0 \).
\end{defn}
Thus, a proto-bialgebroid structure consists of a vector field \( d_{A} \)
on \( \Pi A \), a vector field \( d_{A^{*}} \) on \( \Pi A^{*} \), and two
functions \( \phi \in \Gamma (\bigwedge ^{3}A^{*})\subset \cinf (\Pi A) \)
and \( \psi \in \Gamma (\bigwedge ^{3}A)\subset \cinf (\Pi A^{*}) \). Then
\( \theta =\mu +L^{*}\gamma +\pi ^{*}\phi +L^{*}\bar{\pi }^{*}\psi  \), and
the equation \( \{\theta ,\theta \}=0 \) splits according to the bigrading
into the following five equations:
\begin{equation}
\label{eqn:protobialg}
\begin{array}{c}
\half \{\mu ,\mu \}+\{L^{*}\gamma ,\pi ^{*}\phi \}=0\\
\{\mu ,L^{*}\gamma \}+\{\pi ^{*}\phi ,L^{*}\bar{\pi }^{*}\psi \}=0\\
\half L^{*}\{\gamma ,\gamma \}+\{\mu ,L^{*}\bar{\pi }^{*}\psi \}=0\\
\{\mu ,\pi ^{*}\phi \}=\{\gamma ,\bar{\pi }^{*}\psi \}=0
\end{array}
\end{equation}
 In particular, \( d_{A}\phi =d_{A^{*}}\psi =0 \) and the Schouten brackets
\( \br _{A} \) and \( \br _{A^{*}} \) are defined, but neither \( d_{A} \)
nor \( d_{A^{*}} \) square to zero, nor is \( d_{A} \) a derivation of \( \br _{A^{*}} \).
The defects in all cases are determined by the above relations.

Nevertheless, since the Poisson bracket on \( T^{*}\Pi A \) has total degree
\( -2 \), \( p^{*}\Gamma (A\oplus A^{*}) \) will be closed under both the
Poisson bracket and the derived bracket for any proto-bialgebroid, since elements
of \( p^{*}\Gamma (A\oplus A^{*}) \) have total degree 1. Thus, a slight modification
of Theorem \ref{thm:Courant=3DDerived} to include \( \phi  \) and \( \psi  \),
and repeating the argument of Theorem \ref{thm:Bialg=3D>Courant} yields 

\begin{thm}
Any proto-bialgebroid structure on \( T^{*}\Pi A \) induces a Courant algebroid
structure on the bundle \( A\oplus A^{*} \) given by 
\[
\begin{array}{rcl}
\bform {X+\xi }{Y+\eta } & = & \eta (X)+\xi (Y)\\
(X+\xi )\circ (Y+\eta ) & = & ([X,Y]_{A}+L^{A^{*}}_{\xi }Y-i_{\eta }d_{A^{*}}X-\psi (\xi ,\eta ))+\\
 & + & ([\xi ,\eta ]_{A^{*}}+L_{X}^{A}\eta -i_{Y}d_{A}\xi -\phi (X,Y))\\
\DD f & = & d_{A}f+d_{A^{*}}f
\end{array}\]
 where \( \phi  \) is viewed as a bundle map \( \phi :\bigwedge ^{2}A\rightarrow A^{*} \)
and likewise, \( \psi :\bigwedge ^{2}A^{*}\rightarrow A \).
\end{thm}
We will consider the special case where either \( \phi  \) or \( \psi  \)
is zero, say, \( \psi =0 \). The equations (\ref{eqn:protobialg}) reduce to
\[
\begin{array}{rcl}
\half \{\mu ,\mu \}+\{L^{*}\gamma ,\phi \} & = & 0\\
\{\gamma ,\gamma \} & = & 0\\
\{\mu ,L^{*}\gamma \} & = & 0\\
\{\mu ,\phi \} & = & 0
\end{array}\]
Deciphering these equations we arrive at 

\begin{defn}
A \emph{quasi-Lie bialgebroid} structure on \( (A,A^{*}) \) consists of the
following data: 
\begin{itemize}
\item A Lie algebroid structure on \( A^{*} \)
\item A bundle map \( a:A\rightarrow TM \)
\item A skew-symmetric operation \( \br _{A} \) on \( \Gamma (A) \) 
\item An element \( \phi \in \Gamma (\bigwedge ^{3}A^{*}) \) 
\end{itemize}
satisfying the following properties:
\begin{enumerate}
\item For all \( X,Y\in \Gamma (A), \) \( f\in \cinf (M) \), 
\[
[X,fY]_{A}=f[X,Y]_{A}+(a(X)f)Y\]
 
\item For all \( X,Y\in \Gamma (A) \), 
\[
a([X,Y]_{A})=[a(X),a(Y)]+a_{*}\phi (X,Y)\]
 where \( a_{*} \) is the anchor of the Lie algebroid \( A^{*} \) and \( \phi (X,Y)=i_{X\wedge Y}\phi \in \Gamma (A^{*}) \).
\item For all \( X,Y,Z\in \Gamma (A) \), 
\[
\begin{array}{l}
[[X,Y]_{A},Z]_{A}+[[Y,Z]_{A},X]_{A}+[[Z,X]_{A},Y]_{A} = d_{A^{*}}\phi (X,Y,Z)+\\
+\phi (d_{A^{*}}X,Y,Z)-\phi (X,d_{A^{*}}Y,Z)+\phi (X,Y,d_{A^{*}}Z)
\end{array}\]
 where \( d_{A^{*}} \) is the differential on \( \Gamma (\bigwedge A) \) coming
from the Lie algebroid structure on \( A^{*} \), and \( \phi  \) is viewed
as a bundle map \( \bigwedge ^{4}A\rightarrow A \).
\item \( d_{A}\phi =0 \) where \( d_{A} \) is the differential on \( \Gamma (\bigwedge A^{*}) \)
coming from the structure \( (a,\br _{A}) \) on \( A \).
\end{enumerate}
\end{defn}
Notice that this is completely analogous to Drinfeld's \emph{quasi-Lie bialgebras}
\cite{KS3}\emph{.} Property 3 above is to be interpreted as a homotopy Jacobi
identity for \( \br _{A} \). 

\begin{cor}
A quasi-Lie bialgebroid structure on \( (A,A^{*}) \) gives rise to a Courant
algebroid structure on \( A\oplus A^{*} \).
\end{cor}
Finally, we will look at an important special case of this, \emph{exact Courant
algebroids,} which were recently studied and classified by \v{S}evera \cite{Sev}.
A Courant algebroid \( E \) is called \emph{exact} if the sequence 
\[
0\longrightarrow T^{*}M\stackrel{{\rho ^{*}}}{\longrightarrow }E\stackrel{{\rho }}{\longrightarrow }TM\longrightarrow 0\]
 is exact, where the \emph{co-anchor} \( \rho ^{*}:T^{*}M\rightarrow E \) is
given by 
\[
\bform {\rho ^{*}\xi }{e}=\xi (\rho (e))\]
 for all \( \xi \in T^{*}M \), \( e\in E \) (By Property 4 of Definition \ref{def:quasi-algebroid2},
\( \rho \smalcirc \rho ^{*}=0 \) in any Courant algebroid). Then the image
of \( T^{*}M \) is a Dirac subbundle, and the restriction of \( \circ  \)
to its sections is identically zero, by Lemma \ref{lemma:ideal2}. One then chooses
a ``connection'' on \( E \), i.e. an isotropic splitting \( \sigma :TM\rightarrow E \)
of the above exact sequence. This is not a problem: once we have one isotropic
subbundle \( T^{*}M \), transversal isotropic subbundles are sections of a
bundle over \( M \) whose fiber is an open cell in the Grassmanian of isotropic
subspaces of half dimension in a pseudo-Euclidean space of signature zero; the
fiber is contractible (it is diffeomorphic to the linear space of skew-symmetric
matrices), so sections always exist. The connection \( \sigma  \) identifies
the pseudo-Euclidean vector bundle \( E \) with \( TM\oplus T^{*}M \) with
the canonical inner product. To compute the Courant bracket on \( E \) in this
identification, one looks at the difference 
\[
\sigma (X)\circ \sigma (Y)-\sigma ([X,Y])=\rho ^{*}\phi (X,Y)\]
 where \( X,Y \) are vector fields; this holds because \( \sigma  \) is a
splitting; moreover, using the properties of a Courant algebroid, one immediately
deduces that \( \phi  \) is \( \cinf (M) \)-linear and completely skew-symmetric,
i.e. comes from a 3-form \( \phi \in \Omega ^{3}(M) \), which it is appropriate
to call the ``curvature'' of \( \sigma  \). From the Leibniz-Jacobi identity
for \( \circ  \) (Property 1 of Definition \ref{def:quasi-algebroid2}) one
deduces the ``Bianchi identity'' 
\[
d\phi =0\]
 The Courant bracket becomes 
\[
(X+\xi )\circ (Y+\eta )=[X,Y]+L_{X}\eta -i_{Y}d\xi +\phi (X,Y)\]
 where \( X,Y\in \X (M) \), \( \xi ,\eta \in \Omega ^{1}(M) \) and \( \phi (X,Y)=i_{X\wedge Y}\phi \in \Omega ^{1}(M) \).
Thus, any exact Courant algebroid comes from a quasi-Lie bialgebroid which is
in fact the standard Lie bialgebroid \( (TM,T^{*}M) \) with an additional piece
of data, the closed 3-form \( \phi  \) which twists the standard Courant bracket
(\ref{eqn:standardCA}) on \( TM\oplus T^{*}M \). 

Once a connection \( \sigma  \) is chosen, any other one, \( \sigma ' \),
differs from \( \sigma  \) by the graph of a 2-form \( \omega  \); its curvature
\( \phi ' \) is related to \( \phi  \) simply by 
\[
\phi '=\phi +d\omega \]
 Therefore, the cohomology class \( c=[\phi ]\in H^{3}(M,\R ) \) is independent
of the choice of \( \sigma  \) and completely determines the Courant algebroid
structure on \( E \). It is thus appropriate to call \( c=c(E) \) the \emph{characteristic
class} of \( E \). This classification of exact Courant algebroids is due to
P. \v{S}evera \cite{Sev}. 

\begin{example}
Let \( G \) be a compact semisimple Lie group, with Lie algebra \( \g  \)
and the Killing form \( (\cdot ,\cdot ) \). Then Cartan's \emph{structure tensor
\[
\phi (X,Y,Z)=\frac{1}{12}([X,Y],Z)\]
} is the canonical bi-invariant 3-form on \( G \) that gives a non-trivial
twisting of the standard Courant algebroid structure on \( TG\oplus T^{*}G \).
This Courant algebroid plays a role in the recently developed theory of group-valued
momentum maps \cite{AMM} \cite{AMW}. 

It is also well-known that \( H^{3}(G,\R ) \), which is generated by \( [\phi ] \),
classifies Kac-Moody central extensions of the loop algebra \( L\g  \) \cite{ChPr}.
It is a very interesting question what the above Courant algebroid has to do
with affine Kac-Moody algebras.
\end{example}
\begin{rem}
As a final remark, we note that \v{S}evera's classification of exact Courant
algebroids is completely analogous to the well-known classification of central
extensions 
\[
0\longrightarrow \R \longrightarrow E\stackrel{a}{\longrightarrow }TM\longrightarrow 0\]
 of vector fields by functions. The exact sequence above is known as an \emph{Atiyah
sequence.} \( E \) is then a Lie algebroid, and the kernel of the anchor \( a \)
is the trivial one-dimensional vector bundle. Such Lie algebroids are classified
by \( H^{2}(M,\R ) \); if the characteristic class \( c(E) \) is integral,
the Atiyah sequence integrates to a principal \( U(1) \)-bundle \( P\rightarrow M \)
and \( c(E)=c_{1}(P) \) is the first Chern class of \( P \). We thus recover
the classification of complex line bundles on \( M \). 

Now, the meaning of the integrality of the characteristic class of an exact
Courant algebroid is still unknown. It is a very interesting question related
to the existence of a ``global'' object for a Courant algebroid, like the
principal \( U(1) \)-bundle above, or its gauge groupoid. This was posed as
an open problem in \cite{LWX1}, and there is as yet no solution. \v{S}evera
\cite{Sev} suggests that the answer should come from \emph{Dixmier-Douady gerbes},
but no global object for gerbes is known, either, nor is there a direct correspondence
between Courant algebroids and gerbes. Investigating these and related questions
are a logical continuation of this work.
\end{rem}

\chapter{Poisson Cohomology of \protect\( SU(2)\protect \)-covariant Poisson structures
on \protect\( S^{2}\protect \)\label{chapter:Poisson}}

\newcommand{\C}{\mathbb C}
  
\newcommand{\fr}[1]{\mathfrak {#1 }}
   In this chapter we shall compute the Poisson cohomology of the one-parameter
family of \( SU(2) \)-covariant Poisson structures on the homogeneous space
\( S^{2}=\mathbb CP^{1}=SU(2)/U(1) \), where \( SU(2) \) is endowed with its
standard Poisson-Lie group structure, thus extending the result of Ginzburg
\cite{Gin1} on the Bruhat-Poisson structure which is a member of this family.
As a corollary of our computation, we deduce that these structures are nontrivial
deformations of each other in the direction of the standard rotation-invariant
symplectic structure on \( S^{2} \); another corollary is that these structures
do not admit rescaling.

\section{Poisson-Lie groups and Poisson actions}

Here we briefly recall some basic notions of the theory of Poisson-Lie groups
that we will need. For more details the interested reader should consult \cite{ChPr},
\cite{KorSoi}, or \cite{LuWe}.

\begin{notation}
Let a Lie group \( G \) act on a manifold \( P. \) Then each \( g\in G \)
gives rise to a map \( P\to P \) given by \( p\mapsto gp \). We shall denote
this map as well as its derivatives and their tensor products by the same letter
\( g \) where it does not cause confusion. Likewise, every \( p\in P \) induces
a map \( G\to P \) by \( g\mapsto gp \) which, along with its derivatives,
we shall denote by \( p \) written on the right of the argument. This will
make our notation a lot less cumbersome.
\end{notation}
\begin{defn}
A Poisson structure \( \pi  \) on a Lie group \( G \) is called \emph{multiplicative}
if the group multiplication 
\[
m:G\times G\longrightarrow G\]
 is a Poisson map, where \( G\times G \) is equipped with the product Poisson
structure. The pair \( (G,\pi ) \) is then called a \emph{Poisson-Lie group.} 
\end{defn}
One checks that the multiplicativity condition is equivalent to the identity
\begin{equation}
\label{PLgp}
\pi (gh)=g\pi (h)+\pi (g)h\; \; \forall g,h\in G
\end{equation}
 In particular, one has \( \pi (e)=0, \) so the linearization (intrinsic derivative)
of \( \pi  \) at \( e \) gives a well-defined \emph{cobracket} \( \sigma :\fr{g}\rightarrow \fr{g}\wedge \fr{g} \)
by 
\[
\sigma (X)=(L_{X_{l}}\pi )(e)=\left. \frac{d}{dt}\right| _{t=0}\pi (\exp (tX))\exp (-tX),\]
 where \( X_{l} \) denotes the left-invariant vector field corresponding to
\( X\in \fr{g}. \) The multiplicativity of \( \pi  \) (\ref{PLgp}) then implies
the cocycle property of \( \sigma  \): 
\[
\sigma ([X,Y])=[X,\sigma (Y)]-[Y,\sigma (X)]\]
 On the other hand, the Jacobi identity for \( \pi  \) implies that the adjoint
of \( \sigma , \) \\
\( \sigma ^{*}:\fr{g}^{*}\wedge \fr{g}^{*}\rightarrow \fr{g}^{*} \) also satisfies
Jacobi, i.e. defines a Lie bracket on \( \fr{g}^{*}. \) Thus, \( (\g ,\br ,\sigma ) \)
is a Lie bialgebra (see Section \ref{sec:bialgebras}) called the \emph{tangent
Lie bialgebra} of the Poisson-Lie group \( G \)\emph{.} It can be shown \cite{LuWe}
that if \( G \) is connected, \( \pi  \) is uniquely determined by \( \sigma . \) 

It may happen that the cocycle \( \sigma  \) is a coboundary , that is, there
exists an \( {\bf r}\in \fr{g}\wedge \fr{g} \) such that \( \sigma (X)=-[X,{\bf r}] \)
(always the case if \( \fr{g} \) is semisimple). Such an \( {\bf r} \) is
called a \emph{classical r-matrix.} The Jacobi identity for \( \sigma ^{*} \)
is equivalent to the condition on \( {\bf r} \) that \( [{\bf r},{\bf r}]\in \bigwedge ^{3}\fr{g} \)
be \( ad \) - invariant (the so-called Modified Classical Yang-Baxter Equation).
Here \( \br  \) is the algebraic Schouten bracket of the Lie algebra \( \g  \)
(Example \ref{eg:AlgSchouten}). The multiplicative Poisson structure \( \pi  \)
is given in terms of \( {\bf r} \) by \emph{}
\begin{equation}
\label{SklBr}
\pi (g)={\bf r}g-g{\bf r},
\end{equation}
 and the corresponding Poisson bracket on \( G \) is called the \emph{Sklyanin
bracket.}

\begin{defn}
Let a Poisson-Lie group \( (G,\pi _{G}) \) act on a manifold \( P. \) We say
that a Poisson structure \( \pi _{P} \) on \( P \) is \( G \)\emph{-covariant}
if the action map 
\[
\rho :G\times P\longrightarrow P\]
 is Poisson, where \( G\times P \) is equipped with the product Poisson structure.
The action \( \rho  \) is then called a \emph{Poisson action.} If \( \rho  \)
is transitive, \( (P,\pi _{P}) \) is called a \emph{Poisson homogeneous space.}
\end{defn}
The covariance condition is equivalent to the identity 
\begin{equation}
\label{PAction}
\pi _{P}(gp)=\pi _{G}(g)p+g\pi _{P}(p)\; \; \forall g\in G,\: p\in P
\end{equation}
 Note that \( G \) does \emph{not} act by Poisson transformations unless \( \pi _{G}=0 \). 

\begin{fact}
\cite{KorSoi} If \( (G,\pi _{G}) \) is a Poisson-Lie group, \( H\subset G \)
a Poisson (or even coisotropic) subgroup, then there is a unique Poisson structure
\( \pi _{P} \) on \( P=G/H \) making the canonical projection a Poisson map.
Moreover, \( \pi _{P} \) is \( G \)-covariant. 
\end{fact}
So if \( H \) is coisotropic, \( G/H \) is always a Poisson homogeneous space;
however, the projection of \( \pi _{G} \) is in general not the only \( G \)-covariant
Poisson structure on \( G/H \): adding any \( G \)-invariant bivector field
will give another one provided that the sum satisfies the Jacobi identity.

\section{Description of the Poisson structures}

\subsection{The classical r-matrix and the standard Poisson-Lie structure on \protect\( SU(2)\protect \).}

The constructions below can be carried out for any compact semisimple Lie group,
but we will only consider \( SU(2) \).

Recall that the Lie algebra \( \fr{su}(2) \) of \( 2\times 2 \) skew-hermitian
traceless matrices has a basis 
\[
e_{1}=\half \left( \begin{array}{cc}
i & 0\\
0 & -i
\end{array}\right) ,\; \; \; e_{2}=\half \left( \begin{array}{cc}
0 & 1\\
-1 & 0
\end{array}\right) ,\; \; \; e_{3}=\half \left( \begin{array}{cc}
0 & i\\
i & 0
\end{array}\right) \]
 with the commutation relations \( [e_{\alpha },e_{\beta }]=\epsilon _{\alpha \beta \gamma }e_{\gamma } \),
where \( \epsilon _{\alpha \beta \gamma } \) is the completely skew-symmetric
symbol. The span of \( e_{1} \) is the Cartan subalgebra of \( \mathfrak {a}\in \mathfrak {su}(2) \).
Recall also that 
\[
SU(2)=\left\{ \left. U=\left( \begin{array}{cc}
u & -\bar{v}\\
v & \bar{u}
\end{array}\right) \right| u,v\in \mathbb C,\; \; \; \det U=u\bar{u}+v\bar{v}=1\right\} \]
 identifies \( SU(2) \) with the unit sphere in \( \mathbb C^{2} \). The \emph{standard
r-matrix} \( {\bf r}=e_{2}\wedge e_{3}\in \fr{{su}(2)\wedge \fr{{su}(2)}} \)
defines a multiplicative Poisson structure on \( SU(2) \) by 
\begin{equation}
\label{PoissonSU(2)}
\pi _{SU(2)}(U)={\bf r}U-U{\bf r}
\end{equation}
In coordinates, 
\[
\pi\left( \left( \begin{array}{cc}
u & -\bar{v}\\
v & \bar{u}
\end{array}\right) \right) =\frac{1}{4}\left( \left( \begin{array}{cc}
v & \bar{u}\\
-u & \bar{v}
\end{array}\right) \wedge \left( \begin{array}{cc}
iv & i\bar{u}\\
iu & -i\bar{v}
\end{array}\right) -\left( \begin{array}{cc}
\bar{v} & u\\
-\bar{u} & v
\end{array}\right) \wedge \left( \begin{array}{cc}
-i\bar{v} & iu\\
i\bar{u} & iv
\end{array}\right) \right) =\]
 
\begin{equation}
\label{PoissonSU(2)coord}
=-iv\bar{v}\der{}{u}\wedge \der{}{\bar{u}}+\half \left( iuv\der{}{u}\wedge \der{}{v}+\overline{iuv\der{}{u}\wedge \der{}{v}}\right) +\half \left( iu\bar{v}\der{}{u}\wedge \der{}{\bar{v}}+\overline{iu\bar{v}\der{}{u}\wedge \der{}{\bar{v}}}\right) 
\end{equation}
 The Poisson brackets are 
\[
\begin{array}{cccc}
\{u,\bar{u}\}=-iv\bar{v}, & \{u,v\}=\half iuv, & \{u,\bar{v}\}=\half iu\bar{v}, & \{v,\bar{v}\}=0
\end{array}\]
 It is easy to see that these formulas in fact define a smooth real Poisson
structure on all of \( \C ^{2} \) that restricts to the unit sphere.

\subsection{The Bruhat-Poisson structure on \protect\( \C P^{1}\protect \).}

The r-matrix is invariant under the action of the Cartan subalgebra \( \fr{a} \)
, since 
\[
[e_{1},{\bf r}]=[e_{1},e_{2}\wedge e_{3}]=[e_{1},e_{2}]\wedge e_{3}-e_{2}\wedge [e_{1},e_{3}]=e_{3}\wedge e_{3}+e_{2}\wedge e_{2}=0\]
 Hence, the Poisson tensor (\ref{PoissonSU(2)}) vanishes on the maximal torus
(the diagonal subgroup) \( A=U(1)\subset SU(2) \). In particular, \( U(1) \)
is a Poisson subgroup, and hence \( \pi _{SU(2)} \) descends to the quotient
\( SU(2)/U(1)=S^{3}/S^{1}=(\C ^{2}\setminus 0)/\C ^{\times }=\C P^{1}=S^{2} \).
The resulting Poisson structure \( \pi _{1} \) on \( \C P^{1} \) is called
the \emph{Bruhat-Poisson structure} because its symplectic leaves coincide with
the Bruhat cells in \( \C P^{1} \) \cite{LuWe}: the base point where \( \pi _{1} \)
vanishes, and the complementary open cell where \( \pi _{1} \) is invertible.
It is \( SU(2) \)-covariant since \( \pi _{SU(2)} \) is multiplicative. It
is an easy calculation to deduce from (\ref{PoissonSU(2)coord}) that in the
inhomogeneous coordinate chart \( w=v/u \) covering the base point \( \pi _{1} \)
is given by 
\[
\pi _{1}=-iw\bar{w}(1+w\bar{w})\der{}{w}\wedge \der{}{\bar{w}}\]
 In particular, it has a quadratic singularity at \( w=0 \). The other inhomogeneous
chart \( z=u/v=1/w \) gives coordinates on the open symplectic leaf, in which
\[
\pi _{1}=-i(1+z\bar{z})\der{}{z}\wedge \der{}{\bar{z}}\]
 The corresponding symplectic 2-form is 
\[
\omega _{1}=\frac{idz\wedge d\bar{z}}{1+z\bar{z}}\]
 Notice that this symplectic leaf has infinite volume.

\subsection{The other \protect\( SU(2)\protect \)-covariant Poisson structures on \protect\( S^{2}\protect \).}

The difference between any two \( SU(2) \) - covariant Poisson structures on
\( \C P^{1} \) is an \( SU(2) \) - invariant bivector field (by (\ref{PAction}))
which is Poisson because in two dimensions, any bivector field is. Thus, any
covariant structure is obtained by adding an invariant structure to the Bruhat
structure \( \pi _{1}. \) To see what these structures look like, it is convenient
to embed the Riemann sphere \( \C P^{1} \)as the unit sphere \( S^{2}\subset \R ^{3} \)
by the (inverse of) the stereographic projection. The coordinate transformations
are given by 
\[
\begin{array}{ccr}
x_{1}=\frac{2x}{1+x^{2}+y^{2}} & \; \; \; \;  & x=\frac{x_{1}}{1-x_{3}}\\
x_{2}=\frac{2y}{1+x^{2}+y^{2}} & \; \; \; \;  & y=\frac{x_{2}}{1-x_{3}}\\
x_{3}=\frac{x^{2}+y^{2}-1}{1+x^{2}+y^{2}} & \; \; \; \;  & x^{2}+y^{2}=\frac{1+x_{3}}{1-x_{3}}
\end{array}\]
where \( z=x+iy \). We shall identify \( \R ^{3} \) with \( \fr{su}(2)^{*} \),
with the coadjoint action of \( SU(2) \) by rotations. Then the linear Poisson
structure on \( \R ^{3}=\fr{su}(2)^{*} \) is given by 
\[
-\pi =x_{1}\der{}{x_{2}}\wedge \der{}{x_{3}}+x_{2}\der{}{x_{3}}\wedge \der{}{x_{1}}+x_{3}\der{}{x_{1}}\wedge \der{}{x_{2}}\]
 whose restriction to the unit sphere (a coadjoint orbit), also denoted by \( -\pi  \),
is \( SU(2) \) - invariant and symplectic. Moreover, up to a constant multiple,
\( \pi  \) is the only rotation-invariant Poisson structure on \( S^{2} \):
any other invariant structure is of the form \( \pi '=f\pi  \) for some function
\( f \), but since both \( \pi  \) and \( \pi ' \) are invariant, so is \( f \),
hence \( f \) is a constant. It follows that there is a one-parameter family
of \( SU(2) \) - covariant Poisson structures of the form \( \pi '=\pi _{1}+\alpha \pi  \),
\( \alpha \in \R  \); since \( \pi _{1}=(1-x_{3})\pi  \) (straightforward
calculation), all \( SU(2) \) - covariant structures are of the form 
\[
\pi _{c}=\pi _{1}+(c-1)\pi =(c-x_{3})\pi ,\; \; \; c\in \R \]
 It follows that \( \pi _{c} \) is symplectic for \( |c|>1 \), Bruhat for
\( c=\pm 1 \), while for \( |c|<1 \) \( \pi _{c} \) vanishes on the circle
\( \{x_{3}=c\} \) and is nonsingular elsewhere; \( \pi _{c} \) thus has two
open symplectic leaves (hemispheres) and a ``necklace'' of zero-dimensional
symplectic leaves along the circle. It is these ``necklace'' structures whose
Poisson cohomology we shall compute. Notice that \( \pi _{c} \) and \( \pi _{-c} \)
are isomorphic as Poisson manifolds via \( x_{3}\mapsto -x_{3} \). 

In the original \( \{w,\bar{w}\} \) - coordinates we have 
\begin{equation}
\label{PiStandard}
\pi =-\frac{i}{2}(1+w\bar{w})^{2}\der{}{w}\wedge \der{}{\bar{w}}=\frac{1}{4}(1+x^{2}+y^{2})^{2}\der{}{x}\wedge \der{}{y},
\end{equation}

\begin{eqnarray}
\pi _{c} & =\pi _{1}+(c-1)\pi = & -\frac{i}{2}(1+w\bar{w})((c+1)w\bar{w}+c-1)\der{}{w}\wedge \der{}{\bar{w}}=\nonumber \\
 &  & =\frac{1}{4}(1+x^{2}+y^{2})((c+1)(x^{2}+y^{2})+c-1)\der{}{x}\wedge \der{}{y}\label{PiC} 
\end{eqnarray}
 where \( w=x+iy \).

\subsection{Symplectic areas and modular vector fields.}

Before we proceed to cohomology computations, we shall compute some invariants
of the structures \( \pi _{c} \). For \( |c|>1 \) \( \pi _{c} \) is symplectic,
and the only invariant is the symplectic area. For the other values of \( c \),
the areas of the open symplectic leaves are easily seen to be infinite; instead,
we will compute the modular vector field of \( \pi _{c} \) with respect to
the standard rotation-invariant volume form \( \omega  \) on \( S^{2} \) (the
inverse of \( \pi  \)). By elementary calculations we obtain the following

\begin{lem}
(1) If \( |c|>1, \) the symplectic volume of \( (S^{2},\pi _{c}) \) is given
by 
\[
V(c)=2\pi \ln \frac{c+1}{c-1}\]

(2) For all values of \( c \) the modular vector field with respect to \( \omega  \)
is 
\[
\Delta _{\omega }=x\der{}{y}-y\der{}{x}\]
 
\end{lem}
\begin{cor}
If \( |c|,\: |c'|>1, \) \( \pi _{c} \) and \( \pi _{c'} \) are not isomorphic
unless \( |c|=|c'| \).
\begin{cor}
\label{Cor:ModClass}If \( |c|<1, \) the modular class of \( \pi _{c} \) is
nonzero.
\end{cor}
\end{cor}
\begin{proof}
The modular vector field \( \Delta _{\omega } \) rotates the necklace, hence
cannot be Hamiltonian.
\end{proof}
In fact, the modular class of the Bruhat-Poisson structures \( \pi _{\pm 1} \)
is also nonzero \cite{Gin1}. 

Unfortunately, the modular vector field does not help us distinguish the different
``necklace'' structures. The restriction of \( \Delta _{\omega } \) to the
necklace is independent of \( \omega  \) since changing \( \omega  \) changes
\( \Delta _{\omega } \) by a Hamiltonian vector field which necessarily vanishes
along the necklace, so the period of \( \Delta _{\omega } \) restricted to
the necklace is an invariant, but it has the same value of \( 2\pi  \) for
all \( \pi _{c}. \) When we compute the Poisson cohomology of \( \pi _{c} \)
we will see a different way to distinguish them.

\section{Computation of Poisson cohomology}

For \( |c|>1 \) \( \pi _{c} \) is symplectic, so its Poisson cohomology is
isomorphic to the deRham cohomology of \( S^{2} \); the Poisson cohomology
of the Bruhat-Poisson structure \( \pi _{\pm 1} \) was worked out by Ginzburg
\cite{Gin1}. Here we shall compute the cohomology of the necklace structures
\( \pi _{c} \) for \( |c|<1 \). Our strategy will be similar to Ginzburg's:
first compute the cohomology of the formal neighborhood of the necklace, show
that the result is actually valid in a finite small neighborhood and finally,
use a Mayer-Vietoris argument to deduce the global result. The validity of the
Mayer-Vietoris argument for Poisson cohomology comes from the simple observation
that on any Poisson manifold \( (P,\pi ) \) the differential \( d_{\pi } \)
is functorial with respect to restrictions to open subsets (i.e. a morphism
of the sheaves of smooth multivector fields on \( P \)). 

It will be convenient to introduce another change of coordinates: 
\[
\begin{array}{ccc}
s=\frac{x}{\sqrt{1+x^{2}+y^{2}}} & \; \; \;  & t=\frac{y}{\sqrt{1+x^{2}+y^{2}}}
\end{array}\]
 mapping the \( (x,y) \)-plane to the open unit disk in the \( (s,t) \)-plane.
In the new coordinates \( \pi _{c} \) and \( \pi  \) are given by
\begin{equation}
\label{PiC'}
\pi _{c}=\half (s^{2}+t^{2}-\frac{1-c}{2})\der{}{s}\wedge \der{}{t}
\end{equation}
 
\begin{equation}
\label{PiStandard'}
\pi =\frac{1}{4}\der{}{s}\wedge \der{}{t}
\end{equation}
and the necklace is the circle of radius \( R=\sqrt{\frac{1-c}{2}} \). Observe
that rescaling \( s=\alpha s',\; t=\alpha t' \) (\( \alpha >0 \)) takes \( \pi _{c} \)
with necklace radius \( R \) to \( \pi _{c'} \) with necklace radius \( R'=R/\alpha  \).
But this is only a local isomorphism: it does not extend to all of \( S^{2} \)
since it is not a diffeomorphism of the unit disk. In any case, it shows that
all necklace structures are locally isomorphic, so for local computations we
may assume that \( \pi _{c} \) is given in suitable coordinates by 
\[
\pi _{c}=\half (s^{2}+t^{2}-1)\der{}{s}\wedge \der{}{t}\]

\subsection{Cohomology of the formal neighborhood of the necklace.}

Since \( \pi _{c} \) is rotation-invariant, we can lift the computations in
the formal neighborhood of the unit circle in the \( (s,t) \)-plane to its
universal cover by introducing ``action-angle coordinates'' \( (I,\theta ) \):
\[
s=\sqrt{1+I}\cos \theta \; \; \; t=\sqrt{1+I}\sin \theta \]
 in which \( \pi _{c} \) is linear: 
\[
\pi _{c}=I\der{}{I}\wedge \der{}{\theta }\]
 Of course we will have to restrict attention to multivector fields whose coefficients
are periodic in \( \theta  \). It will be convenient to think of multivector
fields as functions on the supermanifold with coordinates \( (I,\theta ,\xi ,\eta ) \)
where \( \xi  \)''=''\( \partial _{I} \) and \( \eta  \)''=''\( \partial _{\theta } \)
are Grassmann (anticommuting) variables. Then \( \pi _{c}=I\xi \eta  \) is
a function and 
\[
d_{\pi _{c}}=[\pi _{c},\cdot ]=-I\eta \der{}{I}+I\xi \der{}{\theta }-\xi \eta \der{}{\xi }\]
 is a (homological) vector field. Since \( d_{\pi _{c}} \) commutes with rotations,
we can split the complex into Fourier modes
\[
\begin{array}{ccc}
\fr{X}_{n}^{0}=\{f(I)e^{in\theta }\};\;  & \fr{X}_{n}^{1}=\{(f(I)\xi +g(I)\eta )e^{in\theta }\};\;  & \fr{X}_{n}^{2}=\{h(I)\xi \eta e^{in\theta }\},
\end{array}\]
 where \( f(I) \), \( g(I) \) and \( h(I) \) are formal power series in \( I \).
It will be convenient to treat the zero and non-zero modes separately; it will
turn out that the cohomology is concentrated entirely in the zero mode.

\begin{case}
\textbf{The zero mode} (\textbf{\( n=0 \)}) \textbf{}consists of multivector
fields independent of \( \theta  \), so \( d_{\pi _{c}} \) becomes 
\[
\left. d_{\pi _{c}}\right| _{\fr{X}_{0}}=-I\eta \der{}{I}+\eta \xi \der{}{\xi }\]
 which preserves the degree in \( I \) so the complex \( \fr{X}_{0} \) splits
further into a direct product of sub-complexes \( \fr{X}_{0,m} \), \( m\geq 0 \)
according to the degree:
\[
0\rightarrow \fr{X}^{0}_{0,m}\rightarrow \fr{X}^{1}_{0,m}\rightarrow \fr{X}^{2}_{0,m}\rightarrow 0\]
 These complexes are very small (\( \fr{X}_{0,m}^{0} \) and \( \fr{X}_{0,m}^{2} \)
are one-dimensional, while \( \fr{X}_{0,m}^{2} \) is two-dimensional) and their
cohomology is easy to compute. For \( f=cI^{m}\in \fr{X}_{0,m}^{0} \), \( d_{\pi _{c}}f=-cmI^{m}\eta  \),
while for \( X=aI^{m}\xi +bI^{m}\eta \in \fr{X}_{0,m}^{1} \), \( d_{\pi _{c}}X=a(m-1)I^{m}\xi \eta  \).
Therefore, it is clear that for \( m>1 \) the complex is acyclic. On the other
hand, the cohomology of \( \fr{X}_{0,0} \) is generated by \( 1\in \fr{X}_{0,0}^{0} \)
and \( \eta \in \fr{X}_{0,0}^{1} \), while the cohomology of \( \fr{X}_{0,1} \)
is generated by \( I\xi \in \fr{X}_{0,1}^{1} \) and \( I\xi \eta \in \fr{X}_{0,1}^{2} \).
Putting these together we obtain 
\begin{equation}
\label{InvHom}
\begin{array}{ccccl}
H_{0}^{0} & = & \R  & = & \textrm{span}\{1\}\\
H_{0}^{1} & = & \R ^{2} & = & \textrm{span}\{\partial _{\theta },I\partial _{I}\}\\
H_{0}^{2} & = & \R  & = & \textrm{span}\{I\partial _{I}\wedge \partial _{\theta }\}
\end{array}
\end{equation}
 
\begin{case}
\textbf{The non-zero modes \( (n\neq 0) \).} In this case \( d_{\pi _{c}} \)
does not preserve the \( I \)-grading so we'll have to consider all power series
at once. Let 
\[
\begin{array}{cclcc}
f & = & (\sum ^{\infty }_{m=0}f_{m}I^{m})e^{in\theta } & \in  & \fr{X}_{n}^{0}\\
X & = & (\sum ^{\infty }_{m=0}a_{m}I^{m})e^{in\theta }\xi +(\sum ^{\infty }_{m=0}b_{m}I^{m})e^{in\theta }\eta  & \in  & \fr{X}_{n}^{1}\\
B & = & (\sum ^{\infty }_{m=0}c_{m}I^{m})e^{in\theta }\xi \eta  & \in  & \fr{X}_{n}^{2}
\end{array}\]
 Then 
\[
\begin{array}{ccl}
d_{\pi _{c}}f & = & (\sum ^{\infty }_{m=1}inf_{m-1}I^{m})e^{in\theta }\xi +(\sum ^{\infty }_{m=1}mf_{m}I^{m})e^{in\theta }\eta \\
d_{\pi _{c}}X & = & (-a_{0}+\sum ^{\infty }_{m=1}((m-1)a_{m}+inb_{m-1})I^{m})e^{in\theta }\xi \eta 
\end{array}\]
 (and, of course, \( d_{\pi _{c}}B=0 \)). We see immediately that \( d_{\pi _{c}}f=0\Leftrightarrow f=0 \),
hence \( H_{n}^{0}=\{0\} \). Moreover, any \( B \) is a coboundary: 
\[
B=d_{\pi _{c}}\left( (\sum ^{\infty }_{m\neq 1}\frac{c_{m}}{m-1}I^{m})e^{in\theta }\xi +\frac{c_{1}}{in}e^{in\theta }\eta \right) \]
 so \( H_{n}^{0}=\{0\} \) as well. Now, \( X \) is a cocycle if and only if
\[
\begin{array}{cccl}
a_{0} & = & b_{0} & =0\\
b_{m} & = & -\frac{ma_{m+1}}{in}, & m\geq 1
\end{array}\]
 Let \( f_{m}=\frac{a_{m+1}}{in} \) for \( m\geq 0 \), \( f=\sum f_{m}I^{m} \).
Then \( X=d_{\pi _{c}}f \). Hence \( H_{n}^{1} \) is also trivial. So for
\( n\neq 0 \) \( \fr{X}_{n} \) is acyclic.
\end{case}
\end{case}
It follows that the Poisson cohomology of the formal neighborhood of the necklace
is as in (\ref{InvHom}).

\subsection{Justification for the smooth case.}

To see that the cohomology of a finite small neighborhood of the necklace is
the same as for the formal neighborhood we apply an argument similar to Ginzburg's
\cite{Gin1}. For each Fourier mode consider the following exact sequence of
complexes: 
\[
0\rightarrow \fr{X}^{\star }_{n,\textrm{flat }}\rightarrow \fr{X}^{\star }_{n,\textrm{smooth }}\rightarrow \fr{X}^{\star }_{n,\textrm{formal }}\rightarrow 0\]
 where \( \fr{X}^{\star }_{n,\textrm{flat }} \) consists of smooth multivector
fields whose coefficients vanish along the necklace together with all derivatives.
This sequence is exact by a theorem of E. Borel. It suffices to show that the
flat complex is acyclic. But \( \pi _{c}^{\#}:\fr{X}^{\star }_{n,\textrm{flat }}\rightarrow \Omega ^{\star }_{n,\textrm{flat}} \)
is an isomorphism since the coefficient of \( \pi _{c} \) is a polynomial in
\( I \), and every flat form can be divided by a polynomial with a flat result.
Furthermore, the flat deRham complex is acyclic by the homotopy invariance of
deRham cohomology. 

Finally, we observe that a smooth multivector field in a neighborhood of the
necklace (given by a \emph{convergent} Fourier series) is a coboundary if and
only if each mode is, and the primitives can be chosen so that the resulting
series converges, as can be seen from the calculations in the previous subsection
(integration can only improve convergence). Therefore, the Poisson cohomology
of an annular neighborhood \( U \) of the necklace is 
\begin{equation}
\label{LocHom}
\begin{array}{ccccl}
H_{\pi _{c}}^{0}(U) & = & \R  & = & \textrm{span}\{1\}\\
H_{\pi _{c}}^{1}(U) & = & \R ^{2} & = & \textrm{span}\{\partial _{\theta },I\partial _{I}\}\\
H_{\pi _{c}}^{2}(U) & = & \R  & = & \textrm{span}\{I\partial _{I}\wedge \partial _{\theta }\}
\end{array}
\end{equation}
 Notice that the generators of \( H_{\pi _{c}}^{1}(U) \) are the rotation \( \partial _{\theta }=s\partial _{t}-t\partial _{s} \)
(the modular vector field) and the dilation \( I\partial _{I}=\frac{s^{2}+t^{2}-1}{2(s^{2}+t^{2})}(s\partial _{s}+t\partial _{t}) \),
while the generator of \( H^{2}_{\pi _{c}}(U) \) is \( \pi _{c} \) itself,
so in particular \( \pi _{c} \) does not admit rescalings even locally.

\subsection{From local to global cohomology.}

We now have all we need to compute the Poisson cohomology of a necklace Poisson
structure \( \pi _{c} \) on \( S^{2} \). Cover \( S^{2} \) by two open sets
\( U \) and \( V \) where \( U \) is an annular neighborhood of the necklace
as above, and \( V \) is the complement of the necklace consisting of two disjoint
open hemispheres on each of which \( \pi _{c} \) is nonsingular, so that the
Poisson cohomology of \( V \) and \( U\cap V \) is isomorphic to the deRham
cohomology. The short exact Mayer-Vietoris sequence associated to this cover
\[
0\rightarrow \fr{X}^{\star }(S^{2})\rightarrow \fr{X}^{\star }(U)\oplus \fr{X}^{\star }(V)\rightarrow \fr{X}^{\star }(U\cap V)\rightarrow 0\]
 leads to a long exact sequence in cohomology: 
\[
\begin{array}{ccccccccc}
0 & \rightarrow  & H^{0}_{\pi _{c}}(S^{2}) & \rightarrow  & H_{\pi _{c}}^{0}(U)\oplus H^{0}_{\pi _{c}}(V) & \rightarrow  & H^{0}_{\pi _{c}}(U\cap V) & \rightarrow  & \\
 & \rightarrow  & H_{\pi _{c}}^{1}(S^{2}) & \rightarrow  & H_{\pi _{c}}^{1}(U)\oplus H^{1}_{\pi _{c}}(V) & \rightarrow  & H^{1}_{\pi _{c}}(U\cap V) & \rightarrow  & \\
 & \rightarrow  & H_{\pi _{c}}^{2}(S^{2}) & \rightarrow  & H_{\pi _{c}}^{2}(U)\oplus H^{2}_{\pi _{c}}(V) & \rightarrow  & H^{2}_{\pi _{c}}(U\cap V) & \rightarrow  & 0
\end{array}\]
 Now, the first row is clearly exact since a Casimir function on \( S^{2} \)
must be constant on each of the two open symplectic leaves comprising \( V \),
hence constant on all of \( S^{2} \) by continuity. On the other hand, \( H^{1}_{\pi _{c}}(V)=H^{2}_{\pi _{c}}(V)=H^{2}_{\pi _{c}}(U\cap V)=\{0\} \).
Combining this with (\ref{LocHom}), we see that what we have left is  
\[
\begin{array}{cccccccc}
 &  &  &  & \R ^{2} &  & \R ^{2} & \\
 &  &  &  & \Vert  &  & \Vert  & \\
0 & \rightarrow  & H_{\pi _{c}}^{1}(S^{2}) & \rightarrow  & H_{\pi _{c}}^{1}(U)\oplus H^{1}_{\pi _{c}}(V) & \rightarrow  & H^{1}_{\pi _{c}}(U\cap V) & \rightarrow \\
 &  &  &  &  &  &  & \\
 & \rightarrow  & H_{\pi _{c}}^{2}(S^{2}) & \rightarrow  & H_{\pi _{c}}^{2}(U)\oplus H^{2}_{\pi _{c}}(V) & \rightarrow  & 0 & \\
 &  &  &  & \Vert  &  &  & \\
 &  &  &  & \R  &  &  & 
\end{array}\]
Now, on the one hand, we know by Corollary \ref{Cor:ModClass} that \( H_{\pi _{c}}^{1}(S^{2}) \)
is at least one-dimensional; on the other hand, the restriction of the dilation
vector field \( I\partial _{I} \) to \( U\cap V \) is not Hamiltonian: it
corresponds under \( \pi _{c}^{\#} \) to the generator of the first deRham
cohomology of the annulus diagonally embedded into \( U\cap V \) (a disjoint
union of two annuli). It follows that \( H_{\pi _{c}}^{1}(S^{2}) \) is exactly
one-dimensional, while \( H_{\pi _{c}}^{2}(S^{2}) \) is two-dimensional. 

It only remains to identify the generators. \( H_{\pi _{c}}^{1}(S^{2}) \) is
generated by the modular class, while one of the generators of \( H_{\pi _{c}}^{2}(S^{2}) \)
is \( \pi _{c} \) itself, since its class was shown to be nontrivial even locally.
The other generator is the image of \( (I\partial _{I},-I\partial _{I})\in H_{\pi _{c}}^{1}(U\cap V) \)
under the connecting homomorphism. This is somewhat unwieldy since it involves
a partition of unity subordinate to the cover \( \{U,V\} \) which does not
yield a clear geometric interpretation of the generator. Instead, we will show
directly that the standard rotationally invariant symplectic Poisson structure
\( \pi  \) on \( S^{2} \) is nontrivial in \( H_{\pi _{c}}^{2}(S^{2}) \)
and so can be taken as the second generator.

\begin{lem}
The class of the standard \( SU(2) \)-invariant Poisson structure \( \pi  \)
on \( S^{2} \) is nonzero in \( H_{\pi _{c}}^{2}(S^{2}) \).
\end{lem}
\begin{proof}
We will work in coordinates \( (s,t) \) on the unit disk in which \( \pi  \)
and \( \pi _{c} \) are given, respectively by (\ref{PiStandard'}) and(\ref{PiC'}).
Locally \( \pi  \) is a coboundary whose primitive is given by an Euler vector
field \( E=\frac{1}{2(c-1)}(s\partial _{s}+t\partial _{t}) \): it's easy to
check that \( [\pi _{c},E]=\pi  \). But \( E \) does not extend to a vector
field on \( S^{2} \) since it does not behave well ``at infinity'', i.e on
the unit circle in the \( (s,t) \)-plane. Therefore, to prove that \( \pi  \)
is globally nontrivial it suffices to show that there does not exist a Poisson
vector field \( X \) such that \( E+X \) is tangent to the unit circle and
the restriction is rotationally invariant. In fact, it suffices to show that
there is no Hamiltonian vector field \( X_{f} \) such that \( E+X_{f} \) vanishes
on the unit circle (since we can always add a multiple of the modular vector
field to cancel the rotation). Assuming that such an \( f \) exists, we will
have, in the polar coordinates \( s=r\cos \phi  \), \( t=r\sin \phi  \) :
\[
E+X_{f}=\frac{1}{2(c-1)}r\der{}{r}+\frac{1}{2r}(r^{2}-\frac{1-c}{2})\left( \der{f}{\phi }\der{}{r}-\der{f}{r}\der{}{\phi }\right) \]
 Upon restriction to \( r=1 \) this becomes 
\[
\left. \left( E+X_{f}\right) \right| _{r=1}=\left( \frac{1}{2(c-1)}+\frac{c+1}{4}\left. \der{f}{\phi }\right| _{r=1}\right) \left. \der{}{r}\right| _{r=1}+\frac{c+1}{4}\left. \der{f}{r}\right| _{r=1}\left. \der{}{\phi }\right| _{r=1}\]
 In order for this to vanish it is necessary, in particular, that \( \left. \der{f}{\phi }\right| _{r=1} \)
be a nonzero constant which is impossible since \( f \) is periodic in \( \phi . \) 
\end{proof}
We have now arrived at our final result:

\begin{thm}
The Poisson cohomology of a necklace Poisson structure \( \pi _{c} \) on \( S^{2} \)
is given as follows: 
\[
\begin{array}{ccccl}
H_{\pi _{c}}^{0}(S^{2}) & = & \R  & = & \textrm{span}\{1\}\\
H_{\pi _{c}}^{1}(S^{2}) & = & \R  & = & \textrm{span}\{\Delta _{\omega }\}\\
H_{\pi _{c}}^{2}(S^{2}) & = & \R ^{2} & = & \textrm{span}\{\pi _{c},\pi \}
\end{array}\]
 
\end{thm}
\begin{cor}
\( \pi _{c} \) does not admit infinitesimal rescaling.
\begin{cor}
The necklace structures \( \pi _{c} \) and \( \pi _{c'} \) for \( c\neq c' \)
are nontrivial deformations of each other.
\end{cor}
\end{cor}
\begin{proof}
\( \pi _{c'}-\pi _{c} \) is a nonzero multiple of \( \pi  \) but \( \pi  \)
is nontrivial in \( H_{\pi _{c}}^{2}(S^{2}) \).
\end{proof}

\appendix\label{appendix}

\chapter{Poisson manifolds and Poisson cohomology}

\begin{defn}
A \emph{Poisson manifold} is a manifold \( P \) together with an \( \R  \)-bilinear
skew-symmetric operation \( \{\cdot ,\cdot \} \) on \( \cinf (M) \), called
the \emph{Poisson bracket,} satisfying the following properties:
\begin{itemize}
\item The Leibniz rule: \( \forall f,g,h\in \cinf (M) \), 
\[
\{f,gh\}=\{f,g\}h+f\{g,h\}\]
 
\item The Jacobi identity: \( \forall f,g,h\in \cinf (M) \),
\[
\{\{f,g\},h\}+\{\{h,f\},g\}+\{\{g,h\},f\}=0\]
 
\end{itemize}
\end{defn}
Since \( \{\cdot ,\cdot \} \) is skew-symmetric and satisfies the Leibniz rule,
there exists a bivector field \( \pi \in \X ^{2}(M)=\Gamma (\bigwedge ^{2}TM) \)
such that 
\[
\{f,g\}=(df\wedge dg)(\pi )\]
 This bivector field is called the \emph{Poisson structure.} To express the
Jacobi identity in terms of \( \pi  \), recall that the \emph{Schouten bracket}
of multivector fields is defined as the unique extension \( \br  \) of the
commutator bracket of vector fields and the action of vector fields on functions
to \( \X ^{\star }(M)=\Gamma (\bigwedge ^{\star }TM) \) such that:

\begin{enumerate}
\item \( [X,Y]=-(-1)^{pq}[Y,X] \), for \( X\in \X ^{p+1}(M) \), \( Y\in \X ^{q+1}(M) \), 
\item \( [X,f]=X\cdot f \) for \( X\in \X (M) \), \( f\in \cinf (M) \), 
\item If \( X,Y\in \X (M) \), \( [X,Y] \) is the commutator bracket,
\item For \( X\in \X ^{p+1}(M) \), \( [X,\cdot ] \) is a derivation of degree \( p \)
of the exterior multiplication on \( \X ^{\star }(M) \). 
\end{enumerate}
The Schouten bracket satisfies the graded Jacobi identity 
\[
[X,[Y,Z]]=[[X,Y],Z]+(-1)^{pq}[Y,[X,Z]]\]
 for \( X\in \X ^{p+1}(M) \), \( Y\in \X ^{q+1}(M) \), \( Z\in \X ^{r+1}(M) \).
One then checks that the bracket on \( \cinf (M) \) given by a \( \pi \in \X ^{2}(M) \)
satisfies the Jacobi identity if and only if \( \pi  \) satisfies 
\begin{equation}
\label{Jacobi}
[\pi ,\pi ]=0
\end{equation}

A smooth map \( f:(P_{1},\pi _{1})\rightarrow (P_{2},\pi _{2}) \) is called
a \emph{Poisson map} if \( f_{*}\pi _{1}=\pi _{2}. \) The standard constructions
such as Poisson submanifolds and direct products are defined in an obvious manner.
In terms of the Poisson brackets, a submanifold \( N\subset P \) is Poisson
if and only if its vanishing ideal \( I_{N}\subset C^{\infty }(P) \) is a Poisson
ideal; if it is merely a Poisson subalgebra, \( N \) is said to be \emph{coisotropic. }

There are several important geometric objects associated with a Poisson bivector
field \( \pi . \) First, it gives rise to a bundle map \( \pi ^{\#}:T^{*}P\rightarrow TP \)
given by\\
\( <\alpha ,\pi ^{\#}\beta >=<\alpha \wedge \beta ,\pi > \) for any \( \alpha ,\beta \in T^{*}_{p}P. \)
To any function \( f\in C^{\infty }(P) \) one associates its \emph{Hamiltonian
vector field \( X_{f} \)} by 
\[
X_{f}\cdot g=<dg,X_{f}>=\{g,f\}=<dg,\pi ^{\#}df>=[\pi ,f]\cdot g\]
 The image of \( \pi ^{\#} \) is a (generally singular) integrable distribution
on \( P \). Since by definition every Hamiltonian vector field is tangent to
each integral submanifold, it follows easily that the integral submanifolds
are Poisson submanifolds of \( P \) called the \emph{symplectic leaves} of
\( P \) because the restriction of \( \pi  \) to each leaf is nonsingular,
hence symplectic. In general the leaves have different dimensions. 

A Poisson structure \( \pi  \) also gives rise to an operator 
\[
d_{\pi }:\fr{X}^{\star }\longrightarrow \fr{X}^{\star +1}\]
 on multivector fields given by \( d_{\pi }(X)=[\pi ,X] \). The graded Jacobi
identity for the Schouten bracket combined with (\ref{Jacobi}) implies that
\( d^{2}_{\pi }=0, \) making \( \fr{X}^{\star } \) into a complex. The cohomology
of this complex is called the \emph{Poisson cohomology} of \( (P,\pi ) \),
denoted by \( H_{\pi }^{\star }(P). \) The Poisson cohomology in low degrees
has a clear geometric interpretation: \( H_{\pi }^{0}(P) \) is the center of
the Poisson algebra \( C^{\infty }(P), \) consisting of \emph{Casimir functions,}
i.e those whose Hamiltonian vector fields are trivial; \( H_{\pi }^{1}(P) \)
consists of infinitesimal Poisson automorphisms of \( P \) (Poisson vector
fields) modulo inner automorphisms (Hamiltonian vector fields); \( H_{\pi }^{2}(P) \)
consists of nontrivial infinitesimal deformations of \( \pi  \) and \( H_{\pi }^{3}(P) \)
houses obstructions to extending an infinitesimal deformation to a full deformation
(see \cite{Vaisman:Book2}). 

The operator \( \pi ^{\#}:\fr{X}^{\star }\rightarrow \Omega ^{\star } \) intertwines
\( d_{\pi } \) and the deRham differential \( d \), hence induces a map \( \pi ^{\#}:H_{\pi }^{\star }(P)\rightarrow H_{dR}^{\star }(P) \)
which is an isomorphism if \( \pi  \) is symplectic. In general, however, Poisson
cohomology has been notoriously difficult to compute, and there have been but
a handful of successful computations (\textsc{\cite{Gin1},\cite{Xu},\cite{GinLu}). }

For every Poisson manifold there are two Poisson cohomology classes that are
special. The first one is the \emph{modular class}, introduced by Weinstein
\cite{We3}. Given a Poisson manifold \( (P,\pi ) \) with a volume form \( \omega  \),
Weinstein defines an operator \( \Delta _{\omega } \) on \( C^{\infty }(P) \)
that associates to every function the divergence of its Hamiltonian vector field
with respect to \( \omega  \). It turns out that \( \Delta _{\omega } \) is
in fact a vector field called the \emph{modular vector field} of \( \pi  \)
with respect to \( \omega  \). Moreover, \( \Delta _{\omega } \) preserves
\( \pi  \); if \( \omega  \) is replaced by another volume form, \( \Delta _{\omega } \)
is changed by a Hamiltonian vector field. Thus, the class of \( \Delta _{\omega } \)
in \( H_{\pi }^{1}(P) \) is independent of \( \omega  \); it is called the
\emph{modular class} of \( (P,\pi ) \) and measures the obstruction to the
existence of a measure on \( P \) invariant under all Poisson automorphisms.
It is zero for symplectic manifolds due to the existence of the Liouville measure.

The second special class is the class of \( \pi  \) itself in \( H_{\pi }^{2}(P). \)
It is the obstruction to (infinitesimal) rescaling of \( \pi  \). If it vanishes,
\( \pi  \) is called \emph{exact} and there exists a vector field \( X \)
such that \( L_{X}\pi =\pi . \) It is called a \emph{Liouville vector field.}
If \( \pi  \) is symplectic, the class of \( \pi  \) corresponds under \( \pi ^{\#} \)
to the class of the symplectic form \( \omega  \) in the deRham cohomology,
hence it admits infinitesimal rescaling if and only if \( \omega  \) is exact.

\chapter{Supermanifolds}

We recall here the rudiments of the theory of supermanifolds that should suffice
for the understanding of the material in the main text. For a more thorough
introduction the interested reader should consult \cite{Manin}, \cite{Vor2}
or \cite{Leit}.

\section{Algebra}

The basic setting of supermathematics is the category Super of \( \Z _{2} \)-graded
vector spaces \( V=V_{0}\oplus V_{1} \). The \( \Z _{2} \)-grading is called
\emph{parity.} Elements of \( V_{0} \) are called \emph{even,} while elements
of \( V_{1} \) are \emph{odd}; \emph{}the parity of an element is denoted by
a tilde over it. 

If \( V_{0} \) and \( V_{1} \) are finite-dimensional, the dimension of \( V \)
takes values in \( \Z [\Pi ]/(\Pi ^{2}-1) \), the group ring of \( Z_{2} \),
and is denoted by \( \dim V=(\dim V_{0}|\dim V_{1}) \). The \emph{parity reversion
functor} \( \Pi  \) is defined by 
\[
(\Pi V)_{0}=V_{1},\; (\Pi V)_{1}=V_{0}\]
 All he usual universal constructions, such as the direct sum, tensor product,
duality and Hom carry over to the Super category, with a natural assignment
of parity. The notion of an associative algebra in the Super category is the
usual one, except that the multiplication must respect the parity: \( V_{i}V_{j}\subset V_{i+j} \).
The commutator in an associative superalgebra is the supercommutator 
\[
[a,b]=ab-(-1)^{\tilde{a}\tilde{b}}ba\]
 One calls the algebra commutative if this bracket is identically zero.\footnote{
We try not to abuse the prefix ``super'', omitting it whenever it is clear
from the context that we are working in the Super category.
} In general, the sign convention - introducing \( (-1)^{ij} \) whenever two
symbols of parities \( i \) and \( j \) are interchanged - should be used
as a matter of principle in the Super category. Thus, the notions of ``symmetric''
and ``skew-symmetric'' must be modified appropriately. An endomorphism \( D \)
of \( V \) is a \emph{derivation} if 
\[
D(ab)=(Da)b+(-1)^{\tilde{D}\tilde{a}}a(Db)\]
 Derivations of any other kind of a bilinear operation are defined analogously.

One defines a Lie superalgebra (of parity \( \epsilon  \)) to be a vector space
\( V \) with a bilinear skew-symmetric operation \( \br  \) of parity \( \epsilon  \)
(i.e., \( [V_{i},V_{j}]\subset V_{i+j+\epsilon } \)) satisfying the Jacobi
identity: 
\[
[a,[b,c]]=[[a,b],c]+(-1)^{(\tilde{a}+\epsilon )(\tilde{b}+\epsilon )}[b,[a,c]]\]
 The parity reversion functor \( \Pi  \) interchanges the notions of even and
odd Lie superalgebras, so one can always reduce to the even case. If \( V \)
also has a commutative multiplication with respect to which \( ad_{a} \) is
a derivation (of parity \( \tilde{a}+\epsilon  \), of course), it becomes a
Poisson superalgebra, even or odd. An odd Poisson algebra is also called a \emph{Gerstenhaber
algebra.} Various Schouten-like bracket structures encountered in the main text
are Gerstenhaber algebras. It is no longer possible to reduce Gerstenhaber algebras
to even Poisson algebras by parity reversion.

\section{Affine superspaces and superdomains}

A function of odd variables \( \xi ^{1},\ldots ,\xi ^{m} \) is an element of
the free commutative algebra generated by these variables, i.e. the Grassman
algebra \( \bigwedge (\R ^{m})^{*} \): 
\[
f(\xi ^{1},\ldots ,\xi ^{m})=f_{0}+\xi ^{\mu }f_{\mu }+\half \xi ^{\mu _{1}}\xi ^{\mu _{2}}f_{\mu _{1}\mu _{2}}+\cdots +\frac{1}{n!}\xi ^{\mu _{1}}\cdots \xi ^{\mu _{m}}f_{\mu _{1}\ldots \mu _{n}}\]
 where 
\[
f_{\mu _{1}\ldots \mu _{k}}=(-1)^{\sigma }f_{\mu _{\sigma (1)}\ldots \mu _{\sigma (k)}}\]
 for any permutation \( \sigma \in S_{k} \). The variables \( \xi ^{\mu } \)
are to be interpreted as coordinates on the purely odd affine superspace \( \R ^{0|m} \),
which can be thought of as the result of applying the parity reversal functor
\( \Pi  \) to \( R^{m} \). More invariantly, one says that the algebra of
functions on the superspace \( \Pi V \) is the Grassman algebra \( \bigwedge V^{*} \). 

If the coefficients in the above expression are themselves smooth functions
of even variables \( x^{1},\ldots ,x^{n} \) defined on \( \R ^{n} \) or an
open subset \( U_{0}\subset \R ^{n} \), one says that we are given a function
on an \emph{affine superspace} \( \R ^{n|m} \) or a \emph{superdomain \( U^{n|m}\subset \R ^{n|m} \).}
These functions form a supercommutative algebra which is just the tensor product
\( \cinf (U_{0})\otimes \bigwedge (\R ^{m})^{*} \). The domain \( U_{0} \)
is called the \emph{support} of \( U^{n|m} \) and uniquely determines it since
we cannot ``bound'' the odd variables. It is often convenient not to separate
the even and odd variables explicitly but denote them by a collective symbol
\( \{x^{A}\} \) and assign parity to each index \( A \).

If \( (t^{1},\ldots ,t^{p},\tau ^{1},\ldots ,\tau ^{q}) \) is another set of
variables, we can define a substitution 
\begin{equation}
\label{eqn:substit}
\begin{array}{rcl}
x^{a} & = & x^{a}(t,\tau )=x^{a}_{0}(t)+\half \tau ^{\alpha _{1}}\tau ^{\alpha _{2}}x^{a}_{\alpha _{1}\alpha _{2}}(t)+\cdots \\
\xi ^{\mu } & = & \xi ^{\mu }(t,\tau )=\tau ^{\alpha }\xi ^{\mu }_{\alpha }(t)+\sixth \tau ^{\alpha _{1}}\tau ^{\alpha _{2}}\tau ^{\alpha _{3}}\xi ^{\mu }_{\alpha _{1}\alpha _{2}\alpha _{3}}(t)+\cdots 
\end{array}
\end{equation}
 where \( x_{0}^{\alpha }(t),\ldots  \) are smooth functions defined on a domain
\( V_{0}\subset \R ^{p} \), and plug these expression into any \( f=f(x,\xi ) \).
This is not a problem since the \( \tau  \)'s are nilpotent. 

\begin{example*}
\( \sin (t+\tau ^{1}\tau ^{2})=\sin t+\tau ^{1}\tau ^{2}\cos t \) 
\end{example*}
We interpret such substitutions as smooth maps \( V^{p|q}\rightarrow U^{m|n} \).
For example, the inclusion of the support \( U_{0}\hookrightarrow U^{n|m} \)
is defined by setting all the odd variables to zero. More formally, the algebras
\( \cinf (U_{0})\otimes \bigwedge (\R ^{m})^{*} \) form a category with morphisms
defined by the substitutions above, and we simply define the category of superdomains
to be the opposite category. What makes the whole theory of supermanifolds nontrivial
is the possibility of mixing the even and odd variables by allowing nonlinear
terms in (\ref{eqn:substit}). If \( V^{p|q}=U^{n|m} \) and the substitutions
(\ref{eqn:substit}) are invertible, one can think of \( (t,\tau ) \) as giving
a new coordinate system on \( U^{n|m} \). Thus the domain \( U^{n|m} \) is
viewed as the intrinsic geometric object for which we can choose a coordinate
representation at will. This is completely analogous to the ordinary, purely
even case, and leads to the concept of a \emph{supermanifold.}

The derivatives \( \partial /\partial x^{A} \) are defined in the usual manner,
as linear endomorphisms such that one has the Leibniz rule 
\[
\der{(fg)}{x^{A}}=\der{f}{x^{A}}g+(-1)^{\tilde{A}\tilde{f}}f\der{g}{x^{A}}\]
 and 
\[
\der{x^{B}}{x^{A}}=\delta ^{B}_{A}\]
 Then the usual equality of mixed partials holds with appropriate signs,\footnote{
These are left derivatives. There are also right derivatives, satisfying the
Leibniz rule when applied on the right of the argument. The difference between
the left derivative of a function and a right one is only a sign. Right derivatives
are not used in this work.
} and so does the inverse function theorem: the substitution (\ref{eqn:substit})
is locally invertible if and only if its Jacobian matrix is. 

The only difficulty is the absence of a good notion of points: the only numerical
value an odd variable can be assigned is zero, so the values that can be assigned
to functions on a superdomain can only determine its support. For this reason
we try to formulate all our statements in terms of the algebra of functions;
whenever we mention points, we mean ``running points'', i.e., all objects
considered will be allowed to depend, explicitly or implicitly, on any number
of even and odd parameters.

\section{Supermanifolds}

Smooth supermanifolds are glued together out of domains \( U^{n|m} \) in the
same way in which ordinary manifolds are glued out of coordinate domains. To
make this rigorous, one has to use sheaf theory. One considers a locally ringed
space \( M=(M_{0},\OO _{M}) \) where \( M_{0} \) is a topological space and
\( \OO _{M} \) is a sheaf of supercommutative algebras on \( M_{0} \) whose
stalk \( \OO _{x} \) over each point \( x\in M_{0} \) is local. A superdomain
\( U^{n|m}=(U_{0},\cinf (U_{0})\otimes \bigwedge V^{*}) \) is such a space.
A \emph{chart} on \( M \) is, by definition, an isomorphism of locally ringed
spaces \( \phi :V=(V_{0},\OO _{M}|_{V_{0}})\rightarrow U^{n|m} \), where \( V_{0}\subset M_{0} \)
is an open subset. One says that \( M \) is a \emph{supermanifold} if it can
be covered by a countable system of charts, called an \emph{atlas}. Atlases
form a directed set under refinement, and every atlas is contained in a maximal
one. Factoring \( \OO _{M}|_{V_{0}} \) by its nilradical one gets an atlas
on \( M_{0} \) making it a smooth manifold called, naturally, the \emph{support}
of \( M \). 

The charts \( \phi  \) give a system of local coordinates on \( M \), allowing
us to describe supermanifolds and their morphisms in coordinates while making
almost no use of sheaves. Due to the absence of a good notion of points mentioned
above, this is the best course to follow. In this way, all the standard constructions
- products, co-products, fibered products, submanifolds, vector bundles - carry
over to supermanifolds. Thus, a closed submanifold (or, more generally, a singular
subvariety) is locally given by a system of equations: for example, the support
\( M_{0} \) is given by setting all the nilpotents to zero. 

Given a supermanifold \( M \), one defines its tangent bundle by giving, for
each chart \( (V,\{x^{A}\}) \) on \( M \) a chart \( (TV,\{x^{A},\dot{x}^{A}\}) \)
on \( TM \) such that under a change of coordinates \( x=x(x') \) the velocities
transform in the usual way: 
\[
\dot{x}^{A}=\dot{x}^{A'}\der{x^{A}}{x^{A'}}(x')\]
 Similarly, one defines the cotangent bundle \( T^{*}M \).

A basic class of examples of supermanifolds are provided by supermanifolds of
the form \( \Pi A \) where \( A \) is a vector bundle over an ordinary manifold
\( M_{0} \). The structure sheaf of \( \Pi A \) is the sheaf of smooth sections
of \( \bigwedge A^{*} \). Any atlas on \( M_{0} \) gives rise to an atlas
on \( \Pi A \) with the coordinate transformations inherited from the vector
bundle structure on \( A \). This atlas is characterized by the property that
the even coordinates transform independently of the odd ones, while the odd
ones transform linearly. Such atlases are called \emph{simple.} The fundamental
classification theorem of smooth real supermanifolds (proved independently by
Berezin \cite{Ber}, Bachelor \cite{Bach} and Gawedzki \cite{Gaw}) asserts
that any supermanifold \( M \) admits a simple atlas, i.e. globally isomorphic
to one of the form \( \Pi A \). The bundle \( A \) in question is the normal
bundle to the support \( M_{0} \) of \( M \), with the parity in the fibres
reversed. One must emphasize, however, that this isomorphism is strictly non-canonical.
For example, the cotangent bundle \( T^{*}\Pi A \), the central object in the
main text, does not possess a canonical simple atlas (Remark \ref{rem:Gradings}).
The theorem is also false for complex analytic supermanifolds \cite{Manin}.

\section{Vector fields and differential forms}

A vector field on a supermanifold is simply a derivation of its algebra of functions.
Vector fields on a supermanifold \( M=(M_{0},\OO _{M}) \) form a sheaf of left
\( \OO _{M} \)-modules. The Jacobi-Lie bracket of vector fields is just the
commutator of the corresponding derivations: 
\[
[X,Y]f=X(Yf)-(-1)^{\tilde{X}\tilde{Y}}Y(Xf)\]
 In local coordinates, a vector field is given by an expression of the form
\[
X=X^{A}(x)\der{}{x^{A}}=X^{a}(x,\xi )\der{}{x^{a}}+X^{\mu }(x,\xi )\der{}{\xi ^{\mu }}\]
 Unlike the even case, if a vector field \( X \) is odd, the condition 
\[
[X,X]=2X^{2}=0\]
 is nontrivial. When it is satisfied, the vector field \( X \) is called \emph{homological}
because it endows the algebra of functions with the structure of a differential
complex. Vector fields on \( M \) also correspond naturally to fiberwise linear
functions on \( T^{*}M \).

Differential forms on supermanifolds, as objects suitable for integration, are
highly nontrivial \cite{Vor2}; however, for our purposes it suffices to consider
the simplest class of differential forms on \( M \) which are polynomial functions
on \( \Pi TM \). A local coordinate chart \( \{x^{A}\} \) on \( M \) induces
a chart \( \{x^{A},\xi ^{A}\} \) on \( \Pi TM \) where \( \xi ^{A}=dx^{A} \)
have parity \( \tilde{A}+1 \) and transform just as the notation suggests.
Thus a differential form \( M \) is locally a function \( \omega =\omega (x,dx) \),
polynomial in \( dx \). If this restriction is removed, we get the so-called
\emph{Bernstein-Leites pseudoforms.}

\begin{example*}
\( \omega =e^{-(d\xi )^{2}} \) is a pseudoform on \( \R ^{0|1} \).
\end{example*}
The degree of \( \omega  \) as a polynomial in \( dx \) in general differs
from its parity as a function on \( \Pi TM \). On \( \Pi TM \) there is a
canonical homological vector field, the de Rham differential
\[
d=\xi ^{A}\der{}{x^{A}}\]
 Any vector field \( X \) on \( M \) induces a vector field on \( \Pi TM \),
the interior derivative 
\[
i_{X}=(-1)^{\tilde{X}}X^{A}\der{}{\xi ^{A}}\]
 if \( X=X^{A}\der{}{x^{A}} \), of parity \( \tilde{X}+1 \) and the Lie derivative
\( L_{X}=[d,i_{X}] \), of parity \( \tilde{X} \).

\section{Symplectic and Poisson supermanifolds}

A symplectic structure on a supermanifold \( M \) is a two-form \( \omega  \)
(i.e. a quadratic function on \( \Pi TM \)), 
\[
\omega =\half dx^{A}dx^{B}\omega _{AB}(x)=\half dx^{a}dx^{b}\omega _{ab}+dx^{a}d\xi ^{\mu }\omega _{a\mu }+\half d\xi ^{\mu }d\xi ^{\nu }\omega _{\mu \nu }\]
 which is closed (\( d\omega =0 \)) and nondegenerate (the matrix \( \omega _{AB} \)
is invertible). One distinguishes even and odd symplectic supermanifolds, depending
on the parity of \( \omega  \). In the even case, the nondegeneracy of \( \omega  \)
is equivalent to the invertibility of the real matrices \( \omega _{ab} \)
and \( \omega _{\mu \nu } \), after setting the nilpotents to zero. Note that
\( \omega _{ab} \) is skew-symmetric while \( \omega _{\mu \nu } \) is symmetric;
the dimension of \( M \) in this case has to be \( 2n|m \), and the signature
of \( \omega _{\mu \nu } \) is an invariant. In the odd case, the nondegeneracy
is equivalent to the invertibility of the real matrix \( \omega _{a\mu } \),
hence the dimension of \( M \) must be \( n|n \). The Darboux theorem holds
for supermanifolds and asserts that \( \omega  \) locally has the standard
form 
\[
\omega =dp_{a}dq^{a}+\half \sum _{\mu }\pm (d\gamma ^{\mu })^{2}\]
 for \( \tilde{\omega }=0 \), and 
\[
\omega =d\theta _{a}dx^{a}\]
 for \( \tilde{\omega }=1 \). Even and odd symplectic supermanifolds have very
different properties. 

For differential forms there is a natural notion of pullback, in particular,
restriction to submanifolds. As usual, one calls a submanifold \( L \) of \( M \)
\emph{Lagrangian} if the restriction of \( \omega  \) to \( L \) is identically
zero and \( L \) is of the maximal dimension where it is possible. Lagrangian
submanifolds of an even symplectic supermanifold \( M^{2n|m} \) have dimension
\( n|[m/2] \) (they may not even exist if the signature of \( \omega _{\mu \nu } \)
is nonzero), while those of an odd one \( M^{n|n} \) have dimension \( k|n-k \). 

Given a function \( f \) on \( M \), its hamiltonian vector field is defined
by the formula 
\[
i_{X_{f}}\omega =-df\]
 and for a pair of functions \( f,g \) their Poisson bracket is defined by
\[
\{f,g\}=X_{f}g\]
 If \( \tilde{\omega }=0 \), \( (\cinf (M),\{\cdot ,\cdot \}) \) becomes an
even Poisson algebra; if \( \tilde{\omega }=1 \), a Gerstenhaber algebra. 

A typical example of an even symplectic supermanifold is \( T^{*}Q \), where
\( Q \) is a supermanifold, with the standard symplectic structure 
\[
\omega =dx_{A}^{*}dx^{A}\]
 (the momenta \( x^{*}_{A} \) have the same parity as \( x^{A} \)), whereas
a typical odd symplectic supermanifold is \( \Pi T^{*}Q \) with 
\[
\omega =d\theta _{A}dx^{A}\]
 (\( \tilde{\theta }_{A}=\tilde{x}^{A}+1 \)). Functions on \( \Pi T^{*}Q \)
coincide with multivector fields on \( Q \), and the canonical odd Poisson
bracket is nothing but the Schouten bracket of multivector fields. 

One can similarly introduce even or odd Poisson supermanifolds that are not
necessarily symplectic. For example, any Lie algebroid \( A \) gives rise to
an odd Poisson structure on \( \Pi A^{*} \) via the generalized Schouten bracket,
as explained in the main text. 

\begin{rem*}
Finally, we remark that we have left out the part of the supermanifold theory
that is perhaps the most interesting and different from the ordinary manifold
case - the integration theory. This is only because it is not used anywhere
in the main body of this work. The interested reader is strongly advised to
look in the treatise \cite{Vor2} for a thorough exposition.
\end{rem*}

\end{document}